\documentclass[12pt]{amsbook}

\usepackage{euscript}
\usepackage{amsfonts}
\usepackage{amsmath,amscd, amsthm, amssymb}
\usepackage[english, frenchb]{babel}
\usepackage[latin1]{inputenc}
\usepackage{graphicx}
\usepackage[arrow, matrix, curve]{xy}
\usepackage{dsfont}\let\mathbb\mathds
\usepackage{rotating}

\newtheorem{thm}{Théorème}[chapter]
\newtheorem{prop}[thm]{Proposition}
\newtheorem{cor}[thm]{Corollaire}
\newtheorem{lem}[thm]{Lemme}

\newtheorem{rem}[thm]{Remarque}
\newtheorem{Def}[thm]{Définition}

\def\CC{{\mathbb{C}}}

\def\RR{{\mathbb{R}}}
\def\ZZ{{\mathbb{Z}}}
\def\NN{{\mathbb{N}}}
\def\PP{{\mathbb{P}}}

\def\cartesien{%
    \ar@{-}[]+R+<6pt,-1pt>;[]+RD+<6pt,-6pt>%
    \ar@{-}[]+D+<1pt,-6pt>;[]+RD+<6pt,-6pt>%
  }
\newcommand{\immouv}[1][r]
   {\ar@{}[#1] |*[o][F]{\hbox{%
         \vrule width 1.5mm height 0pt depth 0pt%
         \vrule width 0pt height .75mm depth .75mm%
         }}
     \ar@{^{(}->}[#1]}
\newcommand{\demde}[1]{\begin{proof} de #1}
\newcommand{\dem}{\begin{proof}}
\newcommand{\cqfd}{\end{proof}}
\let\cat\mathfrak 
  \newcommand{\UN}[4][r]{%
    \ar@/^1pc/[#1]^{#2}_*=<0.3pt>{}="HAUT"
    \ar@/_1pc/[#1]_{#3}^*=<0.3pt>{}="BAS"
    \ar @{=>} "HAUT";"BAS" ^{#4}
  }
\newcommand{\DEUX}[6][r]{
    \ar@/^2pc/[#1]^{#2}_*=<0.3pt>{}="HAUT"
    \ar@{}    [#1]     ^*=<0.3pt>{}="MILIEUHAUT"
                       _*=<0.3pt>{}="MILIEUBAS"
    \ar[#1]_(0.3){#3}                  
    \ar@/_2pc/[#1]_{#4}^*=<0.3pt>{}="BAS"
    \ar @{=>} "HAUT";"MILIEUHAUT" ^{#5}
    \ar @{=>} "MILIEUBAS";"BAS" ^{#6}
  }   
 \newcommand{\eq}[1][r]
   {\ar@<-3pt>@{-}[#1]
    \ar@<-1pt>@{}[#1]|<{}="gauche"
    \ar@<+0pt>@{}[#1]|-{}="milieu"
    \ar@<+1pt>@{}[#1]|>{}="droite"
    \ar@/^2pt/@{-}"gauche";"milieu"
    \ar@/_2pt/@{-}"milieu";"droite"}
 \newcommand{\incl}[1][r]
  {\ar@<-0.2pc>@{^(-}[#1] \ar@<+0.2pc>@{-}[#1]}

\newcommand{\Ac}{\mathcal{A}}
\newcommand{\Bc}{\mathcal{B}}
\newcommand{\Cc}{\mathcal{C}}
\newcommand{\Dc}{\mathcal{D}}
\newcommand{\Fc}{\mathcal{F}}
\newcommand{\Gc}{\mathcal{G}}
\newcommand{\Hc}{\mathcal{H}}
\newcommand{\Ic}{\mathcal{I}}
\newcommand{\Jc}{\mathcal{J}}
\newcommand{\Kc}{\mathcal{K}}
\newcommand{\Lc}{\mathcal{L}}

\newcommand{\Mm}{\mathcal{M}}

\newcommand{\Pc}{\mathcal{P}}
\newcommand{\Ss}{\mathcal{S}}
\newcommand{\Tc}{\mathcal{T}}

\newcommand{\Sc}{\mathcal{S}}

\newcommand{\CCC}{\mathfrak{C}}
\newcommand{\PPP}{\mathfrak{P}}
\newcommand{\LLL}{\mathfrak{L}}
\newcommand{\GGG}{\mathfrak{G}}
\newcommand{\III}{\mathfrak{I}}
\newcommand{\SSS}{\mathfrak{S}}
\newcommand{\SSSt}{\mathfrak{St}}
\newcommand{\Ff}{\EuScript{F}}
\newcommand{\Ii}{\EuScript{I}}
\newcommand{\Ll}{\EuScript{L}}

\newcommand{\AAA}{\mathbf{A}}
\newcommand{\BBB}{\mathbf{B}}
\newcommand{\ccc}{\mathbf{C}}

\begin{document}

\thispagestyle{empty}

\begin{center}

{UNIVERSIT\'E DE NICE--SOPHIA ANTIPOLIS   --  UFR Sciences}\\

\vspace{0.4cm}

\'Ecole Doctorale Sciences Fondamentales et Appliqu\'ees\\



\vspace{1.cm}

{\Large\bf TH\`ESE} \\
\vspace{0.3cm}
pour obtenir le titre de \\
\vspace{0.1cm}
{\Large\bf Docteur en Sciences}\\
\vspace{0.1cm}

Sp\'ecialit\'e : {\sc Math{\'e}matiques}\\

\vspace{0.8cm}

{pr\'esent\'ee et soutenue par}\\
{\large\bf Delphine DUPONT}\\

\vspace{1.4cm}

{\LARGE{\bf{Exemples de classification du champ des faisceaux pervers}}}\\

\vspace{1.4cm}
{Th\`ese dirig\'ee par \bf{Philippe MAISONOBE}}\\
\vspace{0.2cm}
{soutenue le 4 décembre 2008}

\vspace{0.8cm}

Membres du jury : \\
\vspace{0.4cm}

\begin{tabular}{llll}
A. D'Agnolo & & Professeur à l'Università degli Studi di Padova  & Rapporteur\\
A. Dimca & & Professeur \`a l'Universit\'e de Nice & Examinateur \\
M. Granger & & Professeur \`a l'Universit\'e d'Angers & Rapporteur \\
F. Loeser & & Professeur à l'École normale supérieure & Rapporteur \\
Ph. Maisonobe & & Professeur \`a l'Universit\'e de Nice & Directeur \\
I. Waschkies & & Maître de conférences à l'Université de Nice & Examinateur\\
\end{tabular}

\vspace{1.0cm}

{Laboratoire J.-A. Dieudonn\'e}\\
\vspace{0.1cm}{Universit\'e de Nice}\\
\vspace{0.1cm}
Parc Valrose, 06108 NICE Cedex 2



\end{center}

\newpage
\frontmatter
\newpage
~\\\\\\\\\\\\\\\\\\\\\\\\\\\\\\\\\\\\\\\

\hspace{10cm}
\emph{\`A ma mère.}
\newpage
~\\
\newpage
\begin{center}\begin{bf}
Remerciements 
\end{bf}
\end{center}
Au-delà de l'exercice, qu'il me paraît difficile d'exprimer justement toute ma gratitude. En relisant ces mots je les trouve bien en de\c cà de ce que j'éprouve.

Tout d'abord je tiens à remercier mon directeur, Philippe Maisonobe, sans qui je n'aurais pas commencé cette thèse. Sans lui je n'aurais jamais découvert  les affres et les joies de cette théorie pas si perverse. Je voudrais aussi le remercier du temps qu'il a su me consacrer malgré sa lourde charge  de directeur de laboratoire. Je le remercie aussi de ses nombreuses relectures de ma prose.

Je souhaite tout particulièrement remercier Ingo Waschkies, pour sa disponibilité même à des milliers de kilomètres, pour toutes ces longues discussions si enrichissantes, pour ces quelques ``tu as raison'' lâchés avec conviction après m'avoir obligé  pour le convaincre, à pousser mon raisonnement bien plus loin que je ne l'avais fait, pour toutes ces idées dont il m'a fait si généreusement part. Enfin et ce n'est pas la moindre des choses pour ces re-re-re-re-re-relectures de ce texte. Je tiens encore à le remercier de son indéfectible franchise. Si ces remerciements, la soutenance ou pire le pot ne lui conviennent pas soyez assurés que j'en serais informée. 

Je voudrais bien sûr remercier les trois rapporteurs de cette thèse, Andrea D'Agnolo, Michel Granger et Fran\c cois Loeser, d'avoir pris le temps de se pencher sur ce manuscrit, ainsi qu'Alexandru Dimca qui a accepté de faire partie du jury.

Je tiens aussi à remercier ici Georges Comte pour ses encouragements parfois sévères mais toujours bienveillants. Georges a été un de mes enseignants de licence et je me rappelle notamment d'une phrase lancée durant un td qui m'a, bien plus que n'importe quel discours sur les mathématiques, éclairée et encouragée à continuer dans ce sens. Je venais de comprendre qu'une chose que je connaissais bien, la dérivation, pouvait être vue autrement, comme une application linéaire, et que nous pouvions en déduire des informations su\-pplé\-men\-taires. Il me rétorqua alors, avec un enthousiasme mêlé de la pointe de froideur qui le caractérise : ``Mais c'est \c ca la puissance du formalisme''.

Pendant ces quatre années passées au laboratoire Dieudonné j'ai rencontré nombre de personnes qui sont aujourd'hui devenues mes amis. Quand les temps étaient difficiles, quand les mathématiques se faisaient indéchiffrables, j'avais au moins le plaisir de les y retrouver.  Je pense bien sûr par ordre chronologique à Guillaume, Ti Ha, Maëlle, Marcello, Xavier, Pierre, Patrick, Marc, Fabien, Michel, Hugues et tous les autres.
Je n'oublie, bien sûr, pas Marie, le partage de nos expériences douloureuses a été un véritable réconfort. 

Il me reste encore tant de gens à remercier. Toutes ces personnes rencontrées dans les conférences et qui ont montré un intérêt pour mon travail, tous ces professeurs qui m'ont encouragé et tellement appris, tous ces amis rencontrés durant mes études,  ma famille.

Enfin je tiens à remercier Violaine et Feres de ce qu'ils sont et de ce qu'ils savent être pour moi.

\newpage
~\\

\tableofcontents

\mainmatter
\setcounter{secnumdepth}{3}
\newpage

\chapter*{Introduction}
La catégorie abélienne des faisceaux pervers $\Pc erv_X$ sur un espace topologique $X$ a été introduite par A. Beilinson, J. Bernstien, P. Deligne et O. Gabber dans ~\cite{BBD}. Elle joue un rôle important dans de nombreuses branches des mathématiques, notamment en géométrie algébrique et en théorie des représentations. L'un des intérêts majeur est dû au théorème de Riemann-Hilbert qui établit, dans le cas où $X$ est une variété analytique complexe, l'équivalence entre $\Pc erv_{X}$ et la catégorie des $\Dc$-modules holonômes réguliers. La catégorie des faisceaux pervers est définie  comme sous-catégorie pleine de la catégorie dérivée des faisceaux à cohomologie constructible. Cette définition nécessite donc un langage algébrique très avancé et plusieurs méthodes ont été développées pour en donner, dans le cas où la stratification est fixée, une description élémentaire.\\

Plusieurs descriptions par des catégories de représentations de carquois se font sur un espace topologique particulier muni d'une stratification fixée. Rappelons qu'un carquois est un graphe orienté et qu'une représentation est un foncteur de ce graphe dans la catégorie des espaces vectoriels.\\
Un premier exemple de cette approche est donné par A. Galligo, M. Granger et Ph. Maisonobe qui, dans \cite{GGM}, démontrent l'équivalence de la catégorie des faisceaux pervers sur $\CC^n$ relativement à un croisement normal avec une sous-catégorie pleine, notée $\Cc_{n}$,  de la catégorie de représentations du carquois associé à un hypercube. En utilisant une méthode similaire Ph. Maisonobe démontre, dans \cite{Maiso2},  le même type de résultat avec $X=\CC^2$ stratifié par une courbe plane. On peut aussi citer T. Braden et M. Grinberg qui, en utilisant la théorie des faisceaux pervers microlocaux, décrivent, dans \cite{BG}, la catégorie $\Pc erv_X$ où $X$ est un espace de matrices stratifié par le rang. Citons S. Khoroshkin et A. Varchenko, dans \cite{KV}, ils considèrent $X=\CC^n$ muni de la stratification $\Sigma$ donnée par un arrangement d'hyperplans. Ils définissent un foncteur pleinement fidèle de la catégorie des repésentations du carquois associé à $\Sigma$ dans la catégorie des $\Dc$-modules holonômes réguliers. Un autre exemple de ce type est démontré dans \cite{Nar} par Gudiel Rodr\`iguez et L. Narv\`aez Macarro.\\
D'autres méthodes, plus générales mais du même coup moins ex\-pli\-cites, consistent à décrire la catégorie $\Pc erv_X$ dans le cas où $X$ est un espace de Thom-Mather. Citons par exemple R. MacPherson et K. Vilonen, dans \cite{McV}. C'est d'ailleurs en itérant ce processus que H. Larrouy démontre dans sa thèse \cite{Lar} l'équivalence entre $\Pc erv_{X}$ et une catégorie de représentations de carquois, pour $X$ une variété  torique affine. \\

Les résultats que nous venons de citer ne tiennent pas compte du caractère local des faisceaux pervers. Rappelons, en effet, que les faisceaux pervers sur un espace stratifié forment un champ. Autrement dit connaître un faisceau pervers sur un ouvert $U$ revient à le connaître sur un recouvrement ouvert. On peut donc s'attendre à ce qu'une caractérisation locale de la catégorie des faisceaux pervers par des catégories explicites de représentations de carquois puisse se recoller sur $X$ en une description simple de la catégorie globale des faisceaux pervers. Ceci permettra, par exemple, de décrire les faisceaux pervers sur un espace qui est localement un des espaces cités.

 Récemment, D. Treumann a  commencé à travailler dans cette direction. Il a donné dans \cite{Tr1} et \cite{Tr2} une description  du champ des faisceaux pervers en s'appuyant sur le fait que le champ $\PPP_X$ est un champ constructible (dont la restriction à chaque strate est localement constante). Il démontre que sous certaines conditions topologiques, la catégorie globale $\Pc erv_{X}$ est équivalente à une catégorie de modules sur une algèbre de type fini. Il y arrive en recollant des descriptions locales qui s'inspirent de MacPherson et Vilonen qui ne sont pas des catégories de représentations de carquois explicites. 

Le but de cette thèse est de montrer comment on peut recoller des descriptions élémentaires et explicites de faisceaux pervers en utilisant la théorie des champs. Pour cela on considère des espaces qui sont localement isomorphe à $\CC^n$ stratifié par le croisement normal. C'est le cas des variétés toriques lisses stratifiées par l'action du tore et de $\CC^2$ stratifié par un arrangement générique de droites. Notre modèle local est donc  l'équivalence de catégories démontrées par  A. Galligo, M. Granger et Ph. Maisonobe dans \cite{GGM}. Il s'agit donc de généraliser cette équivalence de catégorie en une équivalence de champs puis de recoller ces équivalences.\\

Le premier chapitre est une présentation des objets avec lesquels nous travaillons. Nous donnons ainsi la définition et quelques propriétés essentielles des catégories de représentations de carquois. Pour décrire le champ des faisceaux pervers nous utilisons le langage des $2$-catégories, c'est pourquoi nous y consacrons un paragraphe. Enfin nous introduisons la $2$-catégorie des champs. Nous l'introduisons comme une généralisation des faisceaux en catégories.  \\

Considérons $X$ un espace topologique muni d'une stratification $\Sigma$. Un faisceau pervers est un complexe de faisceau à cohomologie constructible, c'est à dire ces faisceaux de cohomologie sont localement constant sur chaque strate. Une approche naïve pour coder un faisceau pervers consisterai à considérer ses restrictions à chacune des  strates et des données de recollement à préciser. Effectivement, considérons un faisceau constructible, on montrera dans le chapitre 2 qu'un faisceau sur $X$ est défini de manière unique par ses restrictions à chaque strate et pour chaque couple de strate un morphisme de recollement. Donc un faisceau constructible est codé par une représentation de carquois dont la combinatoire est donnée par la stratification. La situation pour les faisceaux pervers est plus compliqué car on travaille avec des complexes dans la catégorie dérivée. Par contre il s'avère que le champ des faisceaux pervers est un champ constructible pour $\Sigma$ ce qui nous permet de raffiner cette approche naïve en considérant non plus un seul faisceau pervers mais la catégorie $\Pc erv_{X}$. Pour motiver notre classification des champs strictement constructibles dans le chapitre 3, nous nous intéressons tout d'abord dans le deuxième chapitre aux recollements de faisceaux et de faisceaux pervers sur un espace stratifié. \\
La première partie est la démonstration de l'équivalence entre la catégorie des faisceaux sur un espace stratifié $\Sc h_{X}$ et une catégorie $\Ss_{\Sigma}$ dont les objets sont donnés par :
\begin{itemize}
\item pour toute strate, un faisceau sur cette strate,
\item pour tout couple de strate, un morphisme de recollement.
\end{itemize} 
C'est une généralisation du fait qu'un faisceau est défini de manière unique par ses restrictions à un ouvert et à son fermé complémentaire ainsi qu'un morphisme de recollement. Ce résultat et sa preuve sont délibérément formulés dans un langage très formel pour permettre de les adapter aux champs. \\

La deuxième partie de ce chapitre est un résultat sur le recollement de faisceaux pervers qui est une variante explicite d'un résultat démontré par MacPherson et Vilonen dans \cite{McV}. Nous n'utilisons pas cette partie dans la suite de la thèse. \\ 

Le chapitre $3$ est une étude de la $2$-catégorie des champs sur un espace stratifié $\SSS t_{X}$.\\
Dans un premier temps, nous généralisons le théorème démontré dans la première partie du chapitre précédent. On définit une $2$-catégorie, $\SSS_{\Sigma}$, sur le même modèle que $\Ss_{\Sigma}$. Ainsi un objet de cette $2$-catégorie est donné par :
\begin{itemize}
\item pour toute strate, $S_{k}$, un champ $\CCC_{k}$ sur cette strate,
\item pour tout couple de strates $(S_{k}, S_{l})$ tel que $S_{k}\subset \overline{S}_{l} $, un foncteur de champ $F_{lk} : \CCC_{k} \rightarrow i_{k}^{-1}i_{l*}\CCC_{l}$
où $i_{k}$ et $i_{l}$ sont les injections de respectivement $S_{k}$ et $S_{l}$ dans $X$. 
\item pour tout triplet de strates, un isomorphisme $\theta_{klm}$ de foncteurs de champs :
$$\shorthandoff{;:!?}
\xymatrix @!0 @C=1.5cm @R=0.4cm {\CCC_{k} \ar[rrr]^{F_{lk}}  \ar[ddddd]_{F_{mk}} &&& i_{k}^{-1}i_{l*}\CCC_{l} \ar[ddddd]^{i_{k}^{-1}i_{l*}F_{ml}} \\
~\\
& & \ar@{=>}[ld]_\sim^{\theta_{klm}}\\
&~\\
\\
i_{k}^{-1}i_{m*}\CCC_{m} \ar[rrr]_{i^{-1}_k\eta_{lm}} &&& i_{k}^{-1}i_{l*}i_{l}^{-1}i_{m*}\CCC_{m}
}$$
qui vérifient des conditions de commutations.
\end{itemize}
Nous démontrons le théorème suivant.
\begin{thm}
Les deux $2$-catégories $\SSS t_{X}$ et $\SSS_{\Sigma}$ sont $2$-é\-qui\-va\-len\-tes.
\end{thm}
La démonstration du théorème s'inspire de notre démonstration dans le cas des faisceaux (voir chapitre 2), mais devient con\-si\-dé\-ra\-ble\-ment plus technique. Néanmoins l'idée de base reste simple : un champ sur $X$ est $2$-limite projective de ses restrictions sur les strates où le système projectif doit coder les conditions de recollement. Remarquons que, contrairement au cas des faisceaux ou une démonstration élémentaire est possible, dans le cas des champs, l'approche formel est, dans la pratique, incontournable.\\

Dans la deuxième partie de ce chapitre nous introduisons la notion de champs constructibles et strictement constructibles relativement à une stratification fixée. Comme pour les faisceaux, un champ est constructible (resp. strictement constructible) si  la restriction de ce champ à chaque strate est localement constante (resp. constante). Cette notion intervient parce que le champ $\PPP_{X}$ est constructible et dans quelques cas même strictement constructible. On s'intéresse ensuite au cas particulier où $X=\CC^n$ et $\Sigma$ est la stratification, dite du croisement normal, associée à l'ensemble $\{ (z_{1},\ldots,z_{n}) \in \CC^n| z_{1}\ldots z_{n}=0\}$. On démontrera  que dans ce cas le champ $\PPP_{\CC^n}$ est strictement constructible et qu'un champ strictement constructible sur $\CC^n$ est défini de manière unique (à isomorphisme près) par la donnée :
\begin{itemize}
\item pour toute strate $S_{K}$, d'une catégorie,
\item pour tout couple de strates $(S_{k},S_{l})$ telles que $S_{k}\subset \overline{S}_{l}$, d'un foncteur $F_{lk}$,
\item pour tout triplet $(S_{k}, S_{l}, S_{m})$ de strates tel que $S_{k} \subset \overline{S}_{l}\subset \overline{S}_{m}$, d'un isomorphisme de foncteurs.
\end{itemize}
Cette description particulièrement simple est déduit du théorème 0.1 et s'appuie sur le fait que les voisinages tubulaires des strates sont des fibrations triviales. L'argument clef est donnée par la proposition :
\begin{prop}
Si $\CCC_{L}$ est un champ constant sur la strate $S_{L}$ alors le champ $i_{K}^{-1}i_{L*}\CCC_{L}$ est constant, où $i_{K}$ et $i_{L}$ sont les injections des strates $S_{K}$ et $S_{L}$ dans $X$.
\end{prop}
Dans le cas du champ $\PPP_{\CC^n}$, on obtient que la restriction de ce champ aux strates de dimension $k$ est constante de fibre $\Pc erv_{\CC^{n-k}}$, et que la fibre du champ $i_{K}^{-1}i_{L*}i_{L}^{-1}\PPP_{\CC^n}$ est équivalente à la catégorie $\Pc erv_{\CC^{*n-l}\times\CC^{l-k}}$. On a alors la description suivante du champ des faisceaux pervers :
\begin{thm}
Le champ $\PPP_{\CC^n}$ est équivalent à la donnée :
\begin{itemize}
\item pour toute strate $S_{K}$ de dimension $k$, de la catégorie $\Pc erv_{\CC^{n-k}}$,
\item pour tout couple de strates $(S_{K}, S_{L})$ tel que $S_{K}\subset \overline{S}_{L}$, de la restriction :
$$\Pc erv_{\CC^{n-k}} \rightarrow \Pc erv_{\CC^{*n-l}\times\CC^{l-k}} $$
\item pour tout triplet $(S_{k}, S_{l}, S_{m})$ de strates tel que $S_{k} \subset \overline{S}_{l}\subset \overline{S}_{m}$, un isomorphisme de foncteur.
\end{itemize}
\end{thm}

~\\

Dans le quatrième chapitre  nous généralisons l'équivalence de Galligo, Granger et Maisonobe en une équivalence de champ. Dans un premier temps nous rappelons brièvement les résultats démontrés dans \cite{GGM}. Nous rappelons ainsi la définition des catégories de représentations de carquois, $\Cc_{n}$, puis des équivalences $\alpha_{n}$ :
$$\alpha_{n} : \Pc erv_{\CC^n} \longrightarrow \Cc_{n}$$
Dans la deuxième partie nous définissons un champ, $\CCC_{\CC^n}$ de catégories de représentations de carquois équivalent à $\PPP_{\CC^n}$. Pour cela on considère le champ définit par la donnée :
\begin{itemize}
\item pour toute strate $S_{K}$ de dimension $k$, de la catégorie $\Cc_{n-k}$,
\item pour tout couple de strates $(S_{K},S_{L})$ tel que $S_{K}\subset \overline{S}_{L}$ , d'un foncteur  ``oubli'',
\item et enfin pour tout triplet de strates $(S_{K}, S_{L}, S_{M})$ tel que $S_{K}\subset \overline{S}_{L}\subset \overline{S}_{M}$ d'un isomorphisme de foncteur.
\end{itemize}
Puis nous démontrons le théorème suivant :
\begin{thm}
Les objets $\CCC_{\CC^n}$ et $\PPP_{\CC^n}$ sont équivalents.
\end{thm}
Pour le démontrer nous considérons l'équivalence définie par la donnée pour toute strate de dimension $k$ de l'équivalence $\alpha_{n-k}$. Il s'avère techniquement difficile de démontrer que ces équivalences se recollent en un foncteur. 

~\\\\
Enfin les deux derniers chapitres sont deux applications de cette équi\-va\-lence de champs. Mais les stratégies diffèrent sensiblement.\\
Dans le chapitre $5$ nous considérons $\CC^2$ muni de la stratification donnée par un arrangement générique d'hyperplans. Cette stratification est localement un croisement normal, de plus les voisinages tubulaires des strates sont homéomorphes à des produits.  On peut alors appliquer directement la construction faite aux chapitres précédents. L'étude des sections globales du champ ainsi défini nous donne alors une équivalence entre $\Pc erv_{\Sigma}$ et une catégorie de représentation du carquois associé à la stratification. Cette étude nécessite la connaissance de la topologie de chaque strate. \\
Dans le chapitre $6$, on s'intéresse aux variétés toriques lisses stratifiées par l'action du tore. Les variétés affines lisses sont des produits de $\CC$ avec des  $\CC^*$ stratifiés par le croisement normal. On sait donc caractériser $\PPP_{X}$ sur les ouverts affines à l'aide de nos résultats du chapitre 4. Ici on montre que ces caractérisations locales se recollent en une équivalence de champ. La difficulté réside dans l'expression des isomorphismes de recollement torique et dans la compréhension de leur action sur le recollement des champs. Puis on explicite cette équivalence de champ au niveau des sections globales. Ceci nous permet d'établir une équivalence entre $\Pc erv_{X}$ et une catégorie de représentations du carquois associé à l'éventail définissant $X$.

\chapter{Définitions et notations}
Dans ce chapitre nous donnons les principales définitions et notations que nous allons utiliser. 
\section{Carquois}
On donne dans ce paragraphe les définitions de carquois, carquois associé à une stratification et représentations de carquois.
\begin{Def}
    Un carquois $Q$ est la donnée d'un ensemble de sommets $ I(Q):= \{ s_1, s_2, \cdots, s_p \}$ et pour tout couple de sommets non nécessairement différents d'un ensemble de flèches $Arr(s_{\alpha},s_{\beta})$. 
\end{Def}
\noindent
Soit $X$ un espace topologique muni d'une stratification $\displaystyle{S= \bigcup_{\alpha \in A}S_{\alpha}}$. 
\begin{Def}
 On note $Q_{S}$ le carquois suivant : 
\begin{itemize}
\item \`A toute strate, $S_{\alpha}$ on associe un sommet $s_{\alpha}$, 
\item Soient $s_{\alpha}$ et $s_{\beta}$ deux sommets, on a :
$$Arr(s_{\alpha},s_{\beta})=\left\{ \begin{array}{ll} \{a_{\alpha\beta}\}  &\text{si~ $S_{\alpha}$ est une strate incidente à $S_{\beta}$ de codimension $1$} \\
\{a_{\alpha\beta}\}  &\text{si~ $S_{\beta}$ est une strate incidente à $S_{\alpha}$ de codimension $1$} \\
                                                                                             \emptyset & \text{sinon}
                                                              \end{array} \right.$$                               
 
\end{itemize}
\end{Def}
\noindent
Ainsi, entre deux sommets, il y a soit deux flèches orientés de fa\c con opposée, soit aucune flèche.\\\\
\textbf{Exemples :}
\begin{itemize}
\item Pour X=$\CC$, muni de la stratification $S=\{0\}\cup \CC^*$ du croisement normal, $Q_{\CC}$, aussi noté $Q_{1}$, est le carquois suivant :
$$\xymatrix{ \bullet \ar@/^/[r] & \bullet \ar@/^/[l]}$$
\item Pour $X=\CC^2$ muni de la stratification $S$ du croiement normal
 $$S=(\{0\}\times \{0\})\cup( \CC^*\times \{0\} )\cup (\{0\} \times \CC^*\cup \CC^*\times \CC^*)$$
 $Q_{\CC^2}$, aussi noté $Q_{2}$,  est le carquois suivant : 
$$
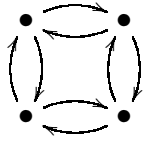
\xymatrix{ \bullet \ar@/^/[r]  \ar@/^/[d] & \bullet  \ar@/^/[l]  \ar@/^/[d]\\
                        \bullet  \ar@/^/[r]  \ar@/^/[u] & \bullet  \ar@/^/[l]  \ar@/^/[u]\\
}$$
\item Pour $X=\CC^n$ muni de la stratification $S$ associée au croisement normal  $\{(z_{1},\ldots ,z_{n}) \in \CC^n |z_{1}z_{2}\cdots z_{n}=0\}$, le carquois $Q_{n}$ est un hypercube de dimension $n$ avec 2 flèches sur chacune des arrêtes. 
\item Pour $X=\CC^2$ muni de la stratification $S'$ donnée par trois droites, $D_{1}$, $D_{2}$ et $D_{3}$ en position générique
$$S'_{\emptyset}=\CC^2 \backslash (\cup_{i} D_{i}),~~ S'_{i}= D_{i} \backslash (\cup_{j\neq i} D_{j}),~~ S'_{ij}= D_{i}\cap D_{j},$$
 le carquois $Q_{S'}$ est le suivant : 
$$
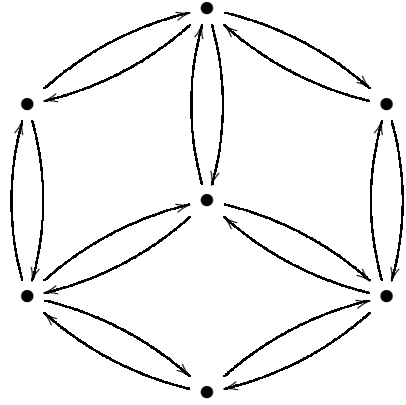
\xymatrix{    & & \bullet \ar@/^/[dll]^{}  \ar@/^/[drr]  \ar@/^/[dd]  \\
                           \bullet  \ar@/^/[urr] \ar@/^/[dd]  & &  & & \bullet   \ar@/^/[ull]  \ar@/^/[dd] \\
                           & & \bullet \ar@/^/[uu]  \ar@/^/[dll]  \ar@/^/[drr]   \\
                           \bullet  \ar@/^/[uu]  \ar@/^/[urr]  \ar@/^/[drr]   & & & & \bullet  \ar@/^/[uu]  \ar@/^/[ull]  \ar@/^/[dll]   \\
                           & & \bullet  \ar@/^/[ull]  \ar@/^/[urr]  \\   }$$ 
\item Pour $X=\PP_{1}$ muni de la stratification $S$ suivante : 
$$S_{\emptyset}= \{[x:y] \mid x\neq 0, y \neq 0\},$$
$$ S_{1}=\{[x:0] \mid x\neq 0\}, ~~    S_{2}=\{[0:y] \mid y\neq 0\}, $$     
le carquois $Q_{\PP_{1}}$ est le suivant : 
$$\xymatrix{
\bullet \ar@/^/[r] & \bullet \ar@/^/[l] \ar@/^/[r] & \bullet \ar@/^/[l]
}$$  
\item Pour $X= \PP_{2}$ stratifié par l'action du tore, ie la stratification associée à $\{[z_{1} : z_{2}:z_{3}] \in \PP_{2}| z_{1}z_{2}z_{3}=0\}$,  le carquois $Q_{\PP_{2}}$ est le même que $Q_{S'}$.       
\end{itemize}
\begin{Def}
Une représentation d'un carquois $Q$ est la donnée
\begin{itemize}
\item pour tout sommet $s_{\alpha}$ de $Q$ d'un espace vectoriel $M_{\alpha}$ de dimension finie
\item pour toute arrête $a \in Arr(s_{\alpha}, s_{\beta})$, d'une application linéaire $M_{a}$ :
$$M_{a} : E_{\alpha}\longrightarrow E_{\beta}$$
\end{itemize}
\end{Def}
\begin{Def}
Soient $R=(\{E_{\alpha}\}, \{M_{a}\})$ et $T=(\{E'_{\alpha}\}, \{M'_{a}\})$ deux représentations d'un même carquois $Q$, un morphisme $f$ de R dans T est la donnée pour tout sommet $s_{\alpha}$ de $Q$ d'un ensemble d'applications linéaires :
$$f_{\alpha}: E_{\alpha} \rightarrow E'_{\alpha} $$
telles que, pour toute arrête $a \in Arr(s_{\alpha},s_{\beta} )$, le diagramme suivant commute : 
$$\xymatrix{
E_{\alpha} \ar[d]^{M_{a}} \ar[r]^{f_{\alpha}}  & E'_{\alpha} \ar[d]^{M'_{a}}\\
E_{\beta}  \ar[r]^{f_{\beta}}  & E'_{\beta} 
}$$
\end{Def}
L'ensemble des représentations d'un carquois muni des morphismes est une catégorie abélienne. Pour un carquois Q donné, on notera $R(Q)$ cette catégorie.\\\\

\section{Faisceaux pervers}
Dans ce paragraphe nous rappelons rapidement  la définition de la catégorie des faisceaux pervers sur un espace stratifié de Thom-Mather. Pour plus de précisions le lecteur pourra se rapporter à \cite{BBD}, \cite{KS} ou \cite{Dim}.\\

Soit $X$ un espace de Thom-Mather, muni de la stratification $\Sigma$, pour les définitions d'espace de Thom-Mather le lecteur peut se reporter à  \cite{T} ou \cite{Ma}. Nous ne définissons la catégorie de faisceaux pervers que pour la perversité moyenne, ainsi nous supposons que tous les espaces sont de dimension pair.\\
Notons $D^b(\Ss h_{X})$ la catégorie dérivée des complexes de faisceaux sur $X$ dont la cohomologie est bornée. 

\begin{Def}
La catégorie $\Pc erv_{X}$ des faisceaux pervers relativement à la stratification $\Sigma$ est la sous-catégorie pleine de $D^b(\Ss h_{X})$ dont les objets sont les complexes de faisceaux $\Fc^\bullet$ satisfaisant les conditions suivantes :
\begin{itemize}
\item $H^k(i_{\Sigma_{j}}^{-1}\Fc^\bullet)$ est un système local de rang fini sur toute strate,
\item  $H^k(i_{\Sigma_{j}}^{-1}\Fc^\bullet)=0$ pour tout $k>\frac{n-\dim(\Sigma_{j})}{2}$,
\item  $H^k(i_{\Sigma_{j}}^{!}\Fc^\bullet)=0$  pour tout $k<\frac{n-\dim(\Sigma_{j})}{2}$.
\end{itemize}
pour toute strate $\Sigma_{j}$ et où $i_{\Sigma_{j}}$ est l'injection $\Sigma_{j} \hookrightarrow X$.
\end{Def}
Notons que dans la littérature la convention suivante existe aussi :
\begin{itemize}
\item $H^k(i_{\Sigma_{j}}^{-1}\Fc^\bullet)=0$ pour tout $k>-\frac{\dim(\Sigma_{j})}{2}$,
\item  $H^k(i_{\Sigma_{j}}^{!}\Fc^\bullet)=0$  pour tout $k<-\frac{\dim(\Sigma_{j})}{2}$.
\end{itemize}

Dans \cite{GGM}, A. Galligo, M. Granger, Ph Maisonobe, donne une formulation équivalente de ces conditions de perversité :

\begin{prop}
Un complexe $\Fc$ de $D^b(\Ss h_{X})$ à cohomologie contructible relativement à $\Sigma$ appartient à $\Pc erv_{\Sigma}$ si et seulement si il vérifie les trois conditions suivantes :
\begin{itemize}
\item $\forall i \notin \{0, 1, \cdots, \frac{n}{2}\} h^{i}(\Fc)=0$\\
\item le support du faisceau $h^{i}(\Fc)$ est contenu dans :
$$\overline{\Sigma}_{n-i}=\bigcup_{0\leq j\leq \frac{n-i}{2}}\Sigma_{j}$$
\item le complexe $(R\Gamma_{\Sigma_{\frac{n-j}{2}}\Fc})\mid_{\Sigma_{\frac{n-j}{2}}}$ est concentré en degré supérieur ou égal à $j$.
\end{itemize}
\end{prop}

\section{$2$-catégories, $2$-limites}
Dans ce chapitre on rappelle la définition d'une $2$-catégorie, d'un $2$-foncteur et d'une $2$-transformation. Puis on donne la définition d'une $2$-limite projective et inductive. Dans la littérature ce que nous nommons $2$-foncteur est souvent appelé pseudo-foncteur, notamment dans \cite{Bo}. De même, les $2$-limites sont souvent définies à partir d'un $2$-foncteur qui a pour $2$-catégorie source une simple catégorie. Nous donnons ici la définition plus générale à partir d'un $2$-foncteur qui a pour source et pour but une $2$-catégorie. 
Ce chapitre s'inspire largement du chapitre 7 de \cite{Bo} et de l'annexe de l'article  \cite{Ingo1}. L'article \cite{St} est aussi une référence standard.

\begin{Def}
Une $2$-catégorie est la donnée :
\begin{itemize}
\item d'une classe d'objets : $\mathbf{A,~B},\ldots$ 
\item d'une classe de $1$-morphismes que nous appelons foncteurs : $$ \xymatrix{\mathbf{A} \ar[r]^\Phi & \mathbf{B}}$$
\item d'une classe de $2$-morphismes que nous appelons morphismes de foncteurs  : \xymatrix{\relax
    \mathbf{A} \UN[r]{\Phi}{\Psi}{\alpha} & \mathbf{B} }
\end{itemize}
satisfaisant les axiomes suivants :
\begin{itemize}
\item[$\bullet$] Si $\Phi$ et $\Psi$ sont deux $1$-morphismes :
$$ \xymatrix{\mathbf{A} \ar[r]^\Phi & \mathbf{B} \ar[r]^\Psi & \mathbf{C}}$$ 
alors la composition :
$$\xymatrix{\mathbf{A} \ar[r]^{\Psi\Phi} & \mathbf{C}} $$
existe et est associative.
\item[$\bullet$] Si $\alpha$ et $\beta$ sont deux $2$-morphismes :
$$\shorthandoff{;:!?}
\xymatrix@!0 @R=2cm @C=4cm {\relax
    \mathbf{A} \DEUX[r]{\Phi}{\Phi'}{\Phi''}{\alpha}{\beta} & 
    \mathbf{B} }$$
alors la composée
$$\shorthandoff{;:!?}
\xymatrix@!0 @R=2cm @C=4cm {\relax
    \mathbf{A} \UN[r]{\Phi}{\Psi}{\beta\circ\alpha} & \mathbf{B} } $$
est définie et est associative. On appel cette composée, composée verticale.
\item[$\bullet$] Si $\alpha$ et $\beta$ sont deux $2$-morphismes :
$$\shorthandoff{;:!?}
\xymatrix@!0 @R=2cm @C=3cm {\relax\mathbf{A} \UN[r]{\Phi}{\Phi'}{\alpha} & \mathbf{B} \UN{\Psi}{\Psi'}{\beta} &\mathbf{C}}$$
alors la composée :
$$\shorthandoff{;:!?}
\xymatrix@!0 @R=2cm @C=4cm {\relax \mathbf{A} \UN{\Psi\Phi}{\Psi'\Phi'}{\beta\bullet\alpha} &\mathbf{C}}$$ est définie et associative. On l'appel la composée horizontale.
\item[$\bullet$] Si $\alpha$, $\alpha'$, $\beta$ et $\beta'$ sont quatre $2$-morphismes :
$$\shorthandoff{;:!?}
\xymatrix@!0 @R=2cm @C=4cm {\relax
\mathbf{A} \DEUX[r]{\Phi}{\Phi'}{\Phi''}{\alpha}{\alpha'} &\mathbf{B}\DEUX[r]{\Psi}{\Psi'}{\Psi''}{\beta}{\beta'} &\mathbf{C}}$$
alors on a l'égalité :
$$(\beta' \circ \beta)\bullet(\alpha'\circ \alpha)=(\beta'\bullet\alpha')\circ(\beta\bullet\alpha)$$
\item[$\bullet$] Pour tout objet $\mathbf{A}$ il existe un $1$-morphisme \xymatrix{\mathbf{A} \ar[r]^{Id_{\mathbf{A}}} & \mathbf{A}} tel que pour tout $1$-morphisme \xymatrix{\mathbf{A}\ar[r]^\Phi& \mathbf{B}} on ait l'égalité :
$$\Phi Id_{\mathbf{A}}=\Phi=Id_{\mathbf{B}}\Phi$$
et pour tout  $1$-morphisme \xymatrix{\mathbf{A}\ar[r]^\Phi& \mathbf{B}} il existe un $2$-morphisme \xymatrix{\Phi \ar@{=>}^{Id_{\Phi}}[r] & \Phi} tel que :
$$\alpha\circ Id_{\Phi} = \Phi= Id_{\Phi}\circ \alpha$$
$$Id_{\Psi}\bullet Id_{\Phi}= Id_{\Psi \Phi} $$
\end{itemize}
\end{Def}
\noindent
\textbf{Remarque.}\\
Si $\cat A$ est une $2$-catégorie et $\mathbf{A}$ et $\mathbf{B}$ sont deux objets de $\cat A$, alors les $1$-morphisme de $\cat A$ entre $\mathbf{A}$ et $\mathbf{B}$ et les $2$-morphismes de $\cat A$ entre de tels $1$-morphismes  forment une catégorie notée $\cat A (\mathbf{A}, \mathbf{B})$.\\
\textbf{Exemples}\\
\begin{itemize}
\item
La catégorie $\mathbf{Cat}$ est une $2$-catégorie dont les objets sont les petites catégories, les $1$-morphismes les foncteurs et les $2$-morphismes les transformations naturelles.
\item Toute catégorie $\mathbf{C}$ peut être vu comme une $2$-catégorie notée $\hat{\mathbf{C}}$ :
$$\hat{\mathbf{C}}~~~ \left\{\begin{array}{l}
                                          \text{Les objets de $\hat{\mathbf{C}}$ sont ceux de $\mathbf{C}$.}\\
                                          \text{Les $1$-morphismes de $\hat{\mathbf{C}}$ sont les flèches de $\mathbf{C}$.}\\
                                          \text{Les $2$-morphismes de $\hat{\mathbf{C}}$ sont les identités.}
                                          \xymatrix{Id_{f} : f\ar@{=>}[r] & f}. 
                                         \end{array} \right.$$
\end{itemize}
\begin{Def}
\'Etant donné deux $2$-catégories $\cat A$ et $\cat B$, un $2$-foncteur $\cat f : \cat A \rightarrow \cat B$ est donné par : 
\begin{itemize}
\item pour tout objet $\mathbf{A}$ de $\cat A$, un objet $\cat f  (\mathbf{A})$ de $\cat B$,
\item pour tout $2$-morphisme $\Phi :\mathbf{A} \rightarrow \mathbf{B}$, un $2$-morphisme de $\cat B$, $$\cat f(\Phi) : \cat f(\mathbf{A}) \rightarrow \cat f(\mathbf{B}),$$
\item pour toute transformation naturelle \xymatrix{ \Phi \ar@{=>}[r]^\alpha & \Psi } de $\cat A$, une transformation naturelle de $\cat B$, \xymatrix{ : \cat f(\Phi) \ar@{=>}[r]^{\cat f(\alpha)} & \cat f(\Psi)},
\item pour tout objet $\AAA$ de $\cat A$, un $1$-isomorphisme $\cat f_{\AAA}$ :
$$\cat  f_{\AAA}: Id_{\cat f(\AAA)} \buildrel\sim\over\longrightarrow\cat f(Id_{\AAA})  $$
\item pour tout couple de $1$-morphismes \xymatrix{\mathbf{A} \ar[r]^\Phi & \mathbf{B} \ar[r]^\Psi & \mathbf{C}} de $\cat A$, un $2$-isomorphisme $\cat f(\Psi, \phi)$ de $\cat B$ :
$$\shorthandoff{;:!?}
\xymatrix@!0 @R=2cm @C=4cm {\relax \cat f( \mathbf{A}) \UN{\cat f(\Psi)\cat f(\Phi)}{\cat f(\Psi \phi)}{\cat f (\Psi, \Phi)} & \cat f(\mathbf{C}) },$$
\end{itemize}
qui vérifient les conditions suivantes :
\begin{itemize}
\item[$\bullet$] $\cat f (Id_{\Phi})=Id_{\cat f(\Phi)} $.
\item[$\bullet$] Le diagramme suivant commute :
$$\shorthandoff{;:!?}
\xymatrix@!0 @R=2cm @C=5cm {\relax
\cat f(\Theta)\cat f(\Psi)\cat f(\Phi) \ar@{=>}[r]^{\cat f(\Theta, \Psi) \bullet Id_{\cat f(\Phi)}} \ar@{=>}[d]_{Id_{\cat f(\Theta)}\bullet \cat f(\Psi, \Phi)} & \cat f(\Theta \Psi)\cat f(\Phi) \ar@{=>}[d]^{\cat f(\Theta\Psi, \Phi)}\\
\cat f(\Theta)\cat f(\Psi\Phi) \ar@{=>}[r]_{\cat f (\Theta, \Psi \Phi)} & \cat f(\Theta\Psi\Phi)
.}$$
\item[$\bullet$] Si $\alpha$ et $\beta$ sont deux transformations naturelles :
$$\shorthandoff{;:!?}
\xymatrix@!0 @R=2cm @C=4cm {\relax
\mathbf{A} \DEUX{\Phi}{\Phi'}{\Phi''}{\alpha}{\beta} &\mathbf{B} }$$
alors $\cat f(\beta \circ \alpha)= \cat f(\beta )\circ \cat f(\alpha)$.
\item Si $\alpha$ et $\beta$ sont deux transformations naturelles : 
$$\shorthandoff{;:!?}
\xymatrix@!0 @R=2cm @C=4cm {\relax
    \mathbf{A} \UN[r]{\Phi}{\Psi}{\beta\circ\alpha} & \mathbf{B} } $$
alors le diagramme suivant commute :
$$\shorthandoff{;:!?}
\xymatrix@!0 @R=2cm @C=4cm {\relax
\cat f(\Psi)\cat f(\Phi) \ar@{=>}[r]^{\cat f(\beta)\bullet\cat f(\alpha)} \ar@{=>}[d]_{\cat f(\Psi, \Phi)} & \cat f(\Psi')\cat f(\Phi') \ar@{=>}[d]^{\cat f(\Psi', \Phi')} \\
\cat f(\Psi \Phi) \ar@{=>}[r]_{\cat f(\beta\bullet\alpha)} & \cat (\Psi'\Phi').
}$$
\item Si $\phi: \AAA \rightarrow \BBB$ est un $1$-morphisme, alors le diagramme suivant commute :
$$\shorthandoff{;:!?}
\xymatrix@!0 @R=3cm @C=8cm {\relax
 \cat f(\phi) \ar[r]^{\cat f(\phi)\bullet \cat f_\AAA} \ar[d]_{\cat f_\BBB\bullet \cat f(\phi)}^\sim \ar[rd]^{Id_{\cat f(\Phi)}} & \cat f(\phi)\cat f(Id_\AAA) \ar[d]_\sim^{\cat f(Id_\AAA,\phi)} \\
\cat f(Id_\BBB)\cat f(\phi) \ar[r]_{\cat f(\phi, Id_\BBB)}^\sim &   \cat f(\phi)=\cat f(Id_\BBB \circ \phi)=\cat f(\phi \circ Id_\BBB)
}$$
\end{itemize}
\end{Def}
\textbf{Remarque.}\\
Pour tout couple $( \mathbf{A}, \mathbf{B})$ d'objets  de $\cat A$ la donnée :
\begin{itemize}
\item pour tout $1$-morphisme $\Phi$ du $1$-morphisme $\cat f(\Phi)$, 
\item pour tout $2$-morphisme $\alpha$ du $2$-morphisme $\cat f(\alpha)$,
\end{itemize}
est un foncteur de la catégorie $\cat A(\mathbf{A}, \mathbf{B})$ dans la catégorie $\cat B\big( \cat f(\mathbf{A}), \cat f(\mathbf{B})\big)$ que l'on note $\cat f_{\mathbf{A}\mathbf{B}}$.\\
On donne maintenant la définition d'une $2$-transformation naturelle. 
\begin{Def}
Soient $\cat f : \cat A \rightarrow \cat B$ et $\cat g : \cat A \rightarrow \cat B $ deux $2$-foncteurs entre les $2$-catégories $\cat A$ et $\cat B$. Une $2$-transformation naturelle        $\alpha :  \cat f \Rightarrow \cat g$ est donnée par :
\begin{itemize}
\item pour tout  objet $\mathbf{A}$ de $\cat A$ un $1$-morphisme $\alpha_{\mathbf{A}} : \cat f (\mathbf{A}) \rightarrow \cat g(\mathbf{A})$,
\item pour tout $1$-morphisme $\phi : \mathbf{A} \rightarrow \mathbf{B}$ un $2$-isomorphisme $\alpha_{\phi}$ :
$$\shorthandoff{;:!?}
\xymatrix@!0 @R=2cm @C=4cm {\relax
   F( \mathbf{A}) \UN[r]{G(\phi)\circ \alpha_{\AAA}}{\alpha_{\BBB}\circ F(\phi)}{\alpha_{\phi}} & G(\mathbf{B}) } $$

\end{itemize}
Ces données doivent de plus vérifiées les conditions suivantes :
\begin{itemize}
\item[$\bullet$] Pour tout $2$-morphisme :
$$\shorthandoff{;:!?}
\xymatrix@!0 @R=2cm @C=4cm {\relax
    \mathbf{A} \UN[r]{\Phi}{\Psi}{\varepsilon} & \mathbf{B} } $$    
le diagramme suivant commute :
$$\xymatrix{
\cat g(\phi)\circ \alpha_\AAA \ar[r]^{\alpha_\phi} \ar[d]_{\cat g(\varepsilon)\bullet Id_{\alpha_\AAA}} & \alpha_\BBB \circ \cat f(\phi) \ar[d]^{Id_{\alpha_\BBB}\bullet \cat f(\varepsilon)}\\
\cat g(\psi)\circ \alpha_\AAA \ar[r]_{\alpha_\psi} & \alpha_\BBB \circ \cat f(\psi)
}$$
\item[$\bullet$] Pour tout objet $\AAA$ de $\cat A$ le diagramme suivant commute :
$$\xymatrix{
& \alpha_\AAA \ar[dl]_{\cat f_\AAA\bullet Id_{\alpha_\AAA}} \ar[dr]^{Id_{\alpha_\AAA}\bullet \cat g_\AAA}\\
\cat f(Id_\AAA) \circ \alpha_\AAA \ar[rr]_{\alpha_{Id_\AAA}} && \alpha_\AAA \circ \cat g(Id_\AAA)
}$$
\item pour tout couple de $2$-morphismes composables $\phi : \AAA \rightarrow \BBB$ et $\psi : \BBB \rightarrow \ccc$, le diagramme suivant commute :
$$\shorthandoff{;:!?}
\xymatrix@!0 @R=2cm @C=4cm {\relax
\cat f(\psi)\cat f(\phi)\alpha_\AAA  \ar[d]^{\cat f(\psi,\phi)} \ar[r]^{Id_{\cat f(\psi)} \bullet \alpha_\phi} & \cat f(\psi)  \alpha_\BBB \cat g(\phi) \ar[r]^{\alpha_\psi\bullet Id_{\cat g(\phi)}}& \alpha_\ccc\cat g(\psi)\cat f(\phi)  \ar[d]^{Id_{\alpha_\ccc}\bullet\cat g(\psi\phi)}\\
\cat f(\psi\phi)\alpha_\AAA \ar[rr]_{\alpha_{\psi\phi}} && \alpha_\ccc\cat g(\psi\phi)
}$$
\end{itemize}
\end{Def}
On donne maintenant la définition d'une $2$-limite. \\
\begin{Def}
Soit $\III$ une $2$-catégorie, et $\cat a: \III\rightarrow \cat B $ un $2$-foncteur.
\begin{itemize}
\item[$\bullet$] Le système $\cat a$ admet une $2$-limite si et seulement si il existe :
\begin{itemize}
\item un objet de $\cat B$ noté $2\varprojlim_{i\in \III}\cat a(i)$ et 
\item une $2$-transformation naturelle notée $\sigma$ entre le $2$-foncteur constant $2\varprojlim_{i\in \III}\cat a(i)$ et le $2$-foncteur $\cat a$,
\end{itemize}
tel que pour tout objet $\BBB$ de $\cat B$ le foncteur :
$$( \sigma \circ) : \text{\underline{Hom}}_{\cat B}(2\varprojlim_{i\in \III}\cat a(i), \BBB) \rightarrow \text{\underline{Hom}}(\cat a, \BBB)$$
est une équivalence de catégorie. 
\item[$\bullet$] Le système $\cat a $ admet une $2$-colimite si et seulement si :
\begin{itemize}
\item un objet de $\cat B$ noté $2\varinjlim_{i\in \III}\cat a(i)$ et 
\item une $2$-transformation naturelle notée $\sigma$ entre le $2$-foncteur $\cat a$ et le $2$-foncteur constant $2\varinjlim_{i\in \III}\cat a(i)$,
\end{itemize}
tel que pour tout objet $\BBB$ de $\cat B$ le foncteur :
$$(\circ \sigma) : \text{\underline{Hom}}_{\cat B}(2\varinjlim_{i\in \III}\cat a(i), \BBB) \rightarrow \text{\underline{Hom}}(\cat a, \BBB)$$
est une équivalence de catégorie. 
\end{itemize}
\end{Def}
C'est en fait la définition la plus générale de la $2$-limite, mais dans la littérature le cas le plus souvent traité est celui où $\cat a $ est un $2$-foncteur de source une petite catégorie. Dans la suite nous n'utiliserons d'ailleurs la $2$-limite que dans ce cas là. \\
Considérons plus en détail la définition d'une $2$-limite dans le cas où $\III$ est une petite catégorie. Soit  $\cat B$ une $2$-catégories et $\cat a: \III \rightarrow \cat B$ un $2$-foncteur. Alors $\cat a$ admet une $2$-limite si et seulement s'il existe :
\begin{itemize}
\item un objet de $\cat B$, $2\varprojlim_{i \in \III}\cat a(i)$;
\item pour tout objet $i$ de $\III$, un $2$-morphisme $$\pi_{i} : 2\varprojlim_{i \in \III}\cat a(i) \rightarrow \cat a(i)~~;$$
\item et, pour tout morphisme de $\III$, $s : i \rightarrow j$, un $1$-morphisme, noté \linebreak$\theta_{s} : \cat a(s)\circ \pi_{i} \buildrel\sim\over\rightarrow \pi_{j}$ :
$$\shorthandoff{;:!?}
\xymatrix @!0 @C=0.7cm @R=0.7cm {
&&&2\varprojlim_{i\in \III}\cat a(i)  \ar[ddddlll]_{\pi_{i}}  \ar[ddddrrr]^{\pi_{j}}\\
~\\
&&&&~\\
&&\ar@{=>}[urr]_{\theta_{s}}^\sim&&\\
\cat a(i) \ar[rrrrrr]_{\cat a(s)}&&&&&&\cat a(j)\\
}$$
\end{itemize}
Ces données doivent vérifiées les conditions suivantes :
\begin{itemize}
\item[$\bullet$] pour tout objet $i$ de $\III$ le diagramme suivant commute :
$$\shorthandoff{;:!?}
\xymatrix@!0 @R=1.5cm @C=2cm {\relax
 \cat \pi_{i} \ar[rr]^{Id_{\pi_{i}}}\ar[rd]_{\cat a_{i}\bullet Id_{\pi_{i}}} && \pi_{i} \\
 &\cat a(Id_{i})\circ \pi_{i} \ar[ru]_{\theta_{Id_{i}}}
 }$$
 \item[$\bullet$] pour tout couple de morphismes composables $s :i\rightarrow j $ et $t : j \rightarrow k$, le diagramme suivant commute :
$$\shorthandoff{;:!?}
\xymatrix@!0 @R=2cm @C=4cm {\relax
\cat a(t)\cat a(s)\pi_{i} \ar[r]^{Id_{\cat a(t)}\bullet \theta_{s}} \ar[d]_{\cat a(t,s)\bullet Id_{\pi_{i}}}  & \cat a(t)\pi_{j} \ar[d]^{\theta_{t}}\\
\cat a(ts)\pi_{i} \ar[r]_{\theta_{ts}} & \pi_{k}
}$$ 
\end{itemize}
Ces conditions sont la traduction du fait que $\sigma$ soit une $2$-transformation naturelle. Mais ces données doivent vérifiées les conditions suivantes, la première traduisant que $(\sigma \circ)$ est essentiellement surjective et la deuxième que $(\sigma \circ )$ est pleinement fidèle.
\begin{itemize}
\item Pour tout objet $\BBB$ de $\cat B$ et toute transformation naturelle $\rho : \BBB\rightarrow \cat a$ il existe un $1$-morphisme $F$ :
$$F : \BBB \longrightarrow 2\varprojlim_{i\in \III}\cat a(i)$$
et pour tout objet $i$ de $\III$ une $2$-équivalence $\varphi_{i} : \sigma_{i}F \buildrel\sim\over\rightarrow \pi_{i}$
$$\shorthandoff{;:!?}
\xymatrix @!0 @C=0.7cm @R=0.7cm {
&&&\BBB  \ar[ddddlll]_{F}  \ar[ddddrrr]^{\rho_{i}}\\
~\\
&&&&\ar@{=>}[dll]_{\theta_{s}}^\sim\\
&&&&\\
2\varprojlim_{i\in \III} \cat a(i) \ar[rrrrrr]_{\pi_{i}}&&&&&&\cat a(i)\\
}$$
qui vérifient les relations de compatibilité suivantes :
$$\varphi_{j}\circ \theta_{s}^\rho= (\theta_{s}\bullet Id_{F})\circ (Id_{\cat a(s)}\bullet \varphi_{i})$$
c'est à dire que le diagramme suivant commute :
$$\shorthandoff{;:!?}
\xymatrix @!0 @C=4cm @R=2cm {
\cat a(s)\circ \rho_{i} \ar[r]^{\theta_{s}^\rho} \ar[d]_{Id_{\cat a(s)} \bullet \phi_{i}} & \rho_{j}\ar[d]^{\varphi_{j}}\\
\cat a(s)\circ \pi_{i}\circ F \ar[r]_{\theta_{s}\bullet Id_{F}} & \pi_{j}\circ F
}$$
\end{itemize}
\begin{prop}
Soit $\III$ une petite catégorie et $\cat a : \III \rightarrow \CCC at$  un $2$-foncteur, alors  le système $\cat a$ admet une $2$-limite et une $2$-colimite.
\end{prop}
Nous ne démontrons pas cette proposition mais nous rappelons la définition explicite  des catégories $2$-colimite.\\
La $2$-colimite est la catégorie dont les objets sont donnés par :
\begin{itemize}
\item pour tout objet $i$ de $\III$ un objet $X_{i}$ de $\cat a(i)$,
\item pour tout morphisme $s : i \rightarrow j$ de $\III$ un $1$-isomorphisme de $\cat B$ :
$$\vartheta_{s}^i : X_{i} \longrightarrow \cat a(s)X_{j}$$
\end{itemize}
tels que les conditions suivantes soient respectées :
\begin{itemize}
\item[$\bullet$] pour tout objet $i$ de $\III$, le diagramme suivant commute :
$$\xymatrix{
X_{i}\ar[rd]_{Id_{X_{i}}} \ar[rr]^{\vartheta^i_{Id_{i}}} &&\cat a(Id_{i})(X_{i})\\
& X_{i}\ar[ru]_{\cat a_{X_{i}}}
}$$
\item[$\bullet$] pour deux morphismes composables $s:i\rightarrow j$, $t : j \rightarrow k$ le diagramme commute :
$$\xymatrix{
X_{i } \ar[d]_{\vartheta_{t\circ s}} \ar[r]^{\vartheta} &\cat a(s)(X_{j}) \ar[d]^{\cat a(s)(\vartheta_{t})}\\
\cat a(t \circ s)(X_{k}) & \cat a(s)\cat a(t)(X_{k}) \ar[l]_{\cat a(s,t)}
}$$
\end{itemize}
Si $\big\{ \{X_{i}\}, \{\vartheta_{s}\}\big\}$ et $\big\{ \{X'_{i}\}, \{\vartheta'_{s}\}\big\}$ sont deux objets de cette catégorie, les morphismes sont donnés par une famille de morphismes $\{\phi_{i}\}$ tels que, pour tout morphisme $s :i\rightarrow j$ le diagramme suivant commute : 
$$\xymatrix{
X_{i} \ar[r]^{\vartheta_{s}} \ar[d]_{\phi_{i}} &\cat a(s)(X_{j})\ar[d]^{\cat a(s)(\phi_{j}}\\
X'_{i} \ar[r]_{\vartheta'_{s}} & \cat a(s)(X'_{j})
}$$

\section{Champs}
Ce paragraphe est une brève introduction à la théorie des champs. Le début de ce texte suit le chapitre 3 de \cite{KS2}, la suite s'inspire de \cite{Ingo1} et \cite{Ingo2}. Les démonstrations ne sont pas données, le lecteur pourra se reporter aux références citées.\\
Soit $X$ un espace topologique. Dans tout ce paragraphe, si $U_{i}$, $U_{j}$ et $U_{k}$ sont des ouverts de $X$, alors on note $U_{ij}$ et $U_{ijk}$ les ouverts :
$$U_{ij}=U_{i}\cap U_{j} \text{~~et~~} U_{ijk}=U_{i}\cap U_{j}\cap U_{k}.$$
\begin{Def}
Un préchamp $\CCC$ sur un espace topologique $X$ est la donnée 
\begin{itemize}
  \item  pour tout ouvert $U$ de $X$, une catégorie $\CCC(U)$
  \item  pour toute inclusion d'ouverts $V \subset U$, un foncteur que l'on appelle foncteur de restriction $\rho_{VU}: \CCC(U) \rightarrow \CCC(V)$, si $U_{1}$ et $U_{2}$ sont des ouverts de $X$, on note $\rho_{12}$ le foncteur $\rho_{U_{1}U_{2}}$,
  \item  pour toute inclusion d'ouverts $W \subset V \subset U$ un isomorphisme de foncteurs $\lambda_{WVU} :\rho_{WV} \circ \rho_{VU} \buildrel\sim\over\longrightarrow \rho_{WU}$, de même si $U_{1}$, $U_{2}$ et $U_{3}$ sont des ouverts de $X$, on note $\lambda_{123}$ le morphisme $\lambda_{U_{1}U_{2}U_{3}}$.
\end{itemize}
Ces données doivent de plus satisfaire les conditions suivantes :
\begin{itemize}
  \item[(i)]  $\rho_{UU}=id_{\CCC(U)}$,
  \item[(ii)] pour toute inclusion d'ouverts $U_1 \subset U_2 \subset U_3 \subset U_4$, le diagramme suivant est commutatif :
$$\shorthandoff{;:!?}
\xymatrix @!0 @R=1,5cm @C=5cm {\relax
               \rho_{12} \circ \rho_{23} \circ \rho_{34} \ar[d]_{\lambda_{123}\bullet Id_{\rho_{34}} } \ar[r]^{Id_{\rho_{12}}\bullet\lambda_{234}}   & \rho_{12} \circ \rho_{24} \ar[d]^{\lambda_{124}}\\
	       \rho_{13} \circ \rho_{34} \ar[r]^{\lambda_{134}} &  \rho_{14}}.$$      
\end{itemize}
\end{Def}
~\\
\noindent
Si $\CCC$ et un préchamp sur $X$, $F$ un objet de $\CCC(X)$, $\theta$ un morphisme de $\CCC(X)$ et $U$ un ouvert de $X$,  on note :
$$F|_{U}= \rho_{UX}(F)\text{~~et~~} \theta|_{U}=\rho_{UX}(\theta)$$
{\bf{ Remarque : }}
La condition (ii) de la définition précédente assure que pour $F, G \in \CCC(X)$, la donnée de : 
$$ Hom_{\CCC(U)} (F|_U, G|_U),$$
pour tout ouvert $U$ de $X$, est naturellement un préfaisceau d'ensemble sur $X$, on le note $\Hc om_\CCC(F,G)$. 
\begin{Def}
On dit qu'un préchamp $\CCC$ satisfait l'axiome ST1, si pour tout ouvert $U$ de $X$ et pour tout $F,G \in \CCC(U)$, le préfaisceau $\Hc om_{\CCC|_{U}} (F, G)$ est un faisceau sur $U$.
\end{Def}

\begin{Def}
On dit qu'un préchamp $\CCC$ satisfait l'axiome ST2 si pour tout ouvert $U \subset X$, pour tout recouvrement ouvert $U= \bigcup_{i \in I} U_i$, pour toute famille $F_i \in \CCC(U_i)$, pour toute famille d'isomorphismes $\theta_{ji}: F_i|_{U_{ji}} \buildrel\sim\over\longrightarrow F_j|_{U_{ji}}$ tels que :
$$\theta_{ij}|_{U_{ijk}} \circ \theta_{jk}|_{U_{ijk}}=\theta_{ik}|_{U_{ijk}},$$
il existe $F \in \CCC(U) $ et une famille d'isomorphismes $\theta_i :F|_{U_i} \buildrel\sim\over\longrightarrow F_i$ telle que :
$$\theta|_{ij} \circ (\theta _j|_{U_{ij}})=\theta_i|_{U_{ij}}.$$
\end{Def}

\begin{Def}
Un champ est un préchamp satisfaisant les axiomes ST1 et ST2.
\end{Def}
\noindent
{\bf {Remarque :}} Si $\CCC$ est un champ et si $F$ est défini comme dans l'axiome ST2, alors $F$ est unique à isomorphisme près. En effet, si $(F', \theta_i')$ est un autre candidat, les isomorphismes $\alpha_i :\theta_i'^{-1} \circ \theta_i :F|_{U_i}  \buildrel\sim\over\longrightarrow F'|_{U_i}$ se recollent en un isomorphisme $\alpha : F \buildrel\sim\over\longrightarrow F'$ par ST1.

\begin{thm}
La donnée, pour tout ouvert, de la catégorie des faisceaux pervers sur cet ouvert et la donnée, pour tout inclusion d'ouverts, du foncteur usuel de restriction des faisceaux pervers est un champ. On le notera $\PPP_X$.
\end{thm}
Dem voir \cite{BBD}
\begin{Def}
Soient $\CCC$ et $\CCC'$ deux préchamps sur X. Un foncteur de préchamps est  donné par  :
\begin{itemize}
  \item[(i)] pour tout ouvert $U \subset X$, un foncteur $\phi_U : \CCC(U) \rightarrow \CCC'(U)$,
  \item[(ii)] pour tout inclusion d'ouverts $U \subset V \subset X$, un isomorphisme de foncteurs $\theta_{VU}: \phi_V \circ \rho_{VU} \buildrel\sim\over\longrightarrow \rho'_{VU} \circ \phi_U$, tel que pour toute inclusion $W \subset V \subset U$, le diagramme suivant commute :
  $$\xymatrix{
               \phi_W \circ \rho_{WV} \circ \rho_{VU} \ar[r]^{Id_{\phi_{W}}\bullet\lambda_{WVU}} \ar[d]_{\theta_{WV}} & \phi_W \circ \rho_{WU} \ar[dd]^{\theta_{WU}} \\
	       \rho'_{WV} \circ \phi_V \circ  \rho_{VU} \ar[d]_{\theta_{VU}} \\
	       \rho'_{WV} \circ \rho'_{VU} \circ \phi_U \ar[r]_{\lambda'_{WVU}} & \rho'_{WU} \circ \phi_U                   }$$
\end{itemize}
\end{Def}

Un foncteur de champs est un foncteur de préchamps entre les deux préchamps sous-jacents aux champs.

\begin{Def}
Soient $\CCC$ et $\CCC'$ deux préchamps sur $X$. Nous allons utiliser les mêmes notations que ci-dessus. Soit $\phi : \CCC \rightarrow \CCC'$ et $\phi' : \CCC \rightarrow \CCC'$ deux foncteurs de préchamps. Un morphisme $f : \phi \rightarrow \phi'$ de foncteurs de champs est la donnée pour tout ouvert $U$ d'un morphisme $f_U : \phi_U \rightarrow \phi'_U$ de foncteurs de $\CCC(U)$ dans $\CCC'(U)$, tels que pour tout inclusion $V \subset U$ et pour tout $F \in \CCC(U)$, le diagramme suivant commute :
$$\shorthandoff{;:!?}
\xymatrix @!0 @R=1,5cm @C=5cm {\relax
             \phi_V(\rho_{VU}F) \ar[r]^{f_V(\rho_{VU}F)} \ar[d]^{\theta_{VU}(F)}) & \phi'_V(\rho_{VU}F) \ar[d]^{\theta'_{VU}(F)})\\
	     \rho_{VU}'(\phi_{U}F) \ar[r]^{\rho_{VU}'(f_{U}F)} & \rho_{VU}'(\phi'_{U}F)}$$
\end{Def}

Ainsi, on a une notion d'équivalence de champs. On dit que deux champs $\CCC$ et $\CCC'$ sont équivalents s'il existe deux foncteurs $\phi$ et $\phi'$ :
$$\begin{array}{cccc}
\phi : &\CCC & \longrightarrow  &\CCC'\\
\phi' : &\CCC' & \longrightarrow  &\CCC\\
\end{array}$$
et deux isomorphismes $f$ et $f'$ :
$$\begin{array}{cccc}
f : &\phi\circ \phi' & \simeq & Id_{\CCC}\\
f : &\phi'\circ \phi & \simeq & Id_{\CCC'}
\end{array}$$
\textbf{Remarque.}\\
Pour un espace topologique fixé, la donnée  des champs sur $X$, des foncteurs de champs sur $X$ munis de la composition naturelle et des morphismes de foncteurs de champs munis des compositions verticales et horizontales naturelles forment une $2$-catégorie.\\

On définit maintenant la fibre d'un préchamp, $\CCC$, en un point $p$ de $X$. Considérons la catégorie, $\Tc_{p}(X)$ des voisinages de $p$. Le préchamp $\CCC$ induit un $2$-foncteur $\alpha_{p} : \Tc_{p}(X)° \rightarrow \Cc \Ac \Tc$ à valeur dans la catégorie des petites catégories. On pose 
$$\begin{displaystyle}\CCC_{p}=2\varinjlim_{U\ni p}\CCC(U)=2\varinjlim_{U\in \Tc_{p}(X)}\alpha_{p}(U)\end{displaystyle}$$
On peut donner une description explicite de cette catégorie. Les objets de $\CCC_{p}$ sont donnés par :
$$Ob~\CCC_{p}=\big\{(U,A) | ~U \text{ est un voisinage de $p$ et } A \in Ob\CCC(U)\big\}$$
Soit $(U,A)$ et $(V,B)$ deux objets de $\CCC_{p}$. Alors 
$$\begin{displaystyle}Hom_{\CCC_{p}}\big((U,A),(V,B)\big)= \varinjlim_{p\in W\subset U\cap V}Hom_{\CCC(W)}\big(A|_{W},B|_{W}\big)\end{displaystyle}$$
En particulier on a la proposition suivante :
\begin{prop}
Soit $\CCC$ un préchamp sur $X$, $p\in X$ un point, $U, V$ deux ouverts de $X$ contenant $p$ et $A \in Ob ~\CCC(U),~B \in Ob ~\CCC(V)$ deux objets. Alors $A$ et $B$ sont isomorphe dans $\CCC_{p}$ si et seulement si ils sont isomorphes dans un voisinage ouvert de $p$.
\end{prop}
Si $A \in Ob ~\CCC(U)$ on note encore $A$ son image dans $\CCC_{p}$ et si $f :A \rightarrow B$ est un morphisme dans $\CCC(U)$, on note $f_{p}$ son image dans $\CCC_{p}$. La raison de cette notation vient de la remarque suivante. \\
\textbf{Remarque : } \\
Considérons un  faisceau d'anneau $\Ac$ et le champ $\Mm od(\Ac)$. Il faut prendre garde que le foncteur naturel suivant n'est pas une équivalence :
$$\Mm od(\Ac)_{p} \longrightarrow \text{Mod}(\Ac_{p})$$
car le morphisme 
$$\Hc om_{\Ac}(\Fc,\Gc)_{p}\longrightarrow \text{Hom}_{\Ac_{p}}(\Fc_{p}, \Gc_{p})$$
n'est pas un isomorphisme en général.
\begin{prop}
 Soit $\CCC$ un préchamp sur $X$. Il existe un champ $\CCC^\dagger$ sur $X$ et un foncteur de préchamps $\eta_\CCC^\dagger : \CCC \rightarrow \CCC^\dagger$ tel que pour tout foncteur de préchamps $F : \CCC \rightarrow \CCC'$ il existe un foncteur de champ $F^\dagger$ tel que le diagramme suivant soit commutatif : 
$$\xymatrix{
\CCC \ar[r]^F \ar[d]_{\eta_\CCC^\dagger} & \CCC'\ar[d]^{\eta_{\CCC'}^\dagger}\\
\CCC^\dagger \ar[r]_{F^\dagger} & \CCC^{'\dagger}
}$$
De plus le foncteur $\eta^\dagger: \CCC \rightarrow \CCC'$ induit une équivalence de catégories sur les fibres en tout point de $X$.
\end{prop}
\begin{cor}\label{equipre}
 Soit $\CCC$ un préchamp et $F: \CCC \rightarrow \CCC'$ un foncteur à valeurs dans un champ. Si $F$ est une équivalence sur les fibres alors $F^\dagger$ est une équivalence de champs. \\
De même, si $\CCC$ est un champ,  $\CCC'$ est un préchamp et $F: \CCC \rightarrow \CCC'$ un foncteur de préchamps qui induit des équivalences sur les fibres, alors le champ associé au préchamp $\CCC'$ est équivalent au champ $\CCC$.
\end{cor}
\begin{lem}
Les foncteurs entre deux champs donn\'es forment un champ.
\end{lem}

On peut maintenant définir l'image directe et l'image inverse d'un champ. Soit $X$ et $Y$ deux espaces topologiques et $f : X \rightarrow Y$ une fonction continue.
\begin{Def}
(i) Soit $\CCC$ un champ sur $Y$. L'image directe de $\CCC$ par $f$, notée $f_{*}\CCC$ est le champ donné par : 
\begin{itemize}
\item pour tout ouvert $U$ de $Y$, la catégorie $\CCC(f^{-1}(U))$,
\item pour toute inclusion d'ouverts $U \subset V$, le foncteur de restriction $\rho_{f^{-1}(U)f^{-1}(V)}$,
\item pour toute inclusion $U \subset V \subset W$, l'isomorphisme de foncteur $\lambda_{f^{-1}(U)f^{-1}(V)f^{-1}(W)}$, 
\end{itemize}
~~~(ii) Soit $\CCC$ un champ sur $X$. L'image inverse de $\CCC$ par $f$, notée $f^{-1}(\CCC)$ est le champ associé au pré-champ donné par : 
\begin{itemize}
\item pour tout ouvert $U$ de $X$,  la catégorie : 
$$ \displaystyle{2\lim_{\substack{\longrightarrow \\ U \subset V}}}\CCC'(V)$$
\item pour toute inclusion de deux ouverts de $X$, le foncteur de restriction donné par la propriété universelle de la $2$-limite, 
\item pour toute inclusion de trois ouverts, l'isomorphisme donné par la propriété universelle.
\end{itemize}
\end{Def}
\begin{prop}
Les $2$-foncteurs 
$$f_{*} : \SSS t_{X} \longrightarrow \SSS t_{Y} ~~~~~~f^{-1} : \SSS t_{Y} \longrightarrow \SSS t_{X}$$
sont $2$-adjoints, $f_{*}$ étant le $2$-adjoint à droite de $f^{-1}$.
\end{prop}
De plus si $g : Y \rightarrow Z$ est une autre application continue on a  les équivalences naturelles de $2$-foncteurs : 
$$g_{*}\circ f_{*}\simeq (g\circ f)_{*}~~~~~~ f^{-1} \circ g^{-1}\simeq(g \circ f)^{-1}$$
\begin{lem}
Si $f$ est une injection, alors le foncteur naturel $\varepsilon$ d'adjonction :
$$\varepsilon  : f^{-1}f_{*} \longrightarrow Id$$
est une équivalence. 
\end{lem}
Comme pour les faisceaux, on peut recoller les champs sur un recouvrement d'ouverts.
Soit $\{U_{i}\}_{i\in I}$ un recouvrement d'ouverts de $X$. On note $(\CCC_{i}, F_{ji}, \theta_{kji})$ les données suivantes : 
\begin{itemize}
\item pour tout $U_{i}$ un champ $\CCC_{i}$,
\item pour  tout couple $(j,i)$ une équivalence de champ sur $U_{ij}$ : 
$$ F_{ji} :\CCC_{i}\mid_{U_{ij}} \longmapsto \CCC_{j} \mid_{U_{ij}} $$
\item pour tout triplet $(i,j,k)$ un isomorphisme de foncteurs sur $U_{ijk}$:
$$\theta_{kji} : F_{kj} \circ F_{ji} \buildrel\sim\over\longrightarrow F_{ki}$$
tels que les diagrammes suivants commutent : 
$$\xymatrix{ F_{lk} \circ F_{kj} \circ F_{ji} \ar[r] \ar[d] & F_{lj} \circ F_{ji} \ar[d]\\
 F_{lk} \circ F_{ki} \ar[r] & F_{ki}}$$
\end{itemize}
\begin{lem}\label{recouvrement}
Supposons que la famille $(\CCC_{i}, F_{ji}, \theta_{kji})$ soit donnée, alors il existe un champ $\CCC$ sur $X$, des équivalences de champs $F_{i} : \CCC \mid_{U_{i}} \buildrel\sim\over\rightarrow \CCC_{i}$ et des isomorphismes de foncteurs : $\theta_{ij} : F_{ij} \buildrel\sim\over\rightarrow F_{i}\circ F_{j}^{-1}$ vérifiant :   $\theta_{ik} \circ \theta_{ijk} = \theta_{ij} \circ \theta_{jk}$. De plus, la donnée $(\CCC, F_{i}, \theta_{ij})$ est unique à équivalence de champ près, cette équivalence est unique à isomorphisme près.
\end{lem}
\dem
On peut définir le champ $\CCC$ comme une $2$-limite projective. En pratique, pour un ouvert $U$ de $X$, ses sections sont données par une famille $\big\{\{S_{i}\}_{i\in I}, \{g_{ij}\}_{(i,j)\in I^2} \big\}$ avec $S_{i} \in \CCC_{i}(U \cap U_{i})$ et où les $g_{ij}$ sont des isomorphismes $g_{ij} : S_{i}|_{U \cap U_{ij}} \buildrel\sim\over\rightarrow S_{j}|_{U\cap U_{ij}}$, tels que pour tout triplet $(i,j,k)$ le diagramme suivant commute :
$$\xymatrix{
S_{i}|_{U\cap U_{ijk}} \ar[rr]^{g_{ij}|_{U\cap U_{ijk}}} \ar[rd]_{g_{ik}|_{U\cap U_{ijk}}}&& S_{j}|_{U\cap U_{ijk}} \\
& S_{k}|_{U\cap U_{ijk}} \ar[ru]_{g_{kl}|_{U\cap U_{ijk}}}
}$$
\cqfd
Soit $j : W \hookrightarrow X$ l'inclusion d'un sous-ensemble quelconque $W$ de $X$ muni de la topologie induite; on notera $\CCC\mid_{W}=j^{-1}\CCC$ et $\Gamma(W, \CCC)=\Gamma(W, \CCC\mid_{W})$. On a la proposition suivante :
\begin{prop}\label{paracompact}
Soit $\CCC$ un champ sur $X$ et $Z$ un fermé admettant une base de voisinage paracompact, le foncteur canonique suivant est une équivalence :
$$\begin{displaystyle}
\Psi : 2\varinjlim_{U\supset Z}\Gamma(U, \CCC) \buildrel\sim\over\longrightarrow \Gamma(Z, \CCC)
\end{displaystyle}$$
où $U$ décrit les voisinages ouverts de $W$.
\end{prop} 
\dem
Cette démonstration est une adaptation de la démonstration, dans le cas des faisceaux, donnée dans \cite{KS} et \cite{God}.\\
Rappelons tout d'abord qu'un espace paracompact est un espace  séparé tel que pour tout recouvrement d'ouverts $\{U_{i}\}_{i\in I}$ de $X$ il existe un sous-recouvrement $\{V_{j}\}_{j\in J}$ localement fini. Si $\{U_{i}\}_{i \in I}$ est localement fini alors il existe un recouvrement $\{V_{i}\}_{i\in I}$ tel que $\overline{V}_{i}\subset U_{i}$ pour tout $i$. Un fermé d'un espace paracompact est paracompact. Un espace métrisable est paracompact, de même  qu'un espace localement compact dénombrable à l'infini. \\
 $\Psi$ est pleinement fidèle par la version faisceautique de cette proposition.\\
Montrons que $\Psi$ est essentiellement surjectif. Soit $S$ une section de $\Gamma(Z, \CCC)$. Par définition de la fibre en un point d'un champ, on sait que pour tout $x\in Z$ il existe un voisinage ouvert $U_{x}$ de $x$, une section $S_{U_{x}}$ de $\CCC(U_{x})$ et un isomorphisme $\Phi_{U_{x}}$ :
$$\Phi_{U_{x}} : S_{U_{x}} |_{U_{x}\cap Z} \buildrel\sim\over\longrightarrow S|_{U_{x}\cap Z}$$
Quitte à restreindre les ouverts $U_{x}$, on peut supposer que $\bigcup_{x \in Z}U_{x}$ est un ouvert paracompact et donc qu'il existe sous recouvrement $\{U_{i}\}_{i \in I}$ localement fini, des sections $S_{i} \in \CCC(U_{i})$ et des isomorphismes $\Phi_{i}$ :
$$\Phi_{i} :  S_{i} |_{U_{i}\cap Z} \buildrel\sim\over\longrightarrow S|_{U_{i}\cap Z}$$
 On peut alors trouver une famille $\{V_{i}\}_{i \in I}$ telle que $\overline{V}_{i}\subset U_{i}$. La famille $\{\overline{V}_{i}\}$ est localement fini et $V_{i}\supset Z$. Soit $x \in \bigcup V_{i}$, on note $I(x)$ l'ensemble $I(x)=\{i\in I | x \in \overline{V}_{i}\}$, considérons alors $W$ l'ensemble des $x \in \bigcup V_{i}$  tels que pour tout couple $(i,j)$ de $I(x)^2$ il existe un isomorphisme $\Phi_{ijx}$ défini sur un voisinage ouvert $W_{x}$ de $x$ :
$$\Phi_{ijx}: S_{i}|_{W_{x}}\buildrel\sim\over\longrightarrow S_{j}|_{W_{x}}$$
 tel que pour tout triplet $(i,j,k) \in I(x)^3$ le diagramme suivant commute : 
$$\xymatrix{
S_{i}|_{W_{x}} \ar[rr]^{\Phi_{ijx}} \ar[rd]_{\Phi_{ikx}}&&S_{j}|_{W_{x}}\ar[ld]^{\Phi_{jkx}}\\
&S_{k}|_{W_{x}}
}$$
et tels que si $y \in W_{x}\cap Z$ on ait :
$$(\Phi_{ijx})_{y}=( \Phi_{U_{x}})_{y}$$
L'ensemble $W$ est ouvert puisque si $y \in W_{x}$ alors $y\in W$. \\
L'ensemble $W$ contient $Z$. En effet si $x \in Z$ et $(i,j)\in I^2(x)$ la composée $\Phi_{ijx}=\Phi_{i,x}^{-1}\circ \Phi_{i,x}$ est bien un isomorphisme donc il existe un voisinage ouvert, $W_{x}$ de $x$ et un isomorphisme $\Phi_{ij,W_{x}}: S_{i}|_{W_{x}}\rightarrow S_{j}|_{W_{x}} $  , de plus pour $k\in I(x)$ on a bien $\Phi_{jkx}\circ \Phi_{ijx}=\Phi_{ikx}$, donc quitte à restreindre un nombre fini de fois $W_{x}$ on a bien les relations de commutations demandées.  \\
Construisons maintenant l'ensemble de données qui va permettre de définir une section dans le champ $\CCC$.\\
Soient $(i,j) \in I$, l'isomorphisme $\Phi_{ij}= \Phi_{i}|_{V_{i}\cap V_{j}\cap Z}\circ \Phi_{j}^{-1}|_{V_{i}\cap V_{j}\cap Z}$ appartient à $\Hc om_{V_{ij}\cap Z}(S_{i}|_{Z}, S_{j}|_{Z})$ ainsi d'après la version faisceautique de cette proposition il existe un voisinage  ouvert, $T_{ij}$, de $V_{i}\cap V_{j}\cap Z $, et un isomorphisme $\Phi_{ij} : S_{i}|_{T_{ij}}\rightarrow S_{j}|_{T_{ij}}$. On considère alors le recouvrement ouverts $\bigcup_{(i,j)\in I^2} T_{ij}\cap V_{i}\cap W$ et la donnée suivante :
\begin{itemize}
\item pour tout ouvert $T_{ij}\cap V_{i}\cap W$, la section $S_{i}|_{T_{ij}\cap V_{i}\cap W}$,
\item pour toute intersection $T_{ij}\cap V_{i} \cap T_{kl}\cap V_{k}\cap W$, l'isomorphisme $\Phi_{ij}\circ \Phi_{jl}$,
\end{itemize}
Comme tout se passe dans $W$ les relations de commutations sont vérifiées. Ainsi comme $\CCC$ est un champ, il existe une section $S' \in \CCC(\bigcup T_{ij}\cap V_{i}\cap W)$ tel que sa restriction à $Z$ soit isomorphe à $S$.
\cqfd
\begin{Def}
Soit $\Cc$ une catégorie, on note $\widetilde{\CCC}_{\Cc}$ le préchamp constant qui à tout ouvert $U$ de $X$ associe $\Cc$. On note $\CCC_{\Cc}$ le champ constant,  c'est le champ associé au préchamp  $\widetilde{\CCC}_{\Cc}$. \\
On dit qu'un champ, $\CCC$ est localement constant s'il existe un recouvrement d'ouverts $\{U_{i}\}_{i \in I}$ de $X$, tel que la restriction de $\CCC$ à chaque $U_{i}$ soit constant.
\end{Def}
\begin{prop}
Soit $\Cc$ une catégorie complète (qui admet toutes les petites limites) et $X$ un espace localement connexe. On note $\LLL(\Cc)$ le champ des faisceaux localement constants à valeurs dans $\Cc$. Alors il existe une équivalence naturelle de champs :
$$\CCC_{\Cc} \buildrel\sim\over\longrightarrow \LLL(\Cc).$$
\end{prop}
Si $\Cc$ est une catégorie et $X$ un espace topologique, on note \linebreak$Rep(\pi_{1}(X), \Cc)$ la catégorie dont les objets sont les foncteurs du groupoïde $\Pi_{1}(X)$ à valeurs dans la catégorie $\Cc$ et dont les morphismes sont les transformations naturelles entre les foncteurs. Autrement dit, si l'on se fixe un point base de X et un système de générateurs $\gamma_{1}, \ldots, \gamma_{n}$  de $\pi_{1}(X)$, la catégorie $Rep(\pi_{1}(X), \Cc)$ est équivalente à la catégorie dont  les objets de sont les familles 
$$(\mathbf{A}, m_{1}, \ldots, m_{n})$$
où $\mathbf{A}$ est un objet de $\Cc$ et $(m_{1}, \ldots, m_{n})$ sont des automorphismes de $\mathbf{A}$ vérifiant les même relations que $\gamma_{1}, \ldots, \gamma_{n}$ et dont les morphismes entre de tels objets $(\mathbf{A}, m_{1}, \ldots, m_{n})$ et $(\mathbf{A}', m'_{1}, \ldots, m'_{n})$ sont les morphismes entre $\mathbf{A}$ et $\mathbf{A}'$ qui commutent aux $m_{i}$.
\begin{cor}
Le champ constant est équivalent au champ défini par :
\begin{itemize}
\item pour tout ouvert $U$ la catégorie $Rep(\pi_{1}(U), \Cc)$,
\item pour tout couple d'ouverts $(V\subset U)$ le foncteur oubli :
$$Rep (\pi_{1}(U), \Cc) \longrightarrow Rep (\pi_{1}(V), \Cc),$$
\item pour toute inclusion d'ouverts $W\subset V\subset U$ l'identité.
\end{itemize}
\end{cor}

\newpage
\chapter{Recollements de faisceaux et de faisceaux pervers}
Dans ce chapitre, nous nous intéressons à deux problèmes indé\-pen\-dants de recollement de faisceaux sur un espace stratifié. 

Dans un premier temps, nous décrivons la catégorie des faisceaux sur $X$ à partir de données sur les strates. Ce résultat nous est très utile dans le chapitre 3.

Dans la deuxième partie nous nous intéressons à un problème de reconstruction explicite de faisceaux pervers à partir de données sur la plus petite strate et sur le complémentaire de celle-ci. C'est une variante d'un théorème de MacPherson et Vilonen, démontré dans \cite{McV}, qui s'inspire des techniques utilisées dans \cite{GGM}. Ce résultat n'est pas utilisé dans la suite de la thèse.

\section{Faisceaux sur un espace stratifié}

Le but de ce paragraphe est de donner une description de la catégorie $\Ss h_{X}$  des faisceaux sur un espace stratifié $X$.  Nous démontrons l'équi\-va\-lence entre $\Ss h_{X}$ et une catégorie dont les objets sont donnés par un faisceau sur chaque strate et des données de recollement formées par des morphismes de faisceaux. C'est ce que nous formalisons ici. L'énoncé et la démonstration principale de ce paragraphe sont volontairement formels et pourraient paraître plus lourds que nécessaire. L'intérêt de cette approche est que  nous envisageons une généralisation de cette proposition aux cas des champs et une démonstration assez formelle peut facilement être adaptée à ce langage. \\
Commen\c cons par le cas où la stratification de l'espace topologique $X$ n'est formée que de deux strates. L'une est ouverte on la note $U$ et l'autre est fermée on la note $Z$. On note $i$ et $j$ leur inclusion dans $X$ :
$$\begin{array}{cccc}
j :& U & \hookrightarrow & X\\
i : & Z & \hookrightarrow &X
\end{array}$$
Dans ce cas on a le carré cartésien :
 $$\xymatrix{ 
    \Fc \ar[r] \ar[d] \cartesien & i_{*}i^{-1}\Fc \ar[d]^{i_{*}f} \\
    j_{*}j^{-1}\Fc \ar[r]_{\eta_{Z}} & i_{*}i^{-1}j_{*}j^{-1}\Fc  }$$
donc la donnée d'un faisceau $\Fc$ équivaut à la donnée du triplet :
$$(\Fc|_{Z}, \Fc|_{U}, f)$$
où $f$ est un morphisme de faisceau :
$$f : \Fc|_{Z}\longrightarrow i^{-1}j_{*}\Fc|_{U}$$
Ce morphisme est la donnée de recollement.\\ 

Soit maintenant $X$ un espace topologique muni d'une stratification $\Sigma$. On note  $\Ss h_X$ la catégorie des faisceaux sur $X$, $S_{k}$ la réunion des strates de dimension $k$, $i_{k}$ l'injection de $S_{k}$ dans $X$ :
$$i_{k} : S_{k} \hookrightarrow X$$
 et $\eta_{k}$ le morphisme d'adjonction :
 $$\eta_{k} : i_{k}^{-1}i_{k*} \longrightarrow Id.$$
Dans ce cas il faut coder les conditions de recollement pour tout couple $(S_{k}, S_{l})$ tel que $S_{k}\subset \overline{S}_{l}$ c'est à dire tel que $k<l$.
\begin{Def}
Soit $\Ss_\Sigma$ la catégorie dont
\begin{itemize}
\item[$\bullet$] les objets sont les familles $(\{ \Ff_{k}\}_{k\leq n}, \{f_{kl}\}_{l<k\leq n})$
où 
\begin{itemize}
\item $\Ff_{k}$ est un faisceau sur  $S_{k}$,
\item $f_{lk}$ est, pour tout couple $(k,l)$ tel que $k <l$,  un morphisme $f_{lk} : \Ff_{k} \rightarrow i_{l*}\Ff_{l} \mid_{S_{k}}$ tel que le diagramme suivant commute pour tout triplet $k$, $l$, $m$ vérifiant $k<l<m$ : 
$$\xymatrix{\Ff_{k} \ar[d]_{f_{mk}} \ar[r]^{f_{lk}} & (i_{l*}\Ff_{l} )\mid_{S_{k}} \ar[d]^{(i_{l*}f_{ml})\mid_{S_{k}}}\\
                       (i_{m*}\Ff_{m}) \mid_{S_{k}}\ar[r] & (i_{l*}i_{l}^{-1}i_{m*}\Ff_{m})\mid_{S_{k}}\\
}$$
où le morphisme du bas est le morphisme naturel. \\
\end{itemize}
\item[$\bullet$] les morphismes entre deux objets $(\{\Ff_{k}\}, \{f_{kl}\})$ et $(\{\Gc_{k}\}, \{g_{kl}\})$   sont donnés par, pour tout $S_{k}$, un morphisme de faisceaux \linebreak $\phi_{k} : \Ff_{k} \rightarrow  \Gc_{k}$ tels que pour tous $k<l$ le diagramme suivant commute : 
$$\xymatrix{\Ff_{k} \ar[d]_{f_{lk}} \ar[r]^{\phi_{l}} &  \Gc_{k} \ar[d]^{g_{kl}}\\
                       (i_{l*}\Ff_{l})\mid_{S_{k}}\ar[r]^{i_{l*}\phi_{l} \mid_{S_{k}}} & (i_{l*} \Gc)_{S_{k}}\\
}$$
\end{itemize}
\end{Def}
Notons que dans le cas de deux strates les objets de $\Ss_{\Sigma}$ sont bien les triplets $(\Fc_{Z}, \Fc_{U}, f)$.\\
Définissons un foncteur de la catégorie $\Ss h_X$ dans la catégorie $\Ss_{\Sigma}$.
\begin{Def}
On note $R_{\Sigma}$ le foncteur de $\Ss h_X$ dans $\Ss_{\Sigma}$ défini par : 
$$\begin{array}{cccc}
R_{\Sigma}: & \Ss h_X & \longrightarrow & \Ss_\Sigma\\
     & \Fc & \longmapsto &  (\{ \Fc \mid_{S_{k}}\}_{k\leq n }, \{\eta_{kl} \}_{k<l\leq n})\\
     & \phi : \Fc \rightarrow  \Gc& \longmapsto & ( \{\phi \mid_{S_{k}}\}_{k \leq n})
\end{array}$$
où les morphismes $f_{lk}$ sont les morphismes d'adjonction suivant : 
$$ \eta_{kl}: \Fc \mid_{S_{k}} \longrightarrow (i_{l*}i_{l}^{-1}\Fc)\mid_{S_{k}}$$
\end{Def}
Les morphismes $f_{lk}$ étant issus d'une transformation naturelle, ils vérifient bien les conditions de commutation.\\\\
Définissons maintenant le foncteur quasi-inverse $Q_{\Sigma}$ de la catégorie $\Ss_{\Sigma}$ dans la catégorie $\Ss h_X$. \\
Notons que si $\Fc$ est un faisceau sur $X$ un section $s \in \Fc(V)$ est déterminée de manière unique par une famille $\{s_{i}\}_{i \in I}$ où $s_{i}$ appartient à $\Fc|_{S_{i}}(V\cap S_{i} )$ compatibles avec les morphismes de recollement. Formellement $\Fc$ est  limite projective d'un système :
$$\Fc =\varprojlim_{\Ii} \Fc|_{S_{i}}$$
où le système projectif tient compte des morphismes de recollement. Précisons ce système projectif.
\begin{Def}
Soit $\Ii$ la catégorie 
\begin{itemize}
\item dont les objets sont les singletons $\{k\}$ avec $k\leq n$ et les couples $(k,l)$ avec $k< l\leq n$,
\item dont les morphismes entre deux objets est donné par :
$$\begin{array}{lccc}
Hom(i,i) &=&\{*\} & \text{pour tout objet $i$ de $\Ii$}\\
Hom((k,l),k)&=&\{*\}\\
Hom((k,l),l)&=&\{*\}\\
\end{array}$$
\end{itemize}
\end{Def}
Pour définir un système projectif dans $\Ss h_X$ il reste à définir un foncteur contravariant de $\Ii$ dans $\Ss h_X$. Soit $\Ff=(\{\Ff_{k}\},\{f_{lk}\})$ un objet de $\Ss_{\Sigma}$.\\
À tout objet de $\Ii$ on associe : 
$$\begin{array}{cclc}
\Ff(k,l)&=& i_{k*}i_{k}^{-1}i_{l*}\Ff_{l} & \text{pour $(k,l) \in \Ic$~} \\
\Ff(l)&=&i_{l*}\Ff_{l}& \text{pour $\{k\} \in \Ic$~} 
\end{array}$$
À tout morphisme $(k,l) \rightarrow k$ avec $k<l$ on associe les morphismes $f_{lk}$, et aux morphismes $(k,l) \rightarrow (l)$ on associe les morphismes $\eta_{kl}$ donnés par l'adjonction : 
$$\eta_{kl} : i_{l*}\Ff_{l}\longrightarrow i_{k*}i_{k}^{-1}i_{l*}\Ff_{l}$$
\begin{Def}
On définit l'image de $\Ff$ par $Q_{\Sigma}$ comme la limite projective du système $\Ff(k,l)$ : 
$$\begin{displaystyle}
Q_{\Sigma}(\Ff)= \lim_{\substack{ \longleftarrow \\ a \in \Ic }} \Ff(a)
\end{displaystyle}$$
Soient maintenant $\Ff'=(\{\Ff'_{k}\}, \{f_{lk}'\})$  un deuxième objet de $\Ss_{\Sigma}$ et  \linebreak$\phi=\{\phi_{k}\} : \Ff \rightarrow \Ff'$ un morphisme entre ces objets. Les $\phi_{k}$ forment alors un morphisme entre les systèmes $\Ff(k,l)$ et $\Ff'(k,l)$. La propriété universelle que vérifie la limite projective définit un morphisme de $Q_{\Sigma}(\Ff)$ dans $Q_{\Sigma}(\Ff')$, c'est le morphisme image de $\phi$ par $Q_{\Sigma}$.
\end{Def}
\noindent
\textbf{Remarque.}\\
On peut donner une expression explicite de la limite projective. L'image de  $\Ff$ par $Q_{\Sigma}$ est alors donnée par :
\begin{itemize}
\item pour tout ouvert $U$ de $X$,
$$Q_{\Sigma}(\Ff)(U)= \{(s_{1}, \cdots, s_{n}), s_{i} \in \Ff(U\cap S_{i}) | f_{lk}(s_{k}) =\eta_{kl}(s_{l}) \}$$
\item pour tout inclusion d'ouverts $V \subset U$, le morphisme produit $\rho_{V\cap S_{1}}\times \ldots \times \rho_{U\cap S_{n}}$. 
\end{itemize}
\begin{thm}\label{description-faisceaux}
Les catégories $\Ss h_X$ et $\Ss_\Sigma$ sont équivalentes et les foncteurs $Q_{\Sigma}$ et $R_{\Sigma}$  sont quasi-inverses l'un de l'autre. 
\end{thm}
\dem
Rappelons tout d'abord que les limites projectives finies commutent à la restriction. Ainsi, si $\Ff=(\{\Ff_{k}\}, \{f_{lk}\})$ est un objet de $\Ss h_X$, on a l'isomorphisme canonique : 
$$i_{j}^{-1}\varprojlim \Fc(k,l) \simeq \varprojlim i_{j}^{-1}\Fc(k,l)$$
Et si on note $\pi_{kl}$ la projection de $\varprojlim \Ff(m,n)$ sur $\Ff(k,l)$ :
$$\pi_{kl} : \varprojlim \Ff(m,n) \longrightarrow \Ff(k,l),$$
et $\pi_{j}$ la projection naturelle :
$$\pi_{j} : \varprojlim i_{j}^{-1}\Ff(m,n) \longrightarrow i_{j}^{-1}\Ff(j),$$
on a le diagramme commutatif :
$$\xymatrix{
i_{j}^{-1}\varprojlim \Fc(k,l) \ar[rr]^\sim \ar[dr]_{i_{j}^{-1}\pi_{jj}} && \varprojlim i_{j}^{-1}\Fc(k,l) \ar[dl]^{\pi_{j}}\\
&i_{j}^{-1}\Fc(j)
}$$

Pour démontrer le théorème \ref{description-faisceaux} nous définissons deux isomorphismes de foncteurs : 
$$\begin{array}{cccc}
R_{\Sigma}\circ Q_{\Sigma} &\longrightarrow& Id_{\Ss h_{X}} \\
 Id_{\Ss_{\Sigma}} & \longrightarrow & Q_{\Sigma} \circ R_{\Sigma}
 \end{array}$$
Chacun de ces isomorphismes est défini par des morphismes donnés par la propriété universelle de la limite projective. Pour démontrer que ces morphismes sont bien des isomorphismes, nous nous appuyons essentiellement sur le fait  la composée suivante :
$$\psi_{j} :\varprojlim i_{j}^{-1}\Fc(l,k) \buildrel\pi_{j}\over\longrightarrow i_{j}^{-1}i_{j*}\Fc_{j} \buildrel\varepsilon_{j}\over\longrightarrow \Fc_{j}$$
est un isomoprhisme. \\
Définissons l'inverse du morphisme $\psi_{j}$. Par définition de la limite projective il faut se donner une famille de morphismes de $\Ff_{j}$ dans $i_{j}^{-1}\Ff(k,l)$. Mais notons  que :
$$\begin{array}{cl} 
i_{j}^{-1}\Ff (k)=0  & \text{pour ~} k<j\\
i_{j}^{-1}\Ff(k,l)=0 & \text{pour ~}k<l\leq j
\end{array}$$
et que comme $i_{j}$ est une injection le morphisme $\varepsilon_{j} : i_{j}^{-1}i_{j*} \longrightarrow Id$ \linebreak est un isomorphisme. Considérons alors la famille de morphismes  : 
$$\left\{
\begin{array}{lccl}
\phi_{j}^{\{j\}}=&(\varepsilon_{j})^{-1}:& \Ff_{j} \rightarrow i_{j}^{-1}i_{j*}\Ff_{j}=i_{j}^{-1}\Ff(j,j)\\
~\\
\phi_{j}^{\{k\}}=&f_{kj}:& \Ff_{j} \rightarrow i_{j}^{-1}i_{k*}\Ff_{k}=i_{j}^{-1}\Ff(k) &\text{pour~} j<k\\
~\\
\phi_{j}^{(k,l)}=&i_{j}^{-1}\eta_{kl}\circ f_{lj} :& \Ff_{j} \rightarrow i_{j}^{-1}i_{l*}\Ff_{l} \rightarrow  i_{j}^{-1}i_{k*}i_{k}^{-1}i_{l*}\Ff_{l}  &\text{pour ~} j\leq k<l
\end{array}
\right.$$
\begin{lem}\label{rest}
Cette famille définit un unique morphisme, noté $\phi_{j}$
$$
\phi_{j } : \Ff_{j} \longrightarrow \varprojlim_{a \in \Ic} i_{j}^{-1} \Ff(a)
$$
tel que $\pi_{kl} \circ \phi_{j}= \phi_{j}^{(k,l)}$.\\
 Ce morphisme est un isomorphisme et son inverse est la composée : 
$$
\psi_{j} :  (\varprojlim i_{j}^{-1} \Ff(k,l))\buildrel\pi_j\over\longrightarrow i_{j}^{-1}i_{j*}\Ff_{j} \buildrel\varepsilon_{j}\over\longrightarrow \Ff_{j}
$$ 

\end{lem}
\dem
Il s'agit de vérifier la compatibilité entre les $\phi_{j}^{(k,l)}$ et le système projectif. Pour les flèches de $(k,l) \rightarrow (l)$ cela est vrai par construction. Pour les flèches $(k,l) \rightarrow k$, telles que $k>j$, cela résulte directement de la définition d'un objet de $\Ss h_{\Sigma}$. Reste à vérifier la compatibilité pour les flèches $(j,l) \rightarrow \{j\}$. Le diagramme suivant commute, car $\varepsilon_{j}$ est une transformation naturelle :
$$\shorthandoff{;:!?}
\xymatrix @!0 @C=3cm @R=2cm {
\Ff_{j} \ar[d]^{f_{lj}} & i_{j}^{-1}i_{j*}\Ff_{j}\ar[d]^{i_{j}^{-1}i_{j*}f_{lj}} \ar[l]_{\varepsilon_{j}}\\
i_{j}^{-1}i_{l*}\Ff_{l} & i_{j}^{-1}i_{j*}i_{j}^{-1}i_{l*}\Ff_{l} \ar[l]_{\varepsilon_{j}}}$$
Donc $\varepsilon_{j} \circ i_{j}^{-1}i_{j*}f_{lj}= f_{lj} \circ \varepsilon_{j}$, en composant par $i_{j}^{-1}\eta_{j}$ à gauche et $(\varepsilon_{j})^{-1}$ à droite on trouve 
$$i_{j}^{-1}\eta_{j} \circ f_{lj}= i_{j}^{-1}\eta_{j}\varepsilon_{j}(i_{j}^{-1}i_{j*}f_{lj})\varepsilon_{j}^{-1}=(i_{j}^{-1}i_{j*}f_{lj})\varepsilon_{j}^{-1}$$
car $i_{j}^{-1} \eta_{j \varepsilon_{j}}=Id$. Ce qui démontre la compatibilité et on obtient $\phi_{j}$ par la propriété universelle de la limite projective.\\

Montrons que $\psi_{j}$ est l'inverse de $\phi_{j}$.\\
On a, par construction, $\pi_j \circ \phi_{j}= (\varepsilon_{j})^{-1}$ ; donc $\psi_{j} \circ \phi_{j}=Id$.\\
Reste à voir que $\phi_{j} \circ \psi_{j} =Id$. Par propriété universelle cela revient à démontrer que $\pi_{kl} \circ \phi_{j}\circ \psi_{j}= \pi_{kl}$. \\
Pour $\{j\} \in \Ic$, on a, par construction, $\pi_j\circ \phi_{j}=\varepsilon_{j}^{-1}$ ; donc :
$$\pi_j\circ \phi_{j} \circ \psi_{j}=(\varepsilon_{j})^{-1} \circ \varepsilon_{j} \circ  \pi_j =  \pi_j$$\\
Pour $j\leq k \leq l$ et $l>j$, on a le diagramme commutatif suivant : 
$$\xymatrix{
\varprojlim i_{j}^{-1} \Ff(k,l) \ar[rr]^{\psi_{j} \circ \phi_{j}} \ar[dr]_{\psi_{j}} &&   \varprojlim i_{j}^{-1} \Ff(k,l) \ar[dd]^{\pi_{ll}} \ar[ddr]^{\pi_{kl}} \\
& \Ff_{j}\ar[ru]^{\phi_{j}} \ar[rd]_{f_{lj}}& \\
&& i_{j}^{-1}\Ff(l) \ar[r]_{i_{j}^{-1}\eta_{kl}} & i_{j}^{-1}\Ff(k,l)
}$$
Le triangle de droite est commutatif par définition des projections d'une limite projective et le triangle du milieu est commutatif par définition de $\phi_{j}$. On a ainsi les égalités :
$$\begin{array}{ccc}
\pi_{kl} \circ \psi_{j} \circ \phi_{j} & =& i_{j}^{-1}\eta_{kl} \circ f_{lj} \circ \psi_{j}\\
\pi_{ll} \circ \psi_{j} \circ \phi_{j} &=& f_{lj} \circ \psi_{j}
\end{array}$$
Or, d'une part, le diagramme suivant est commutatif  :
$$\shorthandoff{;:!?}
\xymatrix @!0 @C=4cm @R=3cm { 
 \varprojlim i_{j}^{-1} \Ff(k,l) \ar[r]^{\pi_j} \ar[d]_{\pi_{ll}} \ar@/^1cm/[rr]^{\psi_{j}} & i_{j}^{-1}i_{j*}\Ff_{j} \ar[d]_{i_{j}^{-1}i_{j*}f_{lj}}\ar[r]^{\varepsilon_{j}} & \Ff_{j} \ar[d]_{f_{lj}}\\
i_{j}^{-1}i_{l*}\Ff_{l} \ar[r]_{i_{j}^{-1}\eta_{j}} & i_{j}^{-1}i_{j*}i_{j}^{-1}i_{l*}\Ff_{l}  \ar[r]_{\varepsilon_{j}} & i_{j}^{-1}i_{l*}\Ff_{l}
}$$
et, d'autre part, $\varepsilon_{j} \circ i_{j}^{-1}\eta_{j}=Id$. On a donc, pour $l>j$, 
$$f_{lj}\circ \psi_{j}= \pi_{ll}.$$
Ainsi d'une part on a :
$$\pi_{ll}\circ \phi_{j}\circ \psi_{j}=f_{lj} \circ \varepsilon_{j} \circ \pi_j=\pi_{ll}$$
et d'autre part en composant par $i_{j}^{-1} \eta_{kl}$ on obtient : 
$$\pi_{kl} \circ \phi_{j} \circ \psi_{j}= \pi_{kl}$$
\cqfd
\noindent
Revenons à la démonstration du théorème \ref{description-faisceaux}. \\
Montons que $R_{\Sigma}\circ Q_{\Sigma}$ est isomorphe à l'identité. On garde les mêmes notations que ci-dessus, ainsi $\Ff=(\{\Ff_{k}\}, \{f_{lk}\})$ est un objet de $\Ss_{\Sigma}$.  Le lemme \ref{rest} définit un isomorphisme $\phi_{j}$ entre  $\Ff_{j}$ et la restriction de $R_{\Sigma} Q_{\Sigma}(\Ff)$ à $S_{j}$ :
$$\phi_{j} : \Ff_{j} \buildrel\sim\over\longrightarrow R_{\Sigma} Q_{\Sigma}(\Ff) \mid_{S_{j}}$$
Il reste donc à identifier les morphismes $i_{j}^{-1}\eta_{k}$ :
$$i_{j}^{-1}\eta_{k} : i_{j}^{-1}R_{\Sigma} Q_{\Sigma}(\Ff) \longrightarrow i_{j}^{-1}i_{k*}i_{k}^{-1}R_{\Sigma} Q_{\Sigma}(\Ff)$$
 aux morphismes $f_{kj}$. Il s'agit donc de montrer que le diagramme suivant commute :
$$\shorthandoff{;:!?}
\xymatrix @!0 @C=5cm @R=2cm { 
i_{j}^{-1}\varprojlim\Ff(k,l) \ar[r]\eq[d] & i_{j}^{-1}i_{k*}i_{k}^{-1}\varprojlim \Ff(k,l)\eq[d] \\
\varprojlim i_{j}^{-1}\Ff(k,l) \ar[d]_{\psi_{j}}\ar[r]_{\varprojlim i_{j}^{-1}\eta_{k}} & \varprojlim i_{j}^{-1}i_{k*}i_{k}^{-1}\Ff(k,l)\ar[d]^{i_{j}^{-1}i_{k*}\psi_{k}}\\
\Ff_{j} \ar[r]_{f_{kj}} & i_{j}^{-1}i_{k*}\Ff_{k}
}$$
Le carré du haut commute par définition des morphismes d'adjonction. 
Et comme $\phi_{j}$ est l'inverse de $\psi_{j}$ il suffit donc de démontrer l'égalité : 
$$i_{j}^{-1}i_{k*}\psi_{k}\circ i_{j}^{-1}\eta_{k} \circ \phi_{j}= f_{kj}$$
Or on a le diagramme commutatif suivant : 
$$\shorthandoff{;:!?}
\xymatrix @!0 @C=4cm @R=2cm { 
 &\varprojlim i_{j}^{-1}\Ff(l,m) \ar[r] \ar[d]^{\pi_{kk}} & \varprojlim i_{j}^{-1}i_{k*}i_{k}^{-1}\Fc(l,m)\ar[d]^{\pi_{k,k}}\ar[ddr]^{i_{j}^{-1}i_{k*}\psi_{k}}\\
& i_{j}^{-1}\Ff(k) \ar[r]_{i_{j}^{-1}\eta_{k}} & i_{j}^{-1}i_{k*}i_{k}^{-1}\Ff(k)\ar[dr]_{i_{j}^{-1}i_{k*}\varepsilon_{k}}\\
\Ff_{j}\ar[ur]_{f_{kj}} \ar[rrr]_{f_{kj}}\ar[uur]^{\phi_{j}} &&& i_{j}^{-1}i_{k*}\Ff_{k}
}$$
En effet les triangles commutent par définition de $\phi_{j}$ et de $\psi_{k}$.  Le carré du haut commute par définition de la limite projective et le carré du bas car  $i_{j}^{-1}i_{k*}\varepsilon_{k}\circ i_{j}^{-1}\eta_{k}=Id$. On a donc bien l'égalité cherchée.\\

Montrons maintenant que $Q_{\Sigma}\circ R_{\Sigma}$ est isomorphe à l'identité. \\
Soit $\Fc$ un faisceau sur $X$. On rappelle que d'après la définition de $R_{\Sigma}$, on a 
$$R_{\Sigma}(\Fc)=(\{i_{k}^{-1}\Fc\},\{ \eta_{kl}\})$$
où $\eta_{kl}$ sont les morphismes : 
$$\eta_{kl} : i_{k}^{-1} \Fc \longrightarrow i_{k}^{-1}i_{l*}i_{l}^{-1}\Fc$$
Notons que, pour tout $(k,l)$, le diagramme suivant est commutatif . 
$$\shorthandoff{;:!?}
\xymatrix @!0 @C=2cm @R=1.5cm {
\Fc \ar@{.>}[dr]^{\phi} \ar[dddr]_{\eta_{l}} \ar[rrrd]^{\eta_{k}} \\
 & Q_{\Sigma}R_{\Sigma}(\Fc) \ar[dd]^{\pi_{l}} \ar[rr]_{\pi_{k}} &&i_{k*}i_{k}^{-1} \Fc \ar[dd]^{i_{k*}i_{k}^{-1}\eta_{l}}\\
 ~\\
  & i_{l*}i_{l}^{-1}\Fc \ar[rr]_{\eta_{kl}} &&  i_{k*}i_{k}^{-1}i_{l*}i_{l}^{-1}\Fc
}$$
Ainsi d'après la propriété universelle, il existe un unique morphisme $\phi$ tel que le diagramme commute. Considérons sa restriction à $S_{j}$. Comme la restriction commute aux limites projectives finies on a en particulier le diagramme commutatif suivant :
$$\shorthandoff{;:!?}
\xymatrix @!0 @C=4.5cm @R=3cm {
i_{j}^{-1}\Fc \ar[d]_{i_{j}^{-1}\phi} \ar[rd]^{i_{j}^{-1}\eta_{j}} & \\ 
i_{j}^{-1}Q_{\Sigma} R_{\Sigma} (\Fc) \ar[r]_{\pi_j} &i_{j}^{-1} i_{j*}i_{j}^{-1}\Fc
}$$
Or $i_{j}^{-1} \eta_{j}$ est un isomorphisme d'inverse $\varepsilon_{j}$ de plus 
$$\pi_j= i_{j}^{-1}\eta_{j} \circ \varepsilon_{j}\circ \pi_j= i_{j}^{-1}\eta_{j} \circ \psi_{j}$$ est aussi un isomorphisme. Ce qui démontre que la restriction de $\phi$ à chacune des strates est un isomorphisme.
\cqfd

\section{Recollements de faisceaux pervers}
Soit $X$ un espace topologique de Thom-Mather muni de la stratification $\Sigma$ (pour les définitions voir \cite{T} ou \cite{Ma}). Chaque strate $\Sigma_{i}$ est munie d'un voisinage tubulaire $T_{i}$, d'une projection $\pi_{i} : T_{i} \rightarrow \Sigma_{i}$ qui est une fibration localement triviale et d'une fonction $\rho_{i}$ mesurant la distance à la strate $\Sigma_{i}$. On définit le link $\mathbf{L}$ d'une strate $\Sigma_{i}$ comme suit :
$$\mathbf{L}= \big\{x \in T_{i}| \rho_{i}(x)=\varepsilon (\pi_{i}(x))\big\}$$
où $\varepsilon : \Sigma_{i} \rightarrow \RR$ est une fonction définie positive assez petite.\\
On note $\Sigma_0$ la strate de dimension $d$ minimal. \\
Soit $\Pc erv_{X}$ la catégorie des faisceaux pervers sur $X$ relativement à la stratification $\Sigma$ et $\Pc erv_{X_0}$ la catégorie des faisceaux pervers sur $X_0=X \backslash \Sigma_0$ relativement à la stratification $\Sigma \backslash \Sigma_0$. Dans ce chapitre nous notons $i$ et $j$ les injections : 
$$\begin{array}{cccc}
i : &X_{0}&\hookrightarrow& X \\
j : & \Sigma_{0} &\hookrightarrow & X
\end{array}$$
On note $\Ll_{\Sigma_{0}}$ la catégorie des faisceaux localement constants sur $\Sigma_{0}$.\\

Nous allons reprendre un résultat démontré par MacPherson et Vilonen dans \cite{McV}. Dans cet article ils répondent à la question suivante : Si $\Fc$ est un faisceau pervers sur $X$ quelles données, en plus de $\Fc|_{X_{0}}$,  sont nécessaires pour pouvoir reconstruire $\Fc$ ? Autrement dit, quelles données doit on ajouter à la catégorie $\Pc erv_{X_{0}}$ pour qu'elle soit équivalente à la catégorie $\Pc erv_{X}$ ? \\
Pour cela, ils introduisent la notion de fermé pervers.
Puis ils définissent une catégorie $\Pc$ dont les objets sont les quadruplets $(\Hc, H, u, v)$ où $\Hc$ est un faisceau pervers sur $X_{0}$, $H$ est un faisceau localement constant sur $\Sigma_{0}$ et $u$ et $v$ sont des morphismes de faisceaux qui vérifient la relation suivante :
$$\xymatrix{
F(\Hc) \ar[rr]^{T_{\Hc}} \ar[rd]_{u} && G(\Hc) \\
& H \ar[ru]_{v}
}$$
où $F$ et $G$ sont des foncteurs  de la catégorie $\Pc erv_{X_{0}}$ à valeurs dans la catégorie $\Ll_{\Sigma_{0}}$ dépendant d'un fermé pervers et $T$ est une transformation naturelle entre eux. Enfin ils démontrent que $\Pc$ est équivalente à la catégorie $\Pc erv_{X}$ mais ils ne définissent par directement de couple de foncteurs quasi-inverse l'un de l'autre. \\

Après avoir rappelé brièvement ces résultats, nous nous proposons de donner une définition différente d'un fermé pervers dans le but définir une catégorie $\Pc'$ similaire à la catégorie $\Pc$ et  deux foncteurs quasi-inverses l'un de l'autre entre ces deux catégories.

\subsection{Rappel des résultats des MacPherson et Vilonen}~\\
Dans un premier temps MacPherson et Vilonen construisent une ca\-té\-go\-rie $\Cc(F,G,T)$ à partir de la donnée de deux catégorie $\Ac$ et $\Bc$ de deux foncteurs, $F$ et $G$, de $\Ac$ dans $\Bc$ et dune transformation naturelle, $T$, de $F$ dans $G$. 
\begin{Def}
Soit $\Cc(F,G,T)$ la catégorie 
\begin{itemize}
\item dont les objets sont les familles $(A,B,m,n)$ où $A \in \Ac$, $B \in \Bc$ et $m$ et $n$ sont des morphismes tels que le diagramme suivant commute :
$$\xymatrix{
FA \ar[rr]^{T_{A}}\ar[rd]_{m} && GA\\
& B \ar[ur]_{n}
}$$
\item et dont les morphismes sont les couples $(a,b) \in \Mm or(\Ac) \times \Mm or(\Bc)$ tels que le diagramme suivant commute : 
$$
\shorthandoff{;:!?}
\xymatrix @!0 @R=3pc @C=5pc {
     FA\ar[rr]^{T_{A}} \ar[rd]_m \ar[dd]_{Fa} &&  GA \ar[dd]^{Ga} \\
    & B \ar[ru]_{n} \ar[dd]_{b} \\
    FA'\ar[rr] |!{[ur];[dr]}\hole \ar[rd]_{m'} && GA' \\
    & B' \ar[ru]_{n'} }$$
\end{itemize}
\end{Def}
Puis ils démontrent la proposition suivante : 
\begin{prop}
Si $\Ac$ et $\Bc$ sont abéliennes, si $F$ est exact à droite et si $G$ est exact à gauche alors la catégorie $\Cc(F,G,T)$ est abélienne.
\end{prop}
Ils utilisent cette construction pour définir la catégorie équivalente à la catégorie des faisceaux pervers. Ils introduisent alors la notion de fermé pervers.  \\
On note $\mathbf{L}$ le link de la strate $\Sigma_{0}$.
\begin{Def}
 Soit $\Kc$ un fermé de $\mathbf{L}$, on note $\Ll$  l'ouvert complé- mentaire dans $\mathbf{L}$. Le fermé $\Kc$ est dit pervers si pour tout faisceau pervers $\Fc \in \Pc erv_{X_{0}}$ on a 
\begin{itemize}
 \item $(R^{k}\pi_{0*}\Fc_{\Kc}) =0~~ \forall k \geq n-d$
 \item  $(R^{k}\pi_{0*}\Fc_{\Ll}) =0~~ \forall k < n-d$
\end{itemize}
\end{Def}

Si $\Kc$ est un fermé pervers notons que les foncteurs suivants, à valeur dans la catégorie des faisceaux localement constants sur $\Sigma_{0}$ :
$$\begin{array}{ccccc}
F:&\Pc erv_{X_{0}} & \longrightarrow & \Lc_{\Sigma_{0}} \\
&\Fc & \longmapsto & (R^{n-d-1}\pi_{0*}\Fc_{\Kc})\\
G:&\Pc erv_{X_{0}} & \longrightarrow & \Lc_{\Sigma_{0}}\\
&\Fc & \longmapsto & (R^{-d}\pi_{0*}\Fc_{\Ll})\\
\end{array}$$
sont respectivement exacts à gauche et à droite.
\begin{Def}
Soit $\Pc$ la catégorie $\Cc(F,G,T)$ définie avec \linebreak$\Ac=\Pc erv_{X_{0}}$, $\Bc=\Lc_{\Sigma_{0}}$, $F$ et $G$ les foncteurs définis ci-dessus et $T$ la transformation naturelle : 
$$T_{\Fc}: (R^{n-d-1}\pi_{0*}\Fc_{\Kc}) \longrightarrow  (R^{n-d}\pi_{0*}\Fc_{\Ll}) $$
\end{Def}
\begin{thm}
La catégorie $\Pc$ est équivalente à la catégorie $\Pc erv_{X}$.
\end{thm}
Pour démontrer ce théorème MacPherson et Vilonen définissent un foncteur 
$$\Pc erv_{X} \longrightarrow \Cc(F,G,T).$$
Ils montrent que c'est une équivalence de catégorie en introduisant une troisième catégorie et deux autres foncteurs. Dans le paragraphe suivant, on se propose de définir une nouvelle catégorie $\Cc(F',G',T')$ équivalente à la catégorie $\Pc erv_{X}$ et de le démontrer en définissant deux foncteurs quasi-inverses l'un de l'autre.

\subsection*{Définition de la catégorie $\Pc'$ et des foncteurs quasi-inverses.  }~\\
Pour définir la catégorie $\Pc'$ nous avons besoin d'une définition différente d'un fermé pervers.
\begin{Def}
 Soit $\Kc$ un fermé de $X$, on note $\Ll$ son ouvert complémentaire. Le fermé $\Kc$ est dit pervers si pour tout faisceau pervers $\Fc \in \Pc erv_{X_{0}}$ on a 
\begin{itemize}
 \item $\forall k< n -d$, $R^k\Gamma_{\Kc}Ri_{*}\Fc=0$
 \item $\forall k \geq n-d$, $R^k\Gamma_{\Ll}Ri_{*}\Fc=0$
\end{itemize}
\end{Def}
\noindent
\noindent
Dans tout ce qui suit on fixe $\Kc$ un fermé pervers et $\Ll$ son complémentaire. 
On peut alors définir la catégorie $\Pc'$.
\begin{Def}
On note $\Pc'$ est la catégorie $\Cc(F',G',T')$ où $F'$ et $G'$ sont les foncteurs :
$$\begin{array}{cccc}
F' :& \Pc erv_{X_{0}} & \longrightarrow & \Ll_{\Sigma_{0}}\\
 &\Fc & \longmapsto & (R^{n-d-1}\Gamma_{\Lc}Ri_{*}\Fc)_{\Sigma_{0}} \\
G' &  \Pc erv_{X_{0}} & \longrightarrow & \Ll_{\Sigma_{0}}\\
 &\Fc & \longmapsto & (R^{n-d}\Gamma_{\Kc}Ri_{*}\Fc)_{\Sigma_{0}} \\
\end{array}$$
et $T'$ est la transformation naturelle donnée, pour $\Fc \in \Pc erv_{X_{0}}$, par le morphisme naturel:
$$(R^{n-d-1}\Gamma_{\Lc}Ri_{*}\Fc)_{\Sigma_{0}} \longrightarrow  (R^{n-d}\Gamma_{\Kc}Ri_{*}\Fc)_{\Sigma_{0}}$$
\end{Def}
\noindent
Les deux lemmes suivants seront très utiles dans la suite : 
\begin{lem}\label{ferme}
 Si $\Kc$ est un fermé pervers alors $\Ll \cap \Sigma_0= \emptyset$. \\
\end{lem}
\dem
Supposons qu'il existe un fermé pervers $\Kc$ tel que $\Ll \cap \Sigma_0\neq \emptyset$. Le complexe $\CC_{\Sigma_0}[-n+d]$, où $\CC_{\Sigma_0}$ est le faisceau  supporté par $\Sigma_0$ et constant sur $\Sigma_{0}$ de fibre $\CC$, est un faisceau pervers. Or si $x$ appartient à $\Ll \cap \Sigma_0$ on a :
$$(R^{n-d}\Gamma_\Ll \CC_{\Sigma_0}[-n+d])_x \simeq (\Gamma_{\Ll \cap \Sigma_0} \CC_{\Sigma_0})_x \simeq \CC$$
d'où la contradiction. 
\cqfd
\begin{lem}\label{ferme2}
Pour tout faisceau pervers $\Fc \in \Pc erv_X$ : 
$$R\Gamma_\Kc \Fc\simeq R^{n-d}\Gamma_\Kc \Fc [-n+d]$$
\end{lem}
\dem
Pour démontrer ce lemme il suffit de montrer que $R^k\Gamma_{\Kc} \Fc$ est nul pour $k$ strictement supérieur à $n-d$.
Considérons le triangle distingué suivant : 
$$R\Gamma_{\Kc} \Fc \longrightarrow \Fc \longrightarrow R\Gamma_{\Ll}\Fc$$
D'après les conditions de perversité $\Fc$ est nul en degrés supérieur à $n-d$, la définition d'un fermé pervers assure que $R\Gamma_{\Ll}\Fc$ est nul en degré supérieur à $n-d-1$, ainsi la suite exacte longue associée au triangle distingué est la suivante : 
$$h^{n-d}(\Fc) \rightarrow 0 \rightarrow R^{n-d+1}\Gamma_{\Kc} \Fc \rightarrow 0 \rightarrow 0 \rightarrow R^{n-d+2}\Gamma_{\Kc}\Fc \rightarrow0 \rightarrow$$
Ce qui montre que l'on a bien : 
$$R\Gamma_{\Kc}\Fc \simeq R^{n-d}\Gamma_{\Kc}\Fc [-n+d]$$
\cqfd
Définissons le foncteur $\nu$ de la catégorie $\Pc erv_{X}$ dans $\Pc'$. Soit $\Fc$ un faisceau pervers sur $X$.
Notons tout d'abord que ce diagramme, où chacune des flèches est un morphisme naturel, est commutatif : 
$$\xymatrix{ R\Gamma_\Ll Ri_*i^{-1} \Fc  \ar[rr]^{+1} && R\Gamma_{\Kc}Ri_*i^{-1} \Fc [+1]\\
            R\Gamma_\Ll \Fc \ar[u]  \ar[rr]_{+1} & &R\Gamma_\Kc \Fc [+1] \ar[u]}$$
De plus, d'après le lemme \ref{ferme} l'intersection $\Ll \cap X_0$  est égale à $\Lc$ donc le morphisme naturel suivant est un isomorphisme : 
$$R^{n-d-1} \Gamma_\Ll \Fc \stackrel{\alpha}{\longrightarrow} R^{n-d-1}\Gamma_\Ll Ri_* i^{-1} \Fc $$

Ainsi, si l'on pose $F=(R^{n-d}\Gamma_\Kc \Fc)_{\Sigma_0}$ et $u$ et $v$ respectivement les restrictions en $\Sigma_0$ des morphismes :
$$
\begin{array}{cccc}
U: & R^{n-d-1}\Gamma_\Ll Ri_* i^{-1} \Fc \stackrel{\alpha^{-1}}{\rightarrow} R^{n-d-1} \Gamma_\Ll \Fc &\longrightarrow & R^n\Gamma_{\Kc} \Fc \\
V:& R^{n-d} \Gamma_{\Kc} \Fc &\longrightarrow & R^{n-d}\Gamma_{\Kc} Ri_* i^{-1} \Fc\\
\end{array}
$$
le quadruplet $(i^{-1}\Fc, F, u ,v)$ est bien un objet de $\Pc'$. De plus chacune des opérations est bien fonctorielle.\\
\begin{Def}
 On note $\nu$ le foncteur de la catégorie $\Pc erv_X$ dans $\Pc'$ qui à un faisceau pervers $\Fc$ associe le quadruplet $(i^{-1}\Fc, F, u, v)$ défini ci-dessus.
\end{Def}
Définissons un foncteur, noté $\mu$, de la catégorie $\Pc'$ dans la catégorie dérivée des faisceaux. \\
Soit $(\Gc, G, u, v)$ un objet de $\Pc'$. Cette donnée va nous permettre, grâce à la première partie,  de définir un faisceau $\Bc$ supporté par $\Kc$ et un morphisme de  de faisceaux $\phi$ :
$$\phi :R^{n-d-1} \Gamma_\Lc Ri_*\Gc \longrightarrow \Bc$$
Ce morphisme nous permettra alors de construire un complexe dont nous vérifierons qu'il est pervers. \\\\
Le triplet $((R^{n-d}\Gamma_{K}Ri_*\Gc)\mid_{X_0},G, \underline{v})$ où  $\underline{v}$ est la composition du morphisme $v$ et du morphisme d'adjonction  :
$$\underline{v} : G \stackrel{v}{\longrightarrow} (R^{n-d}\Gamma_{K}Ri_*\Gc)\mid_{\Sigma_0} \stackrel{adj}{\longrightarrow} (i_*i^{-1}R^{n-d}\Gamma_{K}Ri_*\Gc)\mid_{\Sigma_0}$$
est un objet de la catégorie $\Ss_{\Sigma_{0}}=\Ss_{\Sigma'}$ où $\Sigma'$ est la stratification formée des deux strates  $\Sigma_{0}$ et $X_{0}$.\\
Soit $\Bc$ l'image du triplet $((R^{n-d}\Gamma_{K}Ri_*\Gc)\mid_{X_0},G, \underline{v})$ par $Q_{\Sigma_{0}}$ :
$$\Bc=Q_{\Sigma_{0}}((R^{n-d}\Gamma_{K}Ri_*\Gc)\mid_{X_0},G, \underline{v})$$
Soit $ \delta$ le morphisme :
$$
\begin{array}{cccc}
\delta : & R^{n-d-1}\Gamma_\Ll Ri_*\Gc &\longrightarrow  & R^{n-d}\Gamma_{K}Ri_*\Gc
\end{array}
$$
Le diagramme suivant est commutatif : 
$$\shorthandoff{;:!?}
\xymatrix@!0 @R=2cm @C=5cm {
(R^{n-d-1}\Gamma_\Ll Ri_*\Gc)\mid_{\Sigma_0} \ar[r] \ar[d]_{u} \ar[rd]^{j^{-1}\delta} & (i_*i^{-1}R^{n-d-1}\Gamma_\Ll Ri_*\Gc)\mid_{\Sigma_0} \ar[dr]^{j^{-1}i_*i^{-1}\delta}\\
G \ar[r]_v & (R^{n-d}\Gamma_{K}Ri_*\Gc)\mid_{\Sigma_0} \ar[r] & (i_*i^{-1}R^{n-d}\Gamma_{K}Ri_*\Gc)\mid_{\Sigma_0}
}$$ 
En effet le triangle est commutatif par définition et le petit carré l'est car  l'adjonction est une transformation naturelle. Ainsi \linebreak le couple $(u, i^{-1}\delta)$ est un morphisme de $\Ss_{\Sigma_{0}}$ \linebreak entre les triplets   $((R^{n-d-1}\Gamma_\Ll Ri_*\Gc)\mid_{X_0},(R^{n-d-1}\Gamma_\Ll Ri_*\Gc)\mid_{\Sigma_0}, \eta)$ et \linebreak $((R^{n-d}\Gamma_{K}Ri_*\Gc)\mid_{X_0},G, \underline{v})$, où $\eta$ est le morphisme d'adjonction :
$$\eta : (R^{n-d-1}\Gamma_\Ll Ri_*\Gc)\mid_{\Sigma_0} \rightarrow (i_*i^{-1}R^{n-d-1}\Gamma_\Ll Ri_*\Gc)\mid_{\Sigma_0}$$
Comme le foncteur $Q_{\Sigma_{0}}$ est quasi-inverse de $R_{\Sigma_{0}}$ il existe un unique morphisme $\phi$ de $R^{n-d-1}\Gamma_\Ll Ri_* \Gc$ dans $\Bc$ qui a pour image par $R_{\Sigma_{0}}$ le morphisme $(u, i^{-1}\delta)$.
$$\phi : R^{n-d-1}\Gamma_\Ll Ri_* \Gc  \longrightarrow \Bc $$
A partir de $\phi$ nous allons construire un morphisme de complexes entre $\tau^{\leq n-d-1}R\Gamma_\Ll Ri_* \Gc$ et $\Bc [-n+d+1]$. Considérons le foncteur de troncature $\tau^{\leq i}$. On a le morphisme naturel :
$$ \tau^{\leq n-d-1} R\Gamma_\Ll Ri_* \Gc \rightarrow  R^{n-d-1}\Gamma_\Ll Ri_* \Gc[-n+d+1]$$
On note alors $\Phi$ le morphisme de complexes défini comme la composée de $\phi$ décalée de $-n+d+1$ et du morphisme de complexe naturel : 
$$ \tau^{\leq n-d-1} R\Gamma_\Ll Ri_* \Gc \rightarrow  R^{n-d-1}\Gamma_\Ll Ri_* \Gc[-n+d+1] \stackrel{\phi[-n+d+1]}{\longrightarrow} \Bc[-n+d+1] .$$
On définit alors l'image de $(\Gc, G, u, v)$ par $\mu$ comme le mapping cône du morphisme de complexes $\Phi$.
\noindent
\\\\
Définissons maintenant l'image par $\mu$ d'un morphisme de $\Pc'$.  Soit donc $(\Gc', G',u',v')$ un deuxième objet de $\Pc'$ et $(g,w)$ un morphisme entre $(\Gc, G,u,v)$ et $(\Gc', G',u',v')$, on note $\Hc$ et $\Hc'$ les images de $(\Gc, G,u,v)$ et $(\Gc', G',u',v')$ par $\mu$. La stratégie est la même que pour la construction des objets : on définit tout d'abord un morphisme entre $\Bc$ et $\Bc'$ dont on se sert pour définir un morphisme de complexes. Comme le diagramme suivant est commutatif : 
$$\xymatrix{
G\ar[d]^{w} \ar[r]^{v~~~~~~~~~~} & (R^n\Gamma_{\Kc}Ri_*\Gc)|_{\Sigma_0} \ar[d]^{(R^n\Gamma_{\Kc}Ri_*g)|_{\Sigma_0}} \ar[r] &(i_*i^{-1}R^n\Gamma_{\Kc}Ri_*\Gc)|_{\Sigma_0} \ar[d]^{(i_*i^{-1}R^n\Gamma_{\Kc}Ri_*g)|_{\Sigma_0}}\\
G' \ar[r]^{v'~~~~~~~~~~} & (R^n\Gamma_{\Kc}Ri_*\Gc')|_{\Sigma_0}  \ar[r] &(i_*i^{-1}R^n\Gamma_{\Kc}Ri_*\Gc)|_{\Sigma_0} 
}$$
le premier carré l'étant par définition, le deuxième car l'adjonction est une transformation naturelle, le couple $((R^n\Gamma_{\Kc_0}Ri_*g)|_{\Sigma_0}, w)$ est un morphisme de $\Ss_{\Sigma_{0}}$ :
$$((R^{n-d}\Gamma_{K}Ri_*\Gc)\mid_{X_0},G, \underline{v})\rightarrow ((R^{n-d}\Gamma_{K}Ri_*\Gc')\mid_{X_0},G', \underline{v'}) $$ 
On note alors $b$ le morphisme image par $Q_{\Sigma_{0}}$. Or par définition d'un morphisme de $\Ss_{\Sigma_{0}}$ le diagramme suivant commute, il suffit de le vérifier sur les restrictions à $X_{0}$ et $\Sigma_{0}$ : 
$$\shorthandoff{;:!?}
\xymatrix@!0 @R=2cm @C=4cm {
R^{n-d-1}\Gamma_\Ll Ri_* \Gc \ar[d]_{R^{n-d-1}\Gamma_\Ll Ri_*g} \ar[r]^{~~~~~\phi}  & \Bc  \ar[d]^{b} \\
  R^{n-d-1}\Gamma_\Ll Ri_* \Gc' \ar[r]^{~~~~~\phi'}& {\Bc'} \\
}$$
Ainsi par fonctorialité de la troncature ce diagramme commute :
$$\xymatrix{
\tau^{\leq n-d-1}R\Gamma_\Ll Ri_* \Gc \ar[d]_{\tau^{\leq n-d-1}R\Gamma_{\Lc} Ri_*g} \ar[r]^\Phi  & \Bc [-n+d+1] \ar[d]^{b[-n+d+1]} \\
 \tau^{\leq n-d-1} R\Gamma_\Ll Ri_* \Gc' \ar[r]^{\Phi'}& {\Bc'}[-n+d+1]  \\
}$$
Ainsi par définition d'une catégorie triangulé, il existe un morphisme $\tilde{g} : \Hc \rightarrow \Hc'$ tel que $(\tilde{g}, R\Gamma_{\Ll}Ri_{*}g, b)$ soit un morphisme de triangle. Mais pour que cette opération soit fonctorielle il faut montrer que ce morphisme est unique. On utilise alors ce lemme. 
\begin{lem}\label{homo}
Soit $K$ une catégorie additive triangulée, $T$ son foncteur de translation. Considérons le diagramme commutatif suivant entre triangles : 
$$\xymatrix{
A \ar[r]^u\ar@{.>}[d]^\alpha & B\ar[r]^v \ar[d]^\beta & C\ar[r]^w \ar[d]^\gamma &T(A) \ar@{.>}[d]\\
A' \ar[r]^{u'} &B'\ar[r]^{v'} & C' \ar[r]^{w'} & T(A') 
}$$
D'après les axiomes des catégories triangulées, il existe un morphisme $\alpha$ tel que $(\alpha,\beta, \gamma)$ soit un morphisme de triangles, si on a de plus :
$$Hom_K(B,T^{-1}C')=Hom_K(C,B')=0$$
ce morphisme est unique.
\end{lem}
\dem
Voir par exemple \cite{Mais}.
\cqfd
Or comme le complexe $\tau^{\leq n-d-1}R\Gamma_{\Ll}Ri_{*}\Gc$ est concentré en degré inférieur ou égal à $n-d-1$ et comme $\Bc'[-n+d]$ est concentré en degré $n-d$, on a  
$$Hom(\tau^{\leq n-d-1}R\Gamma_{\Ll}Ri_{*}\Gc, \Bc'[-n+d])=0$$
De plus $R\Gamma_{\Ll}Ri_{*}\Gc$ étant concentré en degré inférieur ou égal à $n-d-1$, on a l'isomorphisme :
$$\tau^{\leq n-d-1}R\Gamma_{\Ll}Ri_{*}\Gc'\simeq R\Gamma_{\Ll}Ri_{*}\Gc'$$
D'où
$$\begin{array}{cl}
 &Hom(\Bc[-n+d+1],\tau^{\leq n-d-1}R\Gamma_{\Ll}Ri_{*}\Gc' )\\
 \simeq& Hom(\Bc[-n+d+1],R\Gamma_{\Ll}Ri_{*}\Gc' )\\
\simeq &Hom(i_{\Ll}^{-1}\Bc[-n+d+1], i^{-1}_{\Ll}Ri_{*}\Gc')
\end{array}$$
Le deuxième isomorphisme est donné par l'adjonction des foncteur $Ri_{\Ll*}$ et $i_{\Ll}^{-1}$. Mais comme $\Bc$ est supporté par $\Kc$ on a bien :
$$Hom(\Bc[-n+d+1],\tau^{\leq n-d-1}R\Gamma_{\Ll}Ri_{*}\Gc' )=0$$
\begin{Def}
On note $\mu$ le foncteur de la catégorie $\Pc'$ à valeurs dans la catégorie dérivée des faisceaux qui 
\begin{itemize}
\item à un objet $(\Gc, G, u, v)$ de $\Pc'$ associe le mapping cône de $\Phi$ défini plus haut, 
\item et à un morphisme $(g,w) : (\Gc, G, u, v) \rightarrow (\Gc', G', u', v') $ associe le morphisme défini plus haut.
\end{itemize}
\end{Def}

\begin{lem}
 Le cône de $\Phi$ est un faisceau pervers.
\end{lem}
\dem
On note $\Hc$ le mapping cône du morphisme $\Phi$. Le triangle suivant est bien sûr dinstingué : 
$$R\Gamma_\Ll Ri_*\Gc \buildrel\Phi\over\rightarrow \Bc[-n+d-1] \rightarrow  \Hc \buildrel+1\over\rightarrow$$
La démonstration s'appuie sur la proposition suivante démontrer dans \cite{BBD} : 
\begin{prop}
Soit $Y$ une sous-variété fermé de $X$ de dimension $d$, $U=X\backslash Y$. Désignons par $j: Y\hookrightarrow X$ et $i :U \hookrightarrow X$ les inclusions. Soit $\Fc$ un complexe à cohomologie constructible. Alors $\Fc$ est un faisceau pervers si et seulement si 
\begin{itemize}
\item $i^{-1}\Fc$ est un faisceau pervers sur $U$,
\item$ j^{-1}\Fc$ est concentré en degré inférieur ou égal à $n-d$ et $(R\Gamma_{Y}\Fc)|_{Y}$ est concentré en degré supérieur ou égal à $n-d$.
\end{itemize}
\end{prop}
\dem
\cite{BBD}
\cqfd
\begin{lem}
Le complexe $i^{-1}\Hc$ est un faisceau pervers.
\end{lem}
\dem
Le complexe $i^{-1} \Hc$ est le mapping cône du morphisme $i^{-1} \Phi$ : 
$$i^{-1} \Phi: i^{-1} R\Gamma_{\Ll} Ri_{*} \Gc \longrightarrow i^{-1}\Bc[-n+d-1]$$
Or comme $Q_{\Sigma_{0}}$ et $R_{\Sigma_{0}}$ sont quasi-inverse, le faisceau $i^{-1}\Bc$ est isomorphe au faisceau $i^{-1}R^{n-d}\Gamma_\Kc Ri_{*}\Gc$ et $i^{-1} \phi$ est isomorphe à $i^{-1}\delta $  :
$$\delta :  R^{n-d-1}\Gamma_\Ll Ri_*\Gc \longrightarrow   R^{n-d}\Gamma_{K}Ri_*\Gc$$
de plus le diagramme suivant commute : 
$$\shorthandoff{;:!?}
\xymatrix@!0 @R=1cm @C=5cm{
i^{-1} R^{n-d-1}\Gamma_{\Ll} Ri_{*} \Gc \ar[r]^{i^{-1}\phi} \eq[d]&i^{-1}\Bc \eq[d] \\ 
i^{-1} R^{n-d-1}\Gamma_{\Ll} Ri_{*} \Gc \ar[r]^{i^{-1}\delta} & i^{-1} R^{n-d}\Gamma_{\Kc} Ri_{*} \Gc
}$$
Considérons le triangle distingué : 
$$\xymatrix{
\tau^{\leq n-d-1} R\Gamma_\Ll Ri_* \Gc \ar[r] &  R\Gamma_\Ll Ri_* \Gc \ar[r] &\tau^{\geq n-d} R \Gamma_\Ll Ri_* \Gc 
}$$
Par définition d'un fermé pervers, le complexe  $R\Gamma_\Ll Ri_* \Gc $ est nul en degré supérieur ou égal à $n-d$. Ainsi le complexe $\tau^{\geq n-d} R \Gamma_\Ll Ri_* \Gc $ est nul. Le morphisme naturel suivant est donc un isomorphisme : 
$$\tau^{\leq n-d-1} R\Gamma_\Ll Ri_* \Gc \buildrel\sim\over\longrightarrow R\Gamma_\Ll Ri_* \Gc $$
De plus d'après le lemme \ref{ferme2} le faisceau $R^{n-d}\Gamma_{\Kc} Ri_{*} \Gc$ n'est concentré qu'en degré $n-d$ on a ainsi le diagramme commutatif suivant :
$$\shorthandoff{;:!?}
\xymatrix@!0 @R=1cm @C=5cm{
i^{-1}R\Gamma_{\Ll}Ri_{*} \Gc \eq[d] \ar[r]^{i^{-1}\Phi} & i^{-1}\Bc[-n+d+1] \eq[d]\\
i^{-1}R\Gamma_{\Ll}Ri_{*} \Gc  \ar[r]^{\Delta} & i^{-1}R\Gamma_{\Kc}Ri_{*}\Gc [+1]
}$$
où $\Delta$ est le morphisme (+1) du triangle distingué : 
$$i^{-1}R\Gamma_{\Kc}Ri_{*}\Gc\longrightarrow i^{-1}Ri_{*}\Gc \longrightarrow i^{-1}R\Gamma_{\Kc}Ri_{*}\Gc$$
Ainsi $i^{-1} \Hc$ est isomorphe au cône du morphisme $\Delta$. Donc $i^{-1}\Hc$ est isomorphe à $i^{-1}Ri_{*}\Gc$, lui même isomorphe à $\Gc$. Nous montrerons par la suite que cet isomorphisme  est en fait unique. Le complexe $\Gc$ étant pervers par définition, $i^{-1}\Hc$ l'est aussi.
\cqfd
\begin{lem}
Le complexe $R\Gamma_{\Sigma_{0}}\Hc$ est concentré en degré strictement supérieur à $n-d$.
\end{lem}
\dem
Le complexe $R\Gamma_{\Sigma_{0}}\Hc$ est le mapping cône du morphisme : 
$$R\Gamma_{\Sigma_{0}} \Phi : R\Gamma_{\Sigma_{0}}R\Gamma_{\Ll} Ri_*\Fc \longrightarrow R\Gamma_{\Sigma_{0}}\Bc [-n+d+1]$$
Comme, d'après le lemme \ref{ferme}, l'intersection $\Sigma_{0} \cap \Ll  $ est vide, le complexe $R\Gamma_{\Sigma_{0}} \Hc$ est isomorphe à $R\Gamma_{\Sigma_{0}} \Bc[-n+d+1]$. Mais comme le complexe $\Bc[-n+d+1]$ est concentré en degré $n-d$, le complexe $R\Gamma_{\Sigma_{0}} \Bc[-n+d+1]$ est nul en degré inférieur à $n-d$.
\cqfd
\cqfd

\begin{thm}
Les foncteurs $\nu$ et $\mu$ sont quasi-inverses l'un de l'autre. 
\end{thm}
\dem 
Montrons tout d'abord que $\nu \circ \mu$ est isomorphe à l'identité. Soit $(\Gc, G, u, v)$ un objet de $\Pc'$. On note $\Hc$ l'image de  $(\Gc, G, u, v)$ par $\mu$. 
$$\Hc=\mu((\Gc, G, u, v))$$
On a déjà vu que $i^{-1}\Hc$ est isomorphe à $\Gc$, mais il faut vérifier que cet isomorphisme est fonctoriel.\\
C'est le diagramme commutatif suivant :
$$\shorthandoff{;:!?}
\xymatrix@!0 @R=1cm @C=5cm{
i^{-1}R\Gamma_{\Ll}Ri_{*} \Gc \eq[d] \ar[r]^{i^{-1}\Phi} & i^{-1}\Bc[-n+d+1] \eq[d]\\
i^{-1}R\Gamma_{\Ll}Ri_{*} \Gc  \ar[r]^{\Delta} & i^{-1}R\Gamma_{\Kc}Ri_{*}\Gc 
}$$
qui nous a permis de définir, grâce aux axiomes des catégories triangulées, cet isomorphisme. Les deux isomorphismes de ce diagramme étant fonctoriels, si l'on démontre que le lemme \ref{homo} s'applique, on démontre aussi que l'isomorphisme entre $i^{-1}\Hc$ et $\Gc$ est fonctoriel.\\
Notons $i_\Ll$  l'injection de $\Ll$ dans $X$. On a  :
$$\begin{array}{ccl}
Hom(i^{-1}R \Gamma_\Ll Ri_*\Gc, i^{-1} R \Gamma_\Kc \Gc)& \simeq& 0
\end{array}$$
Car $i^{-1}R \Gamma_\Ll Ri_*\Gc$ est concentré en degré $n-d-1$ et $R \Gamma_\Kc \Gc$ est concentré en degré $n-d$.\\
D'autre part, on a les isomorphismes suivants :
$$\begin{array}{ccl}
Hom(i^{-1}\Bc[n-d], i^{-1} R \Gamma_\Ll \Gc)& \simeq& Hom(\Bc[n-d],  R \Gamma_\Ll \Gc)\\
&\simeq &  Hom(i_{\Ll}^{-1}\Bc[n-d], i_{\Ll}^{-1}\Gc)\\
&=&0
\end{array}$$
Le deuxième isomorphisme est donné par l'adjonction et l'égalité est dut au fait que $ i_{\Ll}^{-1}\Bc$ est nul car $ \Bc$ est supporté par $\Kc$.\\
Donc d'après le lemme \ref{homo}, l'isomorphisme entre $i^{-1} \Hc$ et $\Gc$ est unique. Ces deux complexes sont donc fonctoriellement isomorphes. \\
Montrons que $(R\Gamma_\Kc \Gc)|_{\Sigma_0}$ est fonctoriellement isomorphe à $G$. Le faisceau $R\Gamma_{\Kc} \Gc$ est le cône du morphisme :
$$R\Gamma_\Kc \Phi : R\Gamma_\Kc R\Gamma_\Ll Ri_*\Gc \longrightarrow R\Gamma_\Kc \Bc[-n+d+1]$$ 
Mais comme d'une part $\Ll \cap \Kc = \emptyset$ et d'autre part $\Bc$ est supporté par $\Kc$, le complexe  $R\Gamma_{\Kc} \Gc$ est naturellement isomorphe à $\Bc[-n+d+1]$. Ainsi, on a les isomorphismes fonctoriels (le premier étant donné par le lemme \ref{ferme2}): 
$$(R\Gamma_{\Kc} \Gc)|_{\Sigma_0}\simeq (R^n\Gamma_{\Kc} \Gc)|_{\Sigma_0} [-n+d+1] \simeq \Bc |_{\Sigma_0}[-n+d+1]$$
Mais les foncteurs $R_{\Sigma_{0}}$ et $Q_{\Sigma_{0}}$ étant quasi-inverses l'un de l'autre le faisceau $\Bc|_{\Sigma_0}$ est fonctoriellement isomorphe à $G$.
\\ 
\\Montrons maintenant que $\mu \circ \nu$ est isomorphe à l'identité. Soit donc $\Fc$ un faisceau pervers sur $X$. Par définition des foncteurs $\mu$ et $\nu$, $\Bc$ est l'image par $Q_{\Sigma_{0}}$ du triplet $((R^{n-d} \Gamma_{\Kc_0} Ri_*i^{-1}\Fc)|_{X_0}, (R^{n-d}\Gamma_{\Kc} \Fc)|_{\Sigma_0}, \underline{v})$ :
$$\Bc=Q_{\Sigma_{0}}(((R^{n-d} \Gamma_{\Kc_0} Ri_*i^{-1}\Fc)|_{X_0}, (R\Gamma_{\Kc} \Fc)|_{\Sigma_0}, \underline{v}))$$
où $\underline{v}$ est la composition des morphismes : 
$$ j^{-1} R^{n-d}\Gamma_{\Kc} \Fc \rightarrow j^{-1} R^{n-d} \Gamma_{\Kc}Ri_*i^{-1}\Fc \rightarrow j^{-1}i_*i^{-1}R^{n-d} \Gamma_{\Kc}Ri_*i^{-1}\Fc$$
Le faisceau $\Bc$ est naturellement isomorphe au faisceau $R^{n-d}\Gamma_{\Kc}\Fc$. Pour le montrer construisons un isomorphisme dans $\Ss_{\Sigma_{0}}$ entre  \linebreak le triplet $((R^{n-d} \Gamma_{\Kc} Ri_*i^{-1}\Fc)|_{X_0}, (R^{n-d}\Gamma_{\Kc} \Fc)|_{\Sigma_0}, \underline{v})$ et le triplet \linebreak $((R^{n-d} \Gamma_{\Kc}\Fc)|_{X_0}, (R^{n-d} \Gamma_{\Kc}\Fc)|_{\Sigma_0}, \eta)$, où $\eta$ est le morphisme d'adjonction.  Notons $\gamma$ le morphisme de faisceaux défini par une autre adjonction : 
$$\gamma : R^{n-d}\Gamma_{\Kc}\Fc \rightarrow R^{n-d}\Gamma_\Kc Ri_*i^{-1}\Fc $$
Le morphisme $\underline{v}$ s'écrit alors : 
$$\underline{v} =j^{-1}\eta \circ j^{-1}\gamma \simeq j^{-1}(\eta \circ \gamma) $$
Comme l'adjonction est une transformation naturelle le diagramme suivant est commutatif :
$$\shorthandoff{;:!?}
\xymatrix@!0 @R=2cm @C=6cm{
R^{n-d}\Gamma_\Kc \Fc \ar[r] \ar[d]_{ \eta} & R^{n-d}\Gamma_\Kc \Fc \ar[d]^{\eta\circ i_{*}i^{-1}\gamma}\\
i_*i^{-1}R^{n-d}\Gamma_\Kc \Fc \ar[r]_{i_*i^{-1}\gamma}  & i_*i^{-1}R^{n-d}\Gamma_{\Kc}Ri_*i^{-1} \Fc 
}$$
En appliquant $j^{-1}$ à ce diagramme, on obtient le diagramme :
$$\shorthandoff{;:!?}
\xymatrix@!0 @R=2cm @C=6cm{
j^{-1}R^{n-d}\Gamma_\Kc \Fc \ar[r]^{Id} \ar[d]_{j^{-1} \eta} & j^{-1}R^{n-d}\Gamma_\Kc \Fc \ar[d]^{\underline{v}}\\
j^{-1}i_*i^{-1}R^{n-d}\Gamma_\Kc \Fc \ar[r]_{j^{-1}i_*i^{-1}\gamma}  & j^{-1}i_*i^{-1}R^{n-d}\Gamma_{\Kc_0}Ri_*i^{-1} \Fc 
}$$
Ce qui montre que le couple $(i^{-1}\gamma, Id)$ est un morphisme de $\Ss_{\Sigma_{0}}$.
Mais le morphisme $i^{-1}\gamma$ est en fait un isomorphisme, ainsi le couple $(i^{-1}\gamma, Id)$ est un isomorphisme dans $\Ss_X$, donc $\Bc^\bullet$ et $R\Gamma_{\Kc}\Fc$ sont naturellement isomorphe. De plus le diagramme suivant commute : 
$$\shorthandoff{;:!?}
\xymatrix@!0 @R=1cm @C=4cm{
 R^{n-d-1}\Gamma_\Lc Ri_*i^{-1} \Fc \ar[r]^{~~~~~\phi} \eq[d] & \Bc^{\bullet} \eq[d] \\
R^{n-d-1}\Gamma_\Lc \Fc \ar[r] & R^{n-d}\Gamma_{\Kc} \Fc
}$$
En effet pour le voir il suffit de montrer que ses restrictions à $X_{0}$ et $\Sigma_{0}$ commutent. On a de plus $$Hom(R\Gamma_\Lc Ri_*i^{-1} \Fc, \Bc)=Hom(\Bc^{\bullet}[+1],R\Gamma_\Lc \Fc)=0$$ donc d'après le lemme \ref{homo} le cône de $\Phi$ est fonctoriellement isomorphe à $\Fc$.
\cqfd

\chapter{Champs sur un espace stratifié}
Dans tout ce chapitre $X$ désigne un espace topologique et $S$ une stratification fixée de $X$.
\section{Description}

Dans ce paragraphe nous allons donner l'équivalent du théorème \ref{description-faisceaux} pour les champs. Ainsi nous définissons une $2$-catégorie dont les objets sont formés de champs sur chacune des strates, d'une série de foncteurs de champs et d'isomorphismes entre ces foncteurs. Nous démontrons alors que cette $2$-catégorie est $2$-équivalente à la $2$-catégorie des champs sur $X$. La démonstration est sensiblement la même que dans le cas des faisceaux mais appliquée au langage des $2$-catégories.\\
On notera $S_{k}$ la réunion des strates de dimension $k$ et $i_{k}$ l'inclusion de $S_{k}$ dans $X$. \\
Dans ce paragraphe si $\CCC_{l}$ est un champ sur l'ensemble $S_{l}$  on note $\eta_{kl}$ le foncteur naturel d'adjonction :
$$\eta_{kl}: i_{l*}\CCC_{l} \longrightarrow i_{k*}i_{k}^{-1}i_{l*}\CCC_{l}$$
On note $\SSSt_{X}$ la 2-catégorie des champs sur $X$.
\begin{Def}
Soit $\mathfrak{S}_\Sigma$ la 2-catégorie dont 
\begin{itemize}
\item[$\bullet$] les objets sont donnés par : 
\begin{itemize}
\item[-] pour chaque $S_{k}$, un champ $\CCC_{k}$ sur $S_{k}$,
\item[-] pour tout couple $(k,l)$ tel que $k<l$, un foncteur de champs $F_{lk} : \CCC_{k} \rightarrow i_{k}^{-1}i_{l*} \CCC_{l}$,
\item[-] pour tout triplet $(k,l,m)$ vérifiant $k<l<m$, un isomorphisme de foncteurs $\theta_{klm}$ : 
$$\shorthandoff{;:!?}
\xymatrix @!0 @C=1.5cm @R=0.6cm {\CCC_{k} \ar[rrr]^{F_{lk}}  \ar[ddddd]_{F_{mk}} &&& i_{k}^{-1}i_{l*}\CCC_{l} \ar[ddddd]^{i_{k}^{-1}i_{l*}F_{ml}} \\
~\\
& & \ar@{=>}[ld]_\sim^{\theta_{klm}}\\
&~\\
\\
i_{k}^{-1}i_{m*}\CCC_{m} \ar[rrr]_{i^{-1}_k\eta_{lm}} &&& i_{k}^{-1}i_{l*}i_{l}^{-1}i_{m*}\CCC_{m}
}$$
où $\eta_{lm}$ est le foncteur naturel d'adjonction :
$$\eta_{lm} : i_{m*}\CCC_{m}  \rightarrow i_{l*}i_{l}^{-1}i_{m*}\CCC_{m}$$
tels que pour $m<p$ le diagramme suivant soit commutatif :
$$\shorthandoff{;:!?}
\xymatrix@!0 @C=4cm @R=1cm {  
                       i_{k}^{-1}i_{l*}i_{l}^{-1}i_{m*} F_{pm} \circ i_{k}^{-1}i_{l*} F_{ml} \circ F_{lk} \ar[rr]^{Id \bullet \theta_{klm}} \eq[d]^{}&  &  i_{k}^{-1}i_{l*}i_{l}^{-1}i_{m*} F_{pm} \circ i^{-1}_k\eta_l \circ F_{mk} \ar[ddd] \\
                       i_{k}^{-1}i_{l*}(i_{l}^{-1}i_{m*}F_{pm} \circ F_{ml}) \circ F_{lk} \ar[dd]_{i^{-1}_{k}i_{l*}\theta_{lmp}\bullet Id_{F_{lk}}}&  & \\
                       &&\\
                       i_{k}^{-1} i_{l*}(i_{l}^{-1}\eta_{mp} \circ F_{pk}) \circ F_{lk} \eq[d]&&i_{k}^{-1}\eta_{l}\circ i_{k}^{-1}i_{m*}F_{pm}\circ F_{mk}\ar[ddd]^{Id \bullet \theta_{kmp}} \\
                       i_{k}^{-1} i_{l*}i_{l}^{-1}\eta_{mp} \circ  i_{k}^{-1} i_{l*}F_{pl} \circ F_{lk}\ar[dd]_{Id \bullet \theta_{klp}}  \\
                       \\
                        i_{k}^{-1} i_{l*} i_{l}^{-1}\eta_{mp} \circ i_{k}^{-1}\eta_{lp} \circ F_{pk} \ar[rr]&  & i^{-1}_k\eta_l \circ i_{k}^{-1}\eta_{mp} \circ F_{pk}
                       }$$                  
\end{itemize}

\item[$\bullet$] les foncteurs entre deux tels objets, $(\{\CCC_{k}\}, \{F_{kl}\}, \{\theta_{jkl}\})$ \linebreak$(\{\CCC'_{k}\}, \{F'_{kl}\}, \{\theta'_{klm}\})$, sont donnés par : 
\begin{itemize}
\item[-] pour tout $S_{k} \in S$, un foncteur de champ $G_{k} : \CCC_{k} \rightarrow \CCC_{k}'$,
\item[-] pour tout couple $(i,j)$ tel que $i<j$, un isomorphisme de foncteur : 
$$g_{lk} : F_{lk}' \circ G_{k}  \buildrel\sim\over\longrightarrow  i_{k}^{-1}i_{l*} G_{l} \circ F_{lk}$$
tels que :
$$\shorthandoff{;:!?}
\xymatrix@!0 @C=8cm @R=1,5cm {  
i_{k}^{-1}i_{l*}F'_{ml}\circ F'_{lk}\circ G_{k} \ar[r]^{Id \bullet g_{lk}} \ar[ddd]_{\theta_{mlk}\bullet Id} & i_{k}^{-1}i_{l*}F_{ml}'\circ i_{k}^{-1}i_{l*}G_{l}\circ F_{lk}\ar[d]^{i_{k}^{-1}i_{l*}g_{ml} \bullet Id}\\
& i_{k}^{-1}i_{l*}i_{l}^{-1}i_{m*}G_{m} \circ i_{k}^{-1}i_{l*}F_{ml}\circ F_{lk} \ar[d]^{Id \bullet \theta_{mlk}}\\
&  i_{k}^{-1}i_{l*}i_{l}^{-1}i_{m*}G_{m} \circ i_{k}^{-1}\eta_{km}\circ F_{mk} \eq[d]\\
i_{k}^{-1}\eta_{km}\circ F'_{mk}\circ G_{k} \ar[r]_{Id \bullet g_{mk}} & i_{k}^{-1}\eta_{lm}\circ i_{k}^{-1}i_{m*}G_{m}\circ F_{mk}
}$$
\end{itemize}
\item[$\bullet$]
les morphismes entre deux tels foncteurs, (\{$G_{k}\}, \{g_{kl}\})$  (\{$G'_{k}\}, \{g'_{kl}\})$ sont les données pour chaque $S_{k}$ d'un morphisme de foncteur de champ $\phi_{k} :G_{k} \rightarrow G_{k}'$ tels que, pour tout $k,l$ le diagramme suivant commute : 
$$\shorthandoff{;:!?}
\xymatrix@!0 @C=6cm @R=2cm { 
F'_{kl}\circ G_{k} \ar[r]^{g_{kl}} \ar[d]_{Id \bullet \phi_{k}} & i_{k}^{-1}i_{l*}G_{l}\circ F_{lk}\ar[d]^{i_{k}^{-1}i_{l*}\phi_{l}\bullet Id}\\
F'_{kl}\circ G'_{k} \ar[r]_{g'_{kl}} & i_{k}^{-1}i_{l*}G_{l}' \circ F_{lk}
}$$
\end{itemize}
\end{Def}
\noindent
Définissons maintenant deux $2$-foncteurs entre les deux catégories $\SSS t_{X}$ et $\SSS_{\Sigma}$ : 
\begin{Def}
On note $R_{\Sigma}$ le $2$-foncteur de $\SSSt_{X}$ dans $\SSS_{\Sigma}$ qui à un champ sur $X$ associe la famille de ses restrictions sur chacune de ces strates, les foncteurs naturels d'adjonction et les isomorphismes d'adjonction  : 
$$\begin{array}{cccc}
R_{\Sigma}: & \SSSt _{X}  & \longrightarrow & \SSS_{\Sigma}\\
     & \CCC & \longmapsto &  (\{ \CCC \mid_{S_{k}}\}_{k\leq n }, \{i^{-1}_k\eta_l \}_{k<l\leq n}, \{\lambda_{mlk}\}_{k<l<m\leq n} )\\
     & G : \CCC \rightarrow \CCC' & \longmapsto & ( \{G \mid_{S_{k}}\}_{k \leq n}, \{g_{lk}\}_{k<l\leq n})\\
     & \phi : G \rightarrow G' & \longmapsto & (\{\phi \mid_{S_{k}}\}_{k\leq n})
\end{array}$$
où les foncteurs $\eta_l$ sont les foncteurs d'adjonction suivants : 
$$ \eta_{l}: \CCC  \longrightarrow i_{l*}i_{l}^{-1}\CCC$$
et les $\lambda_{mlk}$ sont les isomorphismes naturels de foncteurs : 
$$\shorthandoff{;:!?}
\xymatrix @!0 @C=3.2pc @R=2pc {\CCC\mid_{S_{k}} \ar[rrr]  \ar[ddd] &&& (i_{l*}i_{l}^{-1}\CCC)\mid_{S_{k}} \ar[ddd] \\
 && \ar@{=>}[ld]_\sim^{\lambda_{klm}}\\
& & &\\
(i_{m*}i_{m}^{-1}\CCC)\mid_{S_{k}} \ar[rrr]_{i^{-1}_k\eta_{lm}} &&& (i_{l*}i_{l}^{-1}i_{m*}i_{m}^{-1}\CCC)\mid_{S_{k}}
}$$
De même les isomorphismes de foncteurs  $g_{lk}$ sont les isomorphismes donnés par l'adjonction :
$$\shorthandoff{;:!?}
\xymatrix @!0 @C=3,2pc @R=2pc {\CCC\mid_{S_{k}} \ar[rrr]^{G\mid_{S_{k}}}  \ar[ddd]_{i_k^{-1}\eta_l} &&& \CCC'\mid_{S_{k}} \ar[ddd]^{i_k^{-1}\eta'_l} \\
 && \ar@{=>}[ld]_\sim^{g_{lk}}\\
& & &\\
(i_{l*}i_{l}^{-1}\CCC)\mid_{S_{k}} \ar[rrr]_{(i_{l*}i_{l}^{-1}G)\mid_{S_{k}}} &&& (i_{l*}i_{l}^{-1}\CCC')\mid_{S_{k}}
}$$
Le fait que ces données soient issues d'une $2-$adjonction assure qu'elles appartiennent bien à $\SSS$.
\end{Def}
Définissons maintenant le deuxième $2$-foncteur $Q_{\Sigma}$ de la $2$-catégorie $\SSS_{\Sigma}$ dans la catégorie $\SSSt_{X}$. \\
Soit  $ \CCC=(\{ \CCC_{k}\}, \{F_{lk} \}, \{\theta_{mlk}\} )$ un objet de $\SSS_{\Sigma}$. Comme dans le cas des faisceaux l'image de $\CCC$ par $Q_{\Sigma}$ est définie par une $2$-limite projective. Il faut donc se donner un $2$-système projectif.
\begin{Def}
On note $\III$ la catégorie,
\begin{itemize}
\item dont les objets sont les singletons $\{j\}$ avec $j\leq n$, les couples $(j,k)$ avec $j<k\leq n$, les triplets $(j,k,l)$ avec $j<k<l\leq n$.
\item dont les foncteurs sont les donnés suivantes :
$$\begin{array}{lccc}
Hom(i,i)&=& \{Id_{i}\}& \text{pour tout objet $i$ de $\III$} \\
Hom((j,k),j)&=&\{s_{jk}^j\}\\
Hom((j,k),k)&=&\{s_{jk}^k\}\\
Hom((j,k,l,),(j,k))&=&\{s_{jkl}^{jk}\} &\\
Hom((j,k,l),j)&=&\{s_{jkl}^{j}\} &\\
Hom((j,k,l,),(j,l))&=&\{s_{jkl}^{jl}\}\\
\end{array}$$
\end{itemize}
\end{Def}
On définit alors un $2$-foncteur $\cat a : \III \rightarrow \SSS t_{X}$. A tout objet $i \in \III$ on associe les champs : 
$$\begin{array}{cclc}
\cat a(j)&=&i_{j*}\CCC_{j}& \\
\cat a (j,k)&=& i_{j*}i_{j}^{-1}i_{k*}\CCC_{l} & \\
\cat a (j,k,l)&=& i_{j*}i_{j}^{-1}i_{k*}i_{k}^{-1}i_{l*}\CCC_{l} & \\
\end{array}$$
Considérons les morphismes de $\III$ :
\begin{itemize}
\item
Pour $j<k$, au morphisme $s_{jk}^k :(j,k) \rightarrow j$   on associe le foncteur 
$$\cat a(s_{jk}^j) =F_{kj},$$
\item Au morphisme $s_{kl}^l : (k,l) \rightarrow l$ on associe le foncteur 
$$\cat a(s_{jk}^j)=\eta_{jk}$$ 
donné par la $2$-adjonction : 
$$\eta_{jk} : i_{k*}\CCC_{k}\longrightarrow i_{j*}i_{j}^{-1}i_{k*}\CCC_{k}.$$
\item
Pour $j<k<l$, au morphisme $s_{jkl}^{jk} : (j,k,l) \rightarrow (j,k)$ on associe le foncteur 
$$\cat a(s_{jkl}^{jk}) =i_{j*}i_{j}^{-1}i_{k*}F_{kl},$$
\item au morphisme $s_{jkl}^{jl}(j,k,l,) \rightarrow (j,l)$  on associe le foncteur 
$$\cat a (s_{jkl}^{jl})=i_{j*}i_{j}^{-1}\eta_{kl}$$
 où $\eta_{kl}$ est le foncteur naturel d'adjonction :
$$\eta_{kl} : i_{l*}\CCC_{l} \longrightarrow i_{k*}i_{k}^{-1}i_{l*}\CCC_{l},$$
\item au morphisme $s_{jkl}^j : (j,k,l)\rightarrow j$,  on associe le foncteur composé : 
$$\cat a(s_{jkl}^j)=i_{j}^{-1}\eta_{kl}\circ i_{j*}F_{lj} $$
\end{itemize}
Si $a$ est un objet de $\Ic$, le $1$-morphisme $\cat a_{a} : Id_{\cat a(a)} \rightarrow \cat a(Id_{a}$ est l'identité.\\
 Enfin d'après la définition d'un $2$-foncteur, il faut maintenant se donner, pour tout couple, $(s,s')$, de foncteurs composables  un isomorphisme de foncteurs 
$$ \cat a(s\circ s')\simeq \cat a(s)\circ \cat a (s')$$
 \begin{itemize}
 \item 
pour $j<k<l$, pour les foncteurs $s_{jl}^j$ et $s_{jkl}^{j}$ l'identité convient,
\item pour les foncteurs $s_{jk}^j$ et $s_{jkl}^{jk}$ on se donne l'isomorphisme $\theta_{jkl}$.
\end{itemize}
\begin{Def}
On définit l'image de $\CCC$ par $Q_{\Sigma}$ comme la $2$-limite projective du système $\CCC(k,l)$ : 
$$\begin{displaystyle}
Q_{\Sigma}(\CCC)= 2\lim_{\substack{ \longleftarrow }} \cat a(k,l)
\end{displaystyle}$$
Soient maintenant $\CCC'=(\{\CCC'_{k}\}, \{F_{lk}'\}, \{\theta'_{mlk}\})$  un deuxième objet de $\SSS_{\Sigma}$ et $G=(\{G_{k}\}, \{g_{lk}\}) : \CCC \rightarrow \CCC'$ un morphisme entre ces objets. Les données $(\{G_{k}\}, \{g_{lk}\})$ forment alors un foncteur entre les systèmes $\CCC(k,l)$ et $\CCC'(k,l)$. La propriété universelle que vérifie la $2$-limite projective définit un foncteur de $Q_{\Sigma}(\CCC)$ dans $Q_{\Sigma}(\CCC')$, c'est le morphisme image de $G$ par $Q_{\Sigma}$.\\
Si maintenant $(\{G'_{k}\}, \{g'_{lk}\})$ est un deuxième foncteur de $\CCC$ dans $\CCC'$ et que $\phi=\{\phi_{k}\}$ est un morphisme entre eux, l'image de $\phi$ par $Q_{\Sigma}$ est le morphisme défini par la propriété universelle.
\end{Def}
\textbf{Remarques}\\
De manière explicite, si $\CCC=(\{\CCC_{k}\}, \{F_{lk}\}, \{\theta_{klm}\})$ est un objet de $\Ss_{\Sigma}$ et $U$ est un ouvert de $X$, les objets de $Q_{\Sigma_{0}}(\CCC)(U)$ sont les familles $(\{S_{k}\}, \{g_{lk}\})$, où $S_{k} \in \CCC_{k}(U\cap S_{k})$ et $g_{lk}$ est un isomorphisme $g_{lk} : F_{lk}(S_{k}) \buildrel\sim\over\rightarrow \eta_{kl}(S_{l})$ \\
Notons que, comme dans le cas des faisceaux, la condition de commutation des isomorphismes $\theta_{klm}$ n'est pas nécessaire pour définir le $2$-foncteur $Q_{\Sigma}$. Elle interviendra dans la démonstration de la $2$-équivalence.
\begin{thm}\label{description-champ}
Les $2$-catégories $\SSSt_X$ et $\SSS_\Sigma$ sont $2$-équivalentes et les $2$-foncteurs $Q_{\Sigma}$ et $R_{\Sigma}$  sont quasi-$2$-inverses l'un de l'autre. 
\end{thm}
Les $2$-limites projectives finies commutent aux restrictions. Ainsi si $\CCC=(\{\CCC_{k}\}, \{F_{lk}\}, \{\theta_{klm}\})$ est  un objet de $\SSS_{\Sigma}$. On a l'équivalence :
$$i_{j}^{-1}2\varprojlim \CCC(l,k) \stackrel{\sim}{\rightarrow} 2\varprojlim i_{j}^{-1} \CCC(l,k)$$
 On note $\pi_{kl}$ et $p_{kl}$ les projections et les isomorphismes de projections : 
$$\shorthandoff{;:!?}
\xymatrix @!0 @C=4cm @R1.5cm { 
  & \CCC(k,k) \ar[dd]^{F_{lk}}\\
  \varprojlim \CCC(m,n)  \ar[ru]^{\pi_{k}} \ar[rd]_{\pi_{kl}}|*{}="A" \ar@{=>}[ru];"A"^{p_{kl}}\\
  & \CCC(k,l)
     }$$
Comme pour les faisceaux nous définissons deux $2$-équivalences de $2$-foncteurs entre les $2$-foncteurs :
$$\begin{array}{ccc}
R_{\Sigma}Q_{\Sigma} &\rightarrow &Id\\
Id & \rightarrow &Q_{\Sigma}R_{\Sigma}
\end{array}$$
elle sont définies à partir de foncteurs et d'isomorphismes de foncteurs donnés par la propriété universelle de la $2$-limite. Mais avant cela montrons que le foncteur suivant, où le deuxième foncteur est le foncteur d'adjonction, est une équivalence. 
$$\Psi_{j} : 2\varprojlim i_{j}^{-1}\CCC(l,k) \rightarrow i_{j}^{-1}i_{j*}\CCC_{j} \rightarrow \CCC_{j}$$
Pour cela on définit un foncteur quasi-inverse $\Phi_{j} :\CCC_{j}\rightarrow i_{j}^{-1}\varprojlim \CCC(a)$. Comme la restriction commute aux $2$-limites projectives finies, pour tout $a \in \III$,  il faut se donner un foncteur de $\CCC_{j}$ dans $i_{j}^{-1}\CCC(a)$. Mais notons que, pour $j<k<l<m$ :
$$\begin{array}{cl}
i_{j}^{-1}\CCC(k)=0 \\
i_{j}^{-1}\CCC(k,l)=0\\
i_{j}^{-1}\CCC(k,l,m)=0\\
\end{array}$$
De plus comme $i_{j}$ est une injection la transformation naturelle \linebreak$\varepsilon_{j} : i_{j}^{-1}i_{j*}\rightarrow Id$ est une équivalence. On choisit $\varepsilon^{-1}_{j}$ un quasi-inverse de $\varepsilon_{j}$ et on fixe $e_{j}$  l'isomorphisme :
$$\begin{array}{ccccccc}
e_{j}: & \varepsilon_{j} \circ \varepsilon_{j}^{-1} &\longrightarrow& Id\\
\end{array}$$
Considérons alors la famille de foncteurs, pour $j\leq k<l<m$ : 
$$\left\{
\begin{array}{lccl}
\Phi_{j}^{j}=(\varepsilon_{j})^{-1}: \CCC_{j} \rightarrow i_{j}^{-1}i_{j*}\CCC_{j}=i_{j}^{-1}\CCC(j,j)\\
~\\
\Phi_{j}^{k}=F_{kj}: \CCC_{j} \rightarrow i_{j}^{-1}i_{k*}\CCC_{k}=i_{j}^{-1}\CCC(k,k)\\
 ~\\
\Phi_{j}^{(k,l)}=i_{j}^{-1}\eta_{kl}\circ F_{lj} : \CCC_{j} \rightarrow i_{j}^{-1}i_{l*}\CCC_{l} \rightarrow  i_{j}^{-1}i_{k*}i_{k}^{-1}i_{l*}\CCC_{l}  \\
~\\
\Phi_{j}^{(k,l,m)}=(i_{j}^{-1}i_{k*}i_{k}^{-1}\eta_{lm}) \circ (i_{j}^{-1}\eta_{km}) \circ F_{mj}\\
\end{array}
\right.$$
Pour définir un foncteur à valeur dans le $2$-système projectif, il faut se donner, pour chaque foncteur $b \rightarrow a$ de la $2$-catégorie $\III$, un isomorphisme $h_{j}^{ab}$ : 
$$\shorthandoff{;:!?}
\xymatrix @!0 @C=0.7cm @R=0.7cm {
&&&\CCC_{j} \ar[ddddlll]_{\Phi_{j}^{a}} \ar[ddddrrr]^{\Phi_{j}^b}\\
~\\
&&&&~\\
&&\ar@{=>}[urr]_{h_{j}^{ab}}^\sim\\
\CCC(a) \ar[rrrrrr]&&&&&&\CCC(b)\\
}$$
\begin{itemize}
\item considérons le foncteur $(k,l) \rightarrow l$ par définition des foncteurs $\Phi_{j}^{(k,l)}$ et $\Phi_{j}^l$ on peut prendre $h_{j}^k=Id$.\\
\item pour le foncteur $(k,l) \rightarrow k$ supposons tout d'abord que $j<k<l$, d'après la définition de $\Phi_{j}^{(k,l)}$ l'isomorphisme $\theta_{jkl}$ convient. \\
Supposons maintenant que $k=j$, nous devons définir un morphisme entre les foncteurs :
$$i_{j}^{-1}i_{j*}F_{kj} \circ \varepsilon_{j}^{-1} \buildrel\sim\over\longrightarrow i_{j}^{-1}\eta_{jk}\circ F_{kj}$$
Comme $\varepsilon_{j}$ est issue d'une $2$-transformation naturelle on a l'isomorphisme $\theta_{j} : F_{kj}\circ \varepsilon_{j} \buildrel\sim\over\rightarrow \varepsilon_{j} \circ (i_{j}^{-1}i_{j*}F_{kj}) $  : 
$$\shorthandoff{;:!?}\relax
\xymatrix @!0 @C=1.5cm @R=0,6cm {\CCC_{j}  \ar[ddddd]_{F_{kj}} &&& i_{j}^{-1}i_{j*}\CCC_{j} \ar[lll]_{\varepsilon_{j}} \ar[ddddd]^{i_{j}^{-1}i_{j*}F_{kj}} \\
~\\
& & \ar@{=>}[ld]_{\sim}^{\theta_{j}}\\
&~\\
\\
i_{j}^{-1} i_{k*}\CCC_{k} &&& i_{j}^{-1}i_{j*}i_{j}^{-1}\i_{k*}\CCC_{k}\ar[lll]^{\varepsilon_{j}}}$$
Notons que la restriction de la $2$-transformation naturelle $\eta_{j}$ à $S_{j}$ est un quasi-inverse de $\varepsilon_{j}$, on note $n_{j}$ l'isomorphisme : 
$$n_{j} : Id \longrightarrow i_{j}^{-1} \eta_{j} \circ \varepsilon_{j}$$
On pose alors : 
$$h_{j}^{jk} = n_{j} \bullet (I_{i_{j}^{-1}\eta_{j}} \circ \theta_{j} \circ I_{\varepsilon_{j}^{-1}}) \bullet e_{j}$$
où $I_{i_{j}^{-1}\eta_{j}}$ et $I_{\varepsilon_{j}^{-1}}$ sont les morphismes identité entre les $2$-transformations naturelles $i_{j}^{-1}\eta_{j}$ et $\varepsilon_{j}^{-1}$, et où $\bullet$ est la composition verticale des morphismes de foncteurs et $\circ$ est la composition horizontale. \\
\item pour le foncteur $(k,l,m) \rightarrow (k,m)$, au vu des définitions de $\Phi_{j}^{km}$ et $\Phi_{j}^{klm}$, l'identité convient. \\
\item pour le foncteur $(k,l,m) \rightarrow (k,l)$, d'après les définitions de $\Phi_{j}^{(k,l)}$ et $\Phi_{j}^{(k,l,m)}$ nous devons définir un morphisme entre les composés de foncteurs : 
$$(i_{j}^{-1}i_{k*}i_{k}^{-1}i_{l*}F_{ml} )\circ (i_{j}^{-1}\eta_{kl} )\circ F_{lj} \buildrel\sim\over\rightarrow (i_{j}^{-1}i_{k*}i_{k}^{-1}\eta_{lm} )\circ i_{j}^{-1}\eta_{km} \circ F_{mj}$$
Mais par définition de $\theta_{jlm}$ et comme $\eta_{kl}$ provient d'une $2$-adjonction on a le  diagramme suivant : 
$$\shorthandoff{;:!?}
\xymatrix @!0 @C=4cm @R=1cm {
& i_{j}^{-1}i_{l*}\CCC_{l} \ar[r]^{i_{j}^{-1}\eta_{kl}} \ar[dddr]_{i_{j}^{-1}i_{l*}F_{ml}} &  i_{j}^{-1}i_{k*}i_{k}^{-1}i_{l*}\CCC_{l} \ar[dddr]^{i_{j}^{-1}i_{k*}i_{k}^{-1}i_{l*}F_{ml}} \\
&\ar@{=>}[d]^\sim_{\theta_{jlm}} & \ar@{=>}[d]^\sim_{i_{j}^{-1}adj}\\
&~& ~\\
\CCC_{j} \ar[uuur]^{F_{lj}} \ar[r]_{F_{mj}} & i_{j}^{-1}i_{m*}\CCC_{m} \ar[r]_{ i_{j}^{-1}i_{k*}i_{k}^{-1}\eta_{lm}}& i_{j}^{-1}i_{k*}i_{k}^{-1}i_{m*}\CCC_{m} \ar[r] &  i_{j}^{-1}i_{k*}i_{k}^{-1}i_{l*}i_{l}^{-1}i_{m*}\CCC_{m}
}$$
En composant correctement ces deux isomorphismes on définit l'isomorphisme cherché.
\end{itemize}
Les relations de commutations des $\theta_{jkl}$ demandées dans la définition de la $2$-catégorie $\SSS_{\Sigma}$ ainsi que le caractère naturel des foncteurs d'adjonction assurent que les conditions de commutations sont vérifiées. Ainsi la propriété universelle définit un foncteur noté $\Phi_{j}$ :
$$\Phi_{j} : \CCC_{j} \longrightarrow i_{j}^{-1}2\varprojlim \CCC(k,l)$$
et pour tout objet $a$ de $\III$ un ismorphisme, noté $\varphi_{j}^{a} $ : 
$$\varphi_{j}^{a} : \pi_{a} \circ \Phi_{j} \buildrel\sim\over\longrightarrow  \Phi_{j}^{a}$$
\begin{lem}\label{phi}
Le foncteur $\Phi_{j}$ est une équivalence et la composée des foncteurs suivants est un quasi-inverse : 
$$\Psi_{j} : 2\varprojlim_{a \in \III} i_{j}^{-1}\CCC(a) \buildrel\pi_{j}\over\longrightarrow i_{j}^{-1}i_{j*}\CCC_{j} \buildrel\varepsilon_{j}\over\longrightarrow \CCC_{j} $$
\end{lem}
\dem
Commen\c cons par montrer que $\Psi_{j}\circ \Phi_{j}$ est isomorphe à l'identité. On a la suite d'isomorphisme :
$$\shorthandoff{;:!?}
\xymatrix @!0 @C=3.5cm @R=1cm {
\Psi_{j} \circ \Phi_{j} \ar@{=}[r] &\varepsilon_{j}\circ \pi_{j}\circ \Phi_{j} \ar[r]^\sim_{I_{\varepsilon_{j}\circ \varphi^j_{j}}} & \varepsilon_{j}\circ (\varepsilon_{j})^{-1}\ar[r]^{~~~~~~\sim}_{~~~~~~e_{j}} & Id
}$$
L'égalité est  donnée par définition de $\Psi_{j}$. Ainsi l'isomorphisme \linebreak $e_{j}\bullet( I_{\varepsilon_{j}}\circ \varphi_{j}^j )$ convient.\\
Montrons maintenant que $\Phi_{j}\circ \Psi_{j}$ est isomorphe à l'identité. Nous allons pour cela définir pour tout $a \in \III$ un isomorphisme noté $l_{j}^{a}$ :
$$l_{j}^{a} : \pi_{a} \circ \Phi_{j} \circ \Psi_{j} \buildrel\sim\over\longrightarrow \pi_{a}$$
Commen\c cons par définir $l_{j}^m$ pour $m\leq n$.
On a le diagramme suivant :
$$\shorthandoff{;:!?}
\xymatrix @!0 @C=1cm @R=1cm {
2\varprojlim \CCC(a)\ar[rrrdd]_{\Psi_{j}} \ar[rrrrrrrr]^{\Phi_{j}\circ \Psi_{j}}&&&&&&&&2\varprojlim \CCC(a) \ar[ddd]^{\pi_{m}}\\
&&&&&&& \ar@{=>}[dl]_{\sim}^{\varphi_{j}^m}\\
& && \CCC_{j} \ar[rrrrruu]^{\Phi_{j}} \ar[rrrrrd]_{F_{mj}}&&&&&\\
&&&&&&&& i_{j}^{-1}i_{m*}\CCC_{m}
}$$
On peut ainsi définir un isomorphisme entre les foncteurs :
\begin{equation}\label{eq1}
\pi_{m}\circ\Phi_{j}\circ \Psi_{j} \buildrel\sim\over\longrightarrow F_{mj}\circ \Psi_{j}
\end{equation}
Considérons maintenant le diagramme suivant : 
$$\shorthandoff{;:!?}
\xymatrix @!0 @C=2cm @R=1.5cm { 
&&&\ar@{=}[d]\\
 2\varprojlim i_{j}^{-1} \CCC(a) \ar[rrr]^{\pi_{j}} \ar[ddd]_{\pi_{m}} \ar@/^1.5cm/[rrrrrr]^{\Psi_{j}} &&& i_{j}^{-1}i_{j*}\CCC_{j} \ar[ddd]_{i_{j}^{-1}i_{j*}F_{mj}}\ar[rrr]^{\varepsilon_{j}} &&& \CCC_{j} \ar[ddd]_{F_{mj}}\\
&&\ar@{=>}[dl]_{\sim}^{p_{jm}}&&&\ar@{=>}[dl]_{\sim}\\
&&&&&&&\\
i_{j}^{-1}i_{m*}\CCC_{m} \ar[rrr]_{i_{j}^{-1}\eta_{j}} \ar@/_1.5cm/[rrrrrr]_{Id}&&& i_{j}^{-1}i_{j*}i_{j}^{-1}i_{m*}\CCC_{m}  \ar[rrr]_{\varepsilon_{j}} &&& i_{j}^{-1}i_{m*}\CCC_{m}\\
&&& \ar@{=>}[u]_{\sim}
}$$
On rappelle que $p_{jm}$ est l'isomorphisme donné par la $2$-limite. L'égalité est donnée par la définition de $\Psi_{j}$. L'isomorphisme de droite et l'isomorphisme vertical sont définis par la $2$-adjonction. En composant correctement ces isomorphismes on définit un isomorphisme entre les foncteurs $F_{mj} \circ \Psi_{j}$ et $\pi_{m} $ :
\begin{equation}\label{eq2}
F_{mj}\circ \Psi_{j} \buildrel\sim\over\longrightarrow \pi_{m}
\end{equation}
On définit alors $l_{j}^m$ comme la composée verticale des isomorphismes (1) et (2).\\
Il faut maintenant définir les isomorphismes $l_{j}^{lm}$ et $l_{j}^{klm}$. Mais la $2$-limite donne les isomorphismes suivants : 
$$\shorthandoff{;:!?}
\xymatrix @!0 @C=1cm @R=1.5cm {
&&&&&2\varprojlim \CCC(a) \ar[dd]_{\pi_{ml}}\ar[ddlllll]_{\pi_{m}} \ar[ddrrrrr]^{\pi_{klm}}\\
&&&\ar@{=>}[r]^\sim_{p_{lm}}&&&\ar@{=>}[r]^\sim_{p_{klm}}&&\\
i_{j}^{-1}i_{m*}\CCC_{m}  \ar[rrrr]_{~~~~~~i_{j}^{-1}\eta_{lm}} &&&&&~~~~i_{j}^{-1}i_{l*}i_{l}^{-1}i_{m*}\CCC_{m}\ar[rrrr]_{~~~~~~~~~i_{j}^{-1}\eta_{klm}}&&&&&~~~~~~~~~~~~~~i_{j}^{-1}i_{k*}i_{k}^{-1}i_{l*}i_{l}^{-1}i_{m*}\CCC_{m}
}$$
Ainsi en composant horizontalement $p_{lm}$ avec les morphismes identité des foncteurs $\Phi_{j}$ et $\Psi_{j}$ on obtient l'isomorphisme :
\begin{equation}
\pi_{lm} \circ \Phi_{j} \circ \Psi_{j} \buildrel\sim\over\longrightarrow i_{j}^{-1}\eta_{lm} \circ \pi_{m} \circ \Phi_{j} \circ \Psi_{j}
\end{equation}
De même en composant horizontalement le morphisme identité du foncteur $i_{j}^{-1}\eta_{lm}$ et l'isomorphisme $l_{j}^m$ on obtient l'isomorphisme : 
\begin{equation}
 i_{j}^{-1}\eta_{lm} \circ \pi_{m} \circ \Phi_{j} \circ \Psi_{j} \buildrel\sim\over\longrightarrow  i_{j}^{-1}\eta_{lm} \circ \pi_{m}
\end{equation}
On définit alors $l_{j}^{lm}$ comme la composée verticale des isomorphismes (3), (4) et de l'inverse de $p_{lm}$. On fait de même pour $l_{j}^{klm}$.\\

Le caractère naturel des foncteurs d'adjonctions, les compatibilités des isomorphismes de projections ainsi que les conditions de commutation demandées dans la définition de la $2$-catégorie $\SSS_{\Sigma}$ assurent que les isomorphismes $l_{j}^{a}$ sont compatibles. Ils définissent ainsi un  isomorphisme entre les foncteurs $\Phi_{j}\circ \Psi_{j}$ et l'identité.
\cqfd
Revenons à la démonstration de la proposition \ref{description-champ}. \\

Montrons tout d'abord que $R_{\Sigma}\circ Q_{\Sigma}$ est équivalent à l'identité. Soit $\CCC=(\{\CCC_{k}\}, \{F_{lk}\}, \{\theta_{klm}\})$ un objet de $\SSS_{\Sigma}$. On a défini dans le lemme \ref{phi} une équivalence entre $\CCC_{j}$ et la restriction de $R_{\Sigma}Q_{\Sigma}(\CCC)$ à $S_{j}$   :
$$\Phi_{j}: \CCC_{j} \longrightarrow i_{j}^{-1}2\varprojlim_{ a\in \III}\CCC(a) $$
Comme dans le cas des faisceaux il faut maintenant identifier le foncteur naturel :
$$i_{j}^{-1}\eta_{k} : i_{j}^{-1}2\varprojlim \CCC(a) \longrightarrow i_{j}^{-1}i_{k*}i_{k}^{-1}2\varprojlim \CCC(a)$$
avec le foncteur $F_{kj}$. C'est à dire qu'il faut construire les isomorphismes suivants : 
$$\shorthandoff{;:!?}\relax
\xymatrix @!0 @C=1.5cm @R=0,6cm {i_{j}^{-1}2\varprojlim \CCC(a)\ar[rrr]^{i_{j}^{-1}\eta_{k}} \ar[ddddd]_{} &&& i_{j}^{-1}i_{k*}i_{k}^{-1}2\varprojlim\CCC(a) \ar[ddddd]^{} \\
~\\
& & \ar@{=>}[ld]_{\sim}^{}\\
&~\\
\\
2\varprojlim i_{j}^{-1} \CCC(a)\ar[rrr] \ar[ddddd]_{\Psi_{j}} &&&2\varprojlim i_{j}^{-1}i_{k*}i_{k}^{-1}\CCC(a) \ar[ddddd]^{i_{j}^{-1}i_{k*}\Psi_{k}}\\
~\\
& & \ar@{=>}[ld]_{\sim}^{}\\
&~\\
\\
\CCC_{j} \ar[rrr]_{F_{kj}} &&& i_{j}^{-1}i_{k*}\CCC_{k}
}$$
Le premier est acquis car la restriction et l'image directe commute à isomorphisme près aux $2$-limites projectives finies. \\
Définissons le second. Comme $\Phi_{j}$ est quasi-inverse de $\Psi_{j}$ pour définir l'isomorphisme cherché il suffit d'en définir un entre les foncteurs :
$$i_{j}^{-1}i_{k*}\Psi_{j} \circ i_{j}^{-1}\eta_{k}\circ \Phi_{j} \buildrel\sim\over\longrightarrow F_{kj}$$
Or on a les isomorphismes suivants :
$$\shorthandoff{;:!?}
\xymatrix @!0 @C=1.5cm @R=0.8cm { 
 &&&2\varprojlim i_{j}^{-1}\CCC(a)\ar[rrr] \ar[ddddd]_{\pi_{k}} &&&2 \varprojlim i_{j}^{-1}i_{k*}i_{k}^{-1}\CCC(a)\ar[ddddd]^{\pi_{k,}}\ar[dddddddddrrr]^{i_{j}^{-1}i_{k*}\Psi_{k}}\\
& \\
&&&&& \ar@{=>}[ld]_{\sim}^{}\\
&&&&&&&&&&~\\
&&\ar@{=>}[dd]_{\sim}&&&&& \ar@{=>}[dd]^\sim\\
&&& i_{j}^{-1}\CCC(k) \ar[rrr]_{i_{j}^{-1}\eta_{k}} &&& i_{j}^{-1}i_{k*}i_{k}^{-1}\CCC(k)\ar[ddddrrr]_{i_{j}^{-1}i_{k*}\varepsilon_{k}}\\
&&&&&\ar@{=>}[dd]_{\sim}&&&&&&&&\\
&&&&\\
&&&&&&&&&&&&&\\
\CCC_{j}\ar[uuuuuuuuurrr]^{\Phi_{j}} \ar[rrrrrrrrr]_{F_{kj}}\ar[rrruuuu]_{F_{kj}} &&&&&&&&& i_{j}^{-1}i_{k*}\CCC_{k}
}$$
Les deux isomorphismes des triangles sont donnés par la définition de $\Phi_{j}$ et $\Psi_{k}$ et enfin l'isomorphisme du bas est la composée de l'isomorphisme naturel d'adjonction :
$$i_{j}^{-1}\eta_{k} \circ i_{j}^{-1}i_{k*}\varepsilon_{k} \simeq Id$$
et de l'identité du foncteur $F_{kj}$. En composant correctement ces isomorphismes on obtient l'isomorphisme cherché.\\

Montrons maintenant que $Q_{\Sigma}R_{\Sigma}$ est équivalent à l'identité.  Soit $\GGG$ un champ sur $X$. Rappelons que d'après la définition de $R_{\Sigma}$ on a : 
$$R_{\Sigma}(\GGG)=(\{i_{k}^{-1} \GGG\}, \{i_{k}^{-1}\eta_{l}\}, \{\lambda_{klm}\})$$
où $\eta_{l}$ est le foncteur naturel : 
$$\eta_{l} : \GGG \longrightarrow i_{l*}i_{l}^{-1}\GGG$$
et $\lambda_{klm}$ est l'isomorphisme naturel : 
$$\shorthandoff{;:!?}
\xymatrix @!0 @C=3.2pc @R=2pc {i_{k}^{-1}\GGG \ar[rrr]^{i_{k}^{-1}\eta_{m}}  \ar[ddd]_{i_{k}^{-1}\eta_{m}} &&& i_{k}^{-1}i_{l*}i_{l}^{-1}\GGG \ar[ddd]^{i_{k}^{-1}i_{l*}i_{l}^{-1}\eta_{m}} \\
 && \ar@{=>}[ld]_\sim^{\lambda_{klm}}\\
& & &\\
 i_{k}^{-1}i_{m*}i_{m}^{-1}\GGG \ar[rrr]_{i^{-1}_k\eta_{lm}} &&&  i_{k}^{-1}i_{l*}i_{l}^{-1}i_{m*}i_{m}^{-1}\GGG
}$$
Considérons la famille de foncteurs et d'isomorphismes de foncteurs constituée : 
\begin{itemize}
\item[$\bullet$] pour tout objet $a \in \III$, d'un foncteur, noté $\Xi_{a}$, de $\GGG$ dans $R_{\Sigma}(\CCC)(a)$ :
\begin{itemize}
\item pour $k\leq n$ c'est le fonceur naturel :
$$\eta_{k} : \GGG \longrightarrow i_{k*}i_{k}^{-1}\GGG$$
\item pour les couples $(k,l)$ avec $k<l\leq n$ on se donne la composée de foncteurs : 
$$ \GGG  \buildrel\eta_{k}\over\longrightarrow i_{k*}i_{k}^{-1}\GGG  \buildrel i_{k*}i_{k}^{-1}\eta_{l}\over\longrightarrow i_{k*}i_{k}^{-1}i_{l*}i_{l}^{-1}\GGG$$
\item pour les triplets $(k,l,m)$ avec $k<l<m\leq n$ on se donne la composée :
$$\shorthandoff{;:!?}
\xymatrix @!0 @C=1.7cm @R=2cm {
\GGG \ar[r]^{\eta_{k}} & i_{k*}i_{k}^{-1}\GGG \ar[rr]^{i_{k*}i_{k}^{-1}\eta_{l}} && i_{k*}i_{k}^{-1}i_{l*}i_{l}^{-1}\GGG \ar[rrr]^{ i_{k*}i_{k}^{-1}i_{l*}i_{l}\eta_{m}} &&& i_{k*}i_{k}^{-1}i_{l*}i_{l}i_{m*}i_{m}^{-1}\GGG
}$$
\end{itemize}
\item[$\bullet$] pour tout foncteur $F :b\rightarrow a$ d'un isomorphisme de foncteur noté $\xi_{F}$  : 
\begin{itemize}
\item pour les foncteurs $(k,l) \rightarrow k$, $(k,l,m) \rightarrow (k,l)$ on se donne l'identité, 
\item pour le foncteur $(k,l) \rightarrow l$ on se donne l'isomorphisme naturel d'ajonction : 
$$\shorthandoff{;:!?}
\xymatrix @!0 @C=3.2pc @R=2pc{\GGG \ar[rrr]  \ar[ddd] &&& i_{k}^{-1}\GGG \ar[ddd] \\
 && \ar@{=>}[ld]_{\sim}\\
& & &\\
 i_{l*}i_{l}^{-1}\GGG \ar[rrr]_{\eta_{kl}} &&&  i_{k*}i_{k}^{-1}i_{l*}i_{l}^{-1}\GGG
}$$
\item et enfin pour le foncteur $(k,l,m) \rightarrow (k,m)$ on se donne la composée horizontale du morphisme identité, du foncteur $\eta_{k}$ et de l'isomorphisme $\lambda_{klm}$.
\end{itemize}
\end{itemize}
Ces données vérifient immédiatement les relations de compatibilité. Ainsi d'après la propriété universelle, il existe un foncteur, noté $\Xi$ :
$$\Xi : \GGG \longrightarrow Q_{\Sigma}R_{\Sigma}(\GGG)$$
et pour tout objet $a \in \III$ un unique isomorphisme noté $\xi_{a}$ :
$$\xi_{a} : \Xi \circ \pi_{a} \buildrel\sim\over\longrightarrow \Xi_{a} $$
Mais notons que, comme la restriction commute aux $2$-limites projectives finies, la restriction de $\Xi$ à $S_{j}$ est en fait le foncteur $\Phi_{j}$ du lemme \ref{phi}  construit en utilisant les données de $R_{\Sigma}(\GGG)$. Ainsi, la retriction de $\Xi$ à chaque strate est une équivalence. En particulier la fibre en chaque point de $X$ est une équivalence. On conclut alors en utilisant le corollaire \ref{equipre}.

\section{Champs constructibles et strictements constructibles}

Dans cette partie on s'intéresse à la notion  de champ constructible et strictement constructible. \\

Dans un premier paragraphe, on rappel les définitions et les premières propriétés de ces notions. La  notion de champs constructibles a été introduite par D. Treumann dans \cite{Tr1}. Elle a été introduite essentiellement parce que le champ des faisceaux pervers est, sur un espace stratifié de Thom-Mather, un champ constructible. 

Puis, dans un second temps, on considère l'espace topologique $\CC^n$ muni de la stratification, dite du croisement normal, associée à l'ensemble $\{(z_{1}, \ldots, z_{n}) \in \CC^n | z_{1}\cdots z_{n}=0\}$. Sur cette espace on peut donner une description particulièrement simple de la $2$-catégorie des faisceaux strictement constructibles relativement à cette stratification. Cela est dut au fait que les voisinages tubulaires des strates sont des produits.

Enfin on s'intéresse au champ $\PPP_{\CC^n}$ des faisceaux pervers sur $\CC^n$ muni de cette stratification. Nous montrons tout d'abord que ce champ est strictement constructible relativement à la stratification du croisement normal. Puis nous démontrons l'équivalence entre $\PPP_{\CC^n}$ et la donnée

\subsection{Définition et premières propriétés}~\\
 Pour les démonstrations des propositions et du théorème, le lecteur pourra se reporter à \cite{Tr1} et \cite{Tr2}.\\
Contrairement au paragraphe précédent $S_{j}$ désigne ici une strate de dimension quelconque, $i_{j}$ désigne l'injection suivante :
$$i_{j}: S_{j} \hookrightarrow X$$
\begin{Def}
Un champ $\CCC$ sur $X$ est dit constructible (resp. strictement constructible) relativement à $S$ si pour tout $i \in I$, la restirction $\CCC \mid_{S_{i}}$ est un champ localement constant (resp. constant)
\end{Def}
\noindent
Tout comme pour les faisceaux constructibles, on a les propositions suivantes : 
\begin{prop}\label{champ-constructible}
Soit $\CCC$ un champ constructible sur $X$ et $V \subset U$ deux ouverts de $X$ tels que l'injection $V \hookrightarrow U$ soit une équivalence stratifiée  d'homotopie. Alors la restriction $\CCC(U) \rightarrow \CCC(V)$ est une équivalence de catégories.
\end{prop}
\begin{prop}
Soit $\CCC$ un champ constructible sur $X$ relativement à une stratification de Thom-Mather et $x$ un point de $X$, alors il existe un ouvert $U$ de $X$ tel que le foncteur naturel suivant soit une équivalence : 
$$\Gamma(U, \CCC) \buildrel\sim\over\longrightarrow \CCC_{x}$$
\end{prop}
\begin{prop}\label{inv}
Soit $X$ et $Y$ deux espaces topologiques stratifiés et $f : X \rightarrow Y$ une application continue telle que l'image inverse d'une strate de $Y$ soit une réunion de strates de $X$, si $\CCC$ est un champ constructible (resp. strictement constructible)  alors $f^{-1}(\CCC)$ est constructible (resp. strictement constructible).
\end{prop}
\dem
L'image inverse par une application continue d'un champ localement constant (resp. constant) est localement constant (resp. constant), donc $f^{-1}(\CCC)\mid_{f^{-1}(\Sigma_{i})}$ est localement constant (resp. constant).
\cqfd
Un exemple particulièrement intéressant de champ constructible ou strictement constructible est le champ des faisceaux pervers relativement à une stratification fixée. Il justifie d'ailleurs à lui seul l'introduction de cette définition. 

Dans le cas général d'un espace topologique stratifié on ne peut rien dire. Mais dans le cas d'un espace de Thom-Mather on a le théorème suivant:
\begin{prop}
Soit $X$ un espace de Thom-Mather, on note $\Sigma$ la stratification de $X$, et $\PPP_{\Sigma}$ le champ des faisceau pervers sur $X$. Le champ $\PPP_{\Sigma}$ est constructible.
\end{prop}
\dem
Voir \cite{Tr1}.
\cqfd
Et, dans le cas où la projection des voisinages tubulaires sur les strates est une fibration trivial, $\PPP_{\Sigma}$ est un champ stritement constructible. Nous ne démontrons cela que dans le cas qui nous intéresse ( $\CC^n$ stratifié par le croisement normal) mais la démonstration est sensiblement la même. \\\\
\textbf{Remarque.}\\
Si $\CCC$ est un champ strictement constructible relativement une  stratification $\Sigma$, on s'attend à ce qu'on puisse en donner, via le théorème \ref{champ-constructible}, une description simple. En effet la restriction à chaque strate étant constante, la donnée d'une catégorie suffit à la décrire. Mais notons que, même si $\CCC_{L}$ est un champ constant, le champ $i_{K}^{-1}i_{L*}\CCC_{L}$ n'est, a priori, pas un champ constant, il est seulement localement constant. Ainsi la donnée des foncteurs de champs :
$$i_{k}^{-1}\CCC \longrightarrow i_{k}^{-1}i_{l*}i_{l}^{-1}\CCC$$ 
ne se réduit pas à la donnée d'un foncteur. \\
Nous allons voir, dans la partie suivante, que dans certains cas cette description peut être encore simplifiée.

\subsection{Champs strictement constructibles sur $\CC^n$ stratifié par le croisement normal}~\\
Considérons maintenant $\CC^n$ muni de la stratification du croisement normal défini  par l'ensemble $\{(z_{1}, \cdots, z_{n})\in \CC^n| z_{1}\cdots z_{n}=0\}$. Les strates sont indexées par les parties de $\{1,\cdots, n\}$. Ainsi on note $S_{K}$, pour $K$ partie de $\{1,\cdots,n\}$, la strate :
$$S_{K}=\prod_{i= 1}^n M_{i} \left \{ \begin{array}{cc}  M_{i}=\{0\} & \text{ si~}i \in K\\
              										M_{i}=\CC^*& \text{ si~} i\notin K
								\end{array} \right.
$$
 Dans ce paragraphe on reprend le théorème \ref{champ-constructible} dans le cas des champs strictement constructibles relativement à cette stratification. Le fait que les voisinages tubulaires de cette stratification soient des produits nous permet  de simplifier notablement l'énoncer de ce théorème. En effet, on montre que la $2$-catégorie des champs strictement constructibles relativement à cette stratification est $2$-équivalente à une $2$-catégorie, $\SSS_{\Sigma}^s$, dont les objets sont donnés par des catégories, des foncteurs et des isomorphismes de foncteurs. \\
Considérons les notations suivantes :\\
Soit $i_{K}$ l'injection de $S_{K}$ dans $\CC^n$ :
$$i_{K}: S_{K} \hookrightarrow \CC^n$$
Pour tout $K \subset \{1, \cdots, n\}$, on note $p_{K}$ le point suivant :
$$p_{K}= (\delta_{1}, \cdots, \delta_{n}) \left\{\begin{array}{ccc}
									     \delta_{i}=0 & i \in K\\
									     \delta_{i}=1 & i\notin K
         							             \end{array}\right.$$

et tout $p \notin K$, on note $\gamma_{Kp}$ le chemin de $S_{K}$ défini par : 
$$ \gamma_{Kp}(t)=(x_{1}(t), \cdots, x_{n}(t)) \text{~avec~} \left\{ \begin{array}{lccc}
                                                   x_{i}(t)= 0 &\text{~si~} i \in K\\
                                                   x_{i}(t)= e^{2i\pi t} & i=p\\
                                                   x_{i}(t)= 1 &\text{~sinon~}
                                                   \end{array} \right.\\$$
La famille $\big\{\gamma_{Kp}\big\}$ forme une famille génératrice du groupe fondamental de la strate $S_{K}$.\\
Avant de définir la $2$-catégorie $\SSS_{\Sigma}^s$ et de démontrer qu'elle est $2$-équivalente à la $2$-catégorie des champs strictement constructibles. Nous allons énoncer quelque lemmes que nous utilisons ultérieurement.
Soit $p$ la projection :
$$p : \CC^{*n-k} \times \CC^{*k-l}\times \{0\}^l \longrightarrow \CC^{*n-k}\times \{0\}^{k-l}\times \{0\}^l$$
et $\CCC_\Cc$ le champ constant de fibre $\Cc$ sur $\CC^{*n-k} \times \CC^{*k-l}\times \{0\}^l$, on note $\GGG$ le champ sur $\CC^{*n-k}\times \{0\}^{k-l}\times \{0\}^l$ défini par :
$$\GGG=p_*\CCC_\Cc$$
\begin{lem}
Le champ $\GGG$ est  équivalent au champ constant sur $\CC^{*n-k}\times \{0\}^{k-l}\times \{0\}^l$ de fibre $Rep(\pi_{1}((\CC^*)^{k-l}), \Cc)$.
\end{lem}
\dem
Notons que, pour $U$ un ouvert de \linebreak$\CC^{*n-k}\times \{0\}^{k-l}\times \{0\}^l$, les sections de $\GGG(U)$ sont données par :
$$\begin{array}{cclc}
\GGG(U) &=&\CCC_{\Cc}(p^{-1}(U))\\
 &=&Rep(\pi_{1}(p^{-1}(U)), \Cc)\\
 &\simeq&Rep(\pi_{1}(U)\times\pi_{1}(\CC^{*k-l}), \Cc)
\end{array}$$
En effet $p^{-1}(U)$ est isomorphe au produit $U\times \CC^{*k-l}$.\\
Soit $\tilde{\CCC}$ le pré-champ constant sur $\CC^{*n-k}\times \{0\}^{k-l}\times \{0\}^l$ de fibre \linebreak$Rep(\pi_{1}((\CC^*)^{k-l}, \Cc)$. 
Considérons le foncteur de champs, de source $\tilde{\CCC}$ et de but $\GGG$, défini par la donnée pour tout ouvert $U$ de  \linebreak$\CC^{*n-k}\times \{0\}^{k-l}\times \{0\}^l$ du foncteur :
$$\begin{array}{cccc}
\tilde{\CCC}(U)= Rep(\pi_{1}((\CC^*)^{k-l}, \Cc) &\longrightarrow & \GGG(U)=Rep(\pi_{1}(U)\times\pi_{1}(\CC^{*k-l}), \Cc)\\
F : \begin{array}{ccc}
       y &\mapsto& F(y) \\
       \gamma & \mapsto & F(\gamma)
       \end{array}& \longmapsto & F' : \begin{array}{ccc}
       (x,y) &\mapsto& F(y) \\
       (\gamma_{1}, \gamma_{2}) & \mapsto & F(\gamma_{2})
       \end{array}
\end{array}$$
Montrons que ce  foncteur est une équivalence sur les fibres. Soit $x=(x_{1}, \ldots, x_{n-k},0, \ldots, 0) \in \CC^{*n-k}\times \{0\}^k$, les produits polydisques ouverts centrés en $x_{i}$ et de $\{0\}^k$ forment une base de voisinage de $x$, or pour des rayons assez petits et comme les boules sont contractiles,  le foncteur associé est isomorphe à l'identité. Ainsi la limite est constante et isomorphe à l'identité. 
\cqfd
\begin{lem}\label{champCCn}
Soit $S_{K}$ et $S_{L}$ deux strates telles que $S_{K} \subset \overline{S}_{L}$ (c'est à dire telles que $K\supset L$), on note $k$ et $l$ les cardinaux de respectivement $K$ et $L$. Soit $\CCC_{L}$ un champ constant sur $S_{L}$ de fibre $\Cc$, alors le champ $i_{K}^{-1}i_{L*} \CCC_{L}$ est constant de fibre équivalente à $ Rep(\pi_{1}((\CC^*)^{k-l}), \Cc)$.  \\
\end{lem} 
\dem 
 On peut supposer sans perte de généralité que :
$$S_{K}= (\CC^*)^{n-k} \times \{0\}^{k} \text{~et~}S_{L}=(\CC^*)^{n-l} \times \{0\}^{l}$$
D'après le lemme précédent il suffit de montrer que $i_{K}^{-1}i_{L*}\CCC_L$ est équivalent au champ $\GGG$. On définit pour cela une équivalence du préchamp $\widetilde{i_{K}^{-1}i_{L*}\CCC_L}$ dans le champ $\GGG$.\\
Considérons les notations  : \\
si $\varepsilon=(\varepsilon_{1}, \ldots, \varepsilon_{n})$ est un n-uplet, on note $\bar{\varepsilon}=(\varepsilon_{1}, \ldots, \varepsilon_{n-k})$,\linebreak $\tilde{\varepsilon}=(\varepsilon_{n-k+1}, \ldots, \varepsilon_{n-l})$ et $\check{\varepsilon}=(\varepsilon_{n-l+1}, \ldots, \varepsilon_{n})$. De même la notation $B_{\tilde{x}}^{\tilde{\varepsilon}}$ désigne le polydisque ouvert centrés en $\tilde{x}$ et de rayon $\tilde{\varepsilon}$.\\
Soit $W$ un ouvert de $S_K$, $W$, est de la forme :
$$W=U\times \{0\}^k$$
où $U$ est un ouvert de $\CC^{*n-k}$. On a alors l'équivalence :
$$\begin{displaystyle}
\widetilde{i_K^{-1}i_{L*}\CCC_L}(W)\simeq 2\varinjlim_{V\supset W}\Gamma(V,i_{L*}\CCC_L)
\end{displaystyle}$$
Comme les boules centrées en zéro forment une base des ouverts contenant zéro, on a l'équivalence :
$$\begin{displaystyle}
    2\varinjlim_{V\supset W}\Gamma(V,i_{L*}\CCC_L) \simeq 2\varinjlim_{(\tilde{\varepsilon},\check{\varepsilon})}\Gamma(U\times B_0^{\tilde{\varepsilon}}\times B_0^{\check{\varepsilon}},i_{L*}\CCC_L)
  \end{displaystyle}$$
On a, par définition du foncteur $i_{L*}$ l'équivalence :
$$\begin{displaystyle}
    2\varinjlim_{(\tilde{\varepsilon},\check{\varepsilon})}\Gamma(U\times B_0^{\tilde{\varepsilon}}\times B_0^{\check{\varepsilon}},i_{L*}\CCC_L) \simeq 2\varinjlim_{\varepsilon}\Gamma(U\times B_0^{\tilde{\varepsilon}*} \times \{0\}^l,\CCC_L)
  \end{displaystyle}$$
Or, pour tout $\tilde{\varepsilon}$, $B_0^{\tilde{\varepsilon}}$ est homotopiquement équivalent, par une équivalence stratifiée, à $\CC^{*k-l}$. Ces équivalences commutant, la $2$-limite est constante et équivalente à  $\Gamma(U\times \CC^{*k-l}\times \{0\}^l)$.  Ainsi on a bien :
$$\widetilde{i_K^{-1}i_{L*}\CCC_L}(W)\simeq p_*\CCC_{\Cc} (W)$$
De plus, si $W'=U' \times \{0\}^{k} \subset W=U\times \{0\}^k$ sont deux ouverts de $S_K$ la restriction étant donnée par la restriction de $U$ à $U'$, le diagramme suivant commute :
$$\xymatrix{
\widetilde{i_K^{-1}i_{L*}\CCC_L}(W) \ar[r]^\sim \ar[d]_{\rho_{W'W}} & p_*\CCC_{\Cc} (W) \ar[d]^{\rho_{W'W}}\\
\widetilde{i_K^{-1}i_{L*}\CCC_L}(W') \ar[r]^\sim & p_*\CCC_{\Cc} (W')
}$$
Ainsi, la donnée, pour tout ouvert $W$ de $S_K$, de l'isomorphisme défini plus haut et, pour tout couple d'ouverts $W' \subset W$, de l'identité, définit une équivalence de préchamp entre $ \widetilde{i_K^{-1}i_{L*}\CCC_L}$ et $p_*\CCC_\Cc$. Ceci démontre en particulier que le préchamp $\widetilde{i_K^{-1}i_{L*}\CCC_L}$ est en fait le champ $i_K^{-1}i_{L*}\CCC_L$ et qu'il est bien équivalent au champ constant de fibre \linebreak$Rep(\pi_1(\CC^{*k-l}), \Cc)$.
\cqfd
\begin{rem}\label{auto}
Notons que si $J$ est une partie de $\{1,\cdots,n\}$ telle que $J\supset K\supset L$ et $j=\vert J\vert$ alors le champ $i_{J}^{-1}i_{K*}i_{K}^{-1}i_{L*}\CCC_{L}$ est encore un champ constant. On le voit en appliquant deux fois le lemme précédent. On obtient de plus que la fibre de ce champ est équivalente à $Rep((\CC^*)^{j-k}, Rep(\pi_{1}(\CC^{*k-l}),\Cc))$. 
 \end{rem}
 On note $\eta$ le foncteur de source $Rep((\CC^*)^{j-l}, \Cc)$ et de but \linebreak$Rep((\CC^*)^{j-k}, Rep(\pi_{1}(\CC^{*k-l}),\Cc))$  qui a un foncteur :
 $$ \begin{array}{ccc}
           (x_{1}, \ldots, x_{j-l}) & \mapsto & F(x_{1}, \ldots, x_{j-l})\\
           (\gamma_{1}, \ldots, \gamma_{j-l}) &\mapsto& F(\gamma_{1}, \ldots, \gamma_{j-l})
 \end{array}$$
associe le foncteur $F'$ de source $\Pi_{1}(\CC^{*j-k})$ et de but $Rep(\pi_{1}(\CC^{*k-l}), \Cc)$ qui à un point $(x_{1}, \ldots, x_{j-k})$ de $\CC^{*j-k}$ associe le foncteur :
$$\begin{array}{cccc}
(x_{j-k+1}, \ldots x_{j-l}) &\mapsto & F(x_{1}, \ldots,x_{j-l})\\
(\gamma_{j-k+1}, \ldots, \gamma_{j-l})& \mapsto & F(Id_{x_{1}}, \ldots, Id_{x_{j-k}}, \gamma_{j-k+1}, \ldots, \gamma_{j-l})
\end{array}$$

\begin{lem}
Le diagramme suivant commute à isomorphisme constant près :
$$\xymatrix{
i_{J}^{-1}i_{L*}\CCC_{L} \ar[r]^{i_{J}^{-1}\eta_{KL}} \eq[d] &i_{J}^{-1}i_{K*}i_{K}^{-1}i_{L*}\CCC_{L} \eq[d]\\
Rep((\CC^*)^{j-l}, \Cc)\ar[r]_{\eta~~~~~~~~~~~} &Rep((\CC^*)^{j-k}, Rep(\pi_{1}(\CC^{*k-l}),\Cc))
}$$
\end{lem}
\dem
Pour démontrer ce lemme, on reprend la démonstration du lemme précédent en construisant étape par étape un isomorphisme de foncteur. Une démonstration similaire est explicité dans la démonstration de la proposition $3.24$.
\cqfd
On  a la proposition suivante :
\begin{prop}\label{stricte}
La $2$-catégorie des champs strictement constructibles sur $\CC^n$ stratifié par le croisement normal est $2$-équivalente à la $2$-catégorie $\SSS^{s}$  dont 
\begin{itemize}
\item[$\bullet$] les objets sont donné par : 
\begin{itemize}
\item[-] pour tout $K\subset \{1, \cdots, n\}$, une catégorie $C_{K}$,
\item[-] pour tout couple $(K,L)$ de parties de $\{1, \ldots, n\}$ tel que \linebreak$ L \subset K$, un foncteur $F_{lk} : C_{K} \rightarrow Rep(\pi_{1}((\CC^*)^{k-l}), C_{L})$,
\item[-] pour tout triplet $(K,L,M)$ de parties de $\{1, \ldots, n\}$ vérifiant $M\subset L\subset K$ un isomorphisme de foncteurs $\lambda_{KLM}$ : 
$$\shorthandoff{;:!?}\relax
\xymatrix @!0 @C=2,5cm @R=0,6cm {C_{K} \ar[rrr]^{F_{LK}}  \ar[ddddd]_{F_{MK}} &&& Rep(\pi_{1}((\CC^*)^{k-l}),C_{L}) \ar[ddddd]^{Rep(\pi_{1}(\CC^{*k-l}),F_{ML})} \\
~\\
& & \ar@{=>}[ld]_\sim\\
&~\\
\\
Rep(\pi_{1}((\CC^*)^{m-k}),C_{M}) \ar[rrr]_{i^{-1}_k\eta_{lm}~~~~~~~~~} &&& Rep(\pi_{1}((\CC^*)^{k-l}),Rep(\pi_{1}((\CC^*)^{l-m}),C_{M}))}$$
\end{itemize}
tels que, pour tout $k>l>m>p$ les deux morphismes que l'on peut définir entre les foncteurs :
$$ Rep(\pi_{1}(\CC^{*l-k}), Rep(\pi_{1}(\CC^{*m-p}),F_{pm}))\circ Rep(\pi_{1}(\CC^{*l-k}), F_{ml}) \circ F_{lk}$$
$$\text{~et~}$$
$$Rep(\pi_{1}(\CC^{*k-m}), F_{pk})$$
soient égaux.
\item[$\bullet$] pour deux objets $(\{C_{K}\}, \{F_{LK}\}, \{\lambda_{KLM}\})$ et $(\{C'_{K}\}, \{F'_{LK}\}, \{\lambda'_{KLM}\})$ un foncteur est donné par : 
\begin{itemize}
\item[-] pour toute strate $S_{K}$, un foncteur $F_{K} : C_{K}\rightarrow C'_{K}$\\
\item[-] pour tout couple de strates $(S_{K},S_{L})$ tel que $\overline{S}_{L} \supset S_{K}$, un isomorphisme de foncteur $f_{KL}$ :
$$f_{KL} : F'_{LK} \circ F_{K} \buildrel\sim\over\longrightarrow Rep(\pi_{1}((\CC^*)^{k-l}),F_{L}) \circ F_{LK}$$
tels que les deux morphismes que l'on peut définir entre les foncteurs :
$$Rep(\pi_{1}(\CC^{*k-l}), F_{ml}')\circ F_{lk}\circ G_{k}$$
$$\text{et}$$
$$Rep(\pi_{1}(\CC^{*k-m}), G_{m})\circ F_{mk}$$
soit égaux.
\end{itemize}
\item[$\bullet$] Les morphismes entre deux tels foncteurs sont les données pour chaque $S_{K}$ d'un morphisme de foncteurs $\Phi_{k}: G_{k}\rightarrow G'_{k}$ tels que le diagramme suivant commute :
$$\xymatrix{
F'_{kl}\circ G_{k} \ar[d]_{Id\bullet\Phi_{k}} \ar[r]^{g_{kl}~~~~~} & Rep(\pi_{1}(\CC^{*k-l}),G_{l})\circ F_{lk} \ar[d]^{Rep(\pi_{1}(\CC^{*k-l}), \Phi_{k})\bullet Id}\\
F'_{kl}\circ G'_{k}\ar[r]_{g'_{kl}~~~~~} & Rep(\pi_{1}(\CC^{*k-l}), G'_{l})\circ F_{lk}
}$$
\end{itemize}

\end{prop}
\dem
Comme d'après les lemmes précédents les champs, les foncteurs de champs et les isomorphismes de foncteurs sont constants, il suffit de se donner leurs fibres en $p_{K}$.
\noindent
On peut même définir une $2$-équivalence de champ :
$$\begin{array}{cccc}
R_{\Sigma}^s: & \SSS _{\CC^n}  & \longrightarrow & \SSS^s\\
     & \CCC & \longmapsto &  (\{ \CCC_{p_{K}}\}_{k\leq n }, \{(\eta_{KL})_{p_{K}} \}, \{\lambda_{MLK})_{p_{K}}\})\\
     & G : \CCC \rightarrow \CCC' & \longmapsto & ( \{G_{p_{K}} \}, \{g_{LK}\})\\
     & \phi : G \rightarrow G' & \longmapsto & (\{(\phi)_{p_{K}}\})
\end{array}$$
Le $2$-foncteur $Q_{\Sigma}^s$ est la composée du foncteur $Q_{\Sigma}$ et du foncteur qui à un objet $(C_{K}, F_{LK}, \lambda_{KLM})$ de $\SSS_{\Sigma}^s$ associe l'objet $(\CCC_{K}, \boldsymbol{F}_{LK}, \boldsymbol{\lambda}_{KLM})$ de $\SSS_{\Sigma}$ formé des champs constants sur $S_{K}$ de fibre $C_{K}$, les foncteurs de fibre $F_{LK}$ ainsi que les isomorphismes constants de fibre $\lambda_{KLM}$.\\
\cqfd
\subsection{Le champs des faisceaux pervers sur $\CC^n$}~\\
Dans ce paragraphe on va montrer que  le champ, $\PPP_{\CC^n}$, des faiscceaux pervers sur $\CC^n$ relativement au croisemet normal est strictement constructible. Ensuite on explicite son image dans la catégorie $\SSS^s$. En particulier on montre que la fibre en $p_{K}$ de $\PPP_{\CC^n}$ est la catégorie $\Pc erv_{\CC^K}$ et que les fibres de $i_{K}^{-1}i_{L*}i_{L}^{-1}\PPP_{\CC^n}$ est équivalent à $\Pc erv_{\CC^{*K\backslash L }\CC^L}$, où $\CC^{*K\backslash L }\CC^L$ est l'ensemble :
\begin{Def}
On note $\CC^{*K\backslash L}\CC^L$ l'ensemble :
$$\Pi_{i=1}^nM_{i}\left\{\begin{array}{cccl}
                                 M_{i}&=& \CC^*&\text{si~} i \in K\backslash L\\
                                 M_{i}&=& \CC &\text{si~} i \in L\\
                                 M_{i}&=& \{1\} & \text{sinon}
                                 \end{array}\right.$$
De manière général on note $Z^KU^LF$ l'ensemble :
$$\Pi_{i=1}^nM_{i}\left\{\begin{array}{cccl}
                                 M_{i}&=& Z&\text{si~} i \in K\\
                                 M_{i}&=& U &\text{si~} i \in L\\
                                 M_{i}&=& F & \text{sinon}
                                 \end{array}\right.$$                                 
\end{Def}
Rappelons que dans notre cas, les voisinages tubulaires des strates sont des fibrations triviales. Pour commencer on étudie donc le comportement du champ des faisceaux pervers par rapport aux projections $X\times Y \rightarrow X$. On a le lemme suivant :
\begin{lem}\label{strictpervers}
Soit $X$ un espace métrique muni d'une stratification $\Sigma= \cup_{i \in I} \Sigma_i$ de Whitney, $B$ un espace métrique contractile et $\Fc$ un faisceau sur $X \times B$ à cohomologie constructible relativement à la stratification $\Sigma'= \cup_{i \in I} \Sigma_i \times B$. Considérons $p$ la première projection. Alors les morphismes naturels :
$$ p^{-1} \circ Rp_* (\Fc) \rightarrow \Fc$$
$$\Fc \rightarrow Rp_*\circ p^{-1} \Fc$$
sont des isomorphismes.
\end{lem}
\dem
Démontrons que ces morphismes sont des isomorphismes sur les fibres. \\
Commen\c cons par le morphisme $\Fc \rightarrow Rp_*\circ p^{-1} \Fc$. Soit $x \in X$, on a :
$$\begin{array}{cccc}
(Rp_{*}p^{-1}\Fc)_{x} &\simeq & \varinjlim_{U\ni x}R\Gamma(U, Rp_{*}p^{-1}\Fc)\\
& \simeq &  \varinjlim_{U\ni x}R\Gamma(U\times B, p^{-1}\Fc)
\end{array}$$
Mais d'après la formule de K\"unneth, on a :
$$\begin{array}{cccc}
R\Gamma(U\times B, p^{-1}\Fc)& \simeq R\Gamma(U, \Fc)\otimes R\Gamma(B, \CCC_{\CC})
\end{array}$$
où $\CCC_{\CC}$ est le faisceau constant sur $B$ de fibre $\CC$. Mais $R\Gamma(B, \CC)=\CC$, on a bien l'équivalence cherchée.\\
Soit maintenant $\Gc$ un faisceau à cohomologie constructible sur $X\times B$ et soit $(x,xy)$ un point de $X\times B$, on a l'isomorphisme suivant :
$$\begin{array}{cccc}
(p^{-1} \circ Rp_* (\Gc))_{(x,y)}&\simeq &\varinjlim_{U\times B'\ni (x,y)}R\Gamma(U\times B', p^{-1} \circ Rp_* (\Gc))\\
& \simeq & \varinjlim_{U\ni x}R\Gamma(U\times B,\Gc)\\
\end{array}$$
Mais comme pour toute boule $B'$ incluse dans $B$, l'injection de $U\times B'$ dans $U\times B$ est une équivalence d'homotopie stratifiée, on a l'équivalence :
$$\begin{array}{cccc}
\varinjlim_{U\ni x}R\Gamma(U\times B,\Gc) &\simeq & \varinjlim_{U\times B'\ni (x, y)}R\Gamma(U\times B',\Gc)\\
& \simeq & (\Gc)_{(x,y)}
\end{array}$$
\cqfd
On en déduit la proposition suivante :
\begin{prop}\label{proj}
Soient $X,Y$ deux espaces topologiques métriques, notons $p$ la projection $p : X \times Y \rightarrow X$,  fixons une stratification \linebreak $\Sigma= \bigcup_{i \in I}\Sigma_i $ de $X$ et considérons $\PPP_X$  le champ des faisceaux pervers sur respectivement $X$  relativement à la stratification $\Sigma$ et $\PPP_{X \times Y}$ le champ des faisceaux pervers sur $X \times Y$ relativement à la stratification $\Sigma=\bigcup_{i \in I} \Sigma_i \times Y$ alors 
$$p^{-1}( \PPP_X) \simeq \PPP_{X \times Y}$$ 

\end{prop}
\dem
On note $\widetilde{p^{-1}\PPP_X}$ le préchamp  défini par la donnée pour tout $U \times V$ ouverts de $X \times Y$ de la catégorie $\PPP_X(U)$. Soit $\tilde{F}$ le foncteur de préchamps $\tilde{F} :\widetilde{p^{-1}\PPP_X} \rightarrow  \PPP_{X\times Y}$ défini pour $U \times V$ un ouvert de $X \times Y$ par le foncteur  : 
$$\begin{array}{cccc}
   F_{U\times V} : & \Pc erv_U & \longrightarrow & \Pc erv_{U \times V}\\
                   & \Fc & \longmapsto & p^{-1}\Fc
  \end{array}$$
où l'on a abusivement noté $p$ la restriction de $p$ à $U \times V$.
L'isomorphisme de restriction est l'isomorphisme naturel.\\
Soit $(x,y)$ un point de $X \times Y$. La fibre en $(x,y)$ du foncteur $\tilde{F}$ est la limite inductive sur les ouverts de la forme $B \times B'$ où $B$ (resp. $B'$) est une boule centrée en $x$ (resp. $y$). Or d'après le lemme précdent $\tilde{F}_{B \times B'}$ est une équivalence pour tout $B \times B'$, donc $\tilde{F}_{(x,y)}$ est une équivalence de catégories. Ainsi, d'après le corollaire \ref{equipre} le champ $p^{-1}\PPP_X$ associé au préchamp $\widetilde{p^{-1}\PPP_X}$ est équivalent au champ $\PPP_{X \times Y}$.
\cqfd

On a alors le théorème suivante :
\begin{thm}\label{pervers-constructible}
Le champ, $\PPP_{\CC^n}$, des faisceaux pervers sur $\CC^n$ relativement au croisement normal est un champ strictement constructible relativement à cette stratification.
\end{thm}
\dem
Démontrons par récurrence.
Pour $n=1$ c'est évident puisque $\PPP \mid_{\CC^*}$ est le champ des faisceaux localement constants et que bien sûr $(\PPP_{1})_0$ est un champ constant.   \\ 

Supposons que $\PPP_{n-1}$ soit strictement constructible.\\
Comme ci-dessus on sait que $(\PPP_{n})_0$ est un champ constant. Soit $\Sigma_i$ une strate de dimension strictement positive, il existe alors $m \in \NN$ tel que $0 \leq m \leq n$, $\Sigma_i \subset \CC^m \times \CC^* \times \CC^{n-m-1}$. On a ainsi :
$$\PPP_n \mid_{\Sigma_i } \simeq (\PPP_n \mid_{\CC^{m} \times \CC^*\times \CC^{n-m-1}} )\mid_{\Sigma_i }$$
Or d'après la proposition \ref{proj} on sait que $\PPP_n \mid_{ \CC^m \times \CC^* \times \CC^{n-m-1}} $ est équivalent au champ $p_m^{-1}(\PPP_{n-1})$, où $p_m$ est la projection : 
$$p_m : \CC^{m} \times \CC^* \times \CC^{n-m-1} \rightarrow \CC^{n-1}$$
Ainsi $\PPP_n \mid_{ \CC^m \times \CC^* \times \CC^{n-m-1}}$ est un champ strictement contructible relativement à la stratification $p_m^{-1}(\Sigma)$ et $\Sigma_i$ est bien une strate de $p_m^{-1}(\Sigma)$ donc $\PPP_n \mid_{\Sigma_i }$ est constant.
\cqfd

On peut maintenant étudier l'image de $\PPP_{\CC^n}$ par $R_{\Sigma}^s$. Cette étude nous permettra par la suite de démontrer l'équivalence entre le champ $\PPP_{\CC^n}$ et un champ que l'on définira. Nous définissons un objet de $\SSS^s$ dont nous démontrons qu'il est équivalent à $R^s_\Sigma(\PPP_{\CC^n})$. On a besoin pour cela de la proposition et des lemmes suivants : 
\begin{prop}
 La catégorie $\Pc erv_{\CC^{*K\backslash L}\CC^L}$ est équivalente à la catégorie $Rep(\pi_1(\CC^{*k-l}),\Pc erv_{\CC^L})$.
\end{prop}
\dem
Notons $\CC^L$ l'ensemble $\CC^L\{1\}$.\\
Comme l'injection de $\CC^L$ dans $\CC^n$ est non caractéristique, la restriction à $\CC^L$, d'un faisceau pervers sur $\CC^n$ est  un faisceau pervers sur $\CC^L$ stratifié par le croisement normal. On a ainsi l'équivalence : 
$$\PPP_{\CC^n}|_{\CC^L}\simeq \PPP_{\CC^L}$$
On peut ainsi se placer sur $\CC^L$ et non $\CC^n$.\\
Notons que l'on a l'équivalence :
$$\Gamma(\CC^{*K\backslash L}, \PPP_{\CC^L})\simeq \Gamma(\CC^{*K\backslash L}\{0\}^L, \PPP_{\CC^L}$$
et comme la restriction à $\CC^{*K\backslash L}\{0\}^L$ du champ $\PPP_{\CC^L}$ est constante, on a bien l'équivalence cherchée.\\
Définissons une équivalence entre ces deux ca\-té\-go\-ries. La définition de cette équivalence s'inspire du cas des faisceaux localement constants. \\
Pour cela on se fixe $(1, \ldots,1)$ comme point base de $(\CC^*)^{k-l}$ et la famille $\{\gamma_{i}\}_{i\leq k-l}$ comme système de générateurs de $\pi_{1}(\CC^{*k-l})$, où les $\gamma_{i}$ sont les lacets : 
$$\begin{array}{cccccccc}
\gamma_{i} :& [0,1] & \longrightarrow & (\CC^*)^{k-l}\\
& t & \longmapsto & (\underbrace{1, \ldots1}_{i-1}, e^{2i\pi t},1, \ldots,1)
\end{array}$$
On définit d'abord un foncteur de la catégorie $\Pc erv_{\CC^{*K\backslash L}\CC^L}$ dans la catégorie $\Pc erv_{\CC^L}$.
$$\begin{array}{cccc}
 \Pc erv_{\CC^{*K\backslash L}\times \CC^L} & \longrightarrow & \Pc erv_{\CC^L} \\
          \Fc & \longmapsto & \Fc\mid_{ \CC^L}
\end{array}$$
Considérons maintenant l'application composée notée $\bar{\gamma_{i}}$ : 
$$\bar{\gamma_{i}} : [0,1]\times \CC^l \buildrel(\gamma_{i},Id)\over\longrightarrow (\CC^*)^{k-l} \times \CC^l \buildrel\sim\over\longrightarrow \CC^{*K\backslash L}\CC^L.$$
Notons que la restriction de $\bar{\gamma}_{i}$ au fermé $[0,1]\times\{0\}$ est un chemin $\gamma_{Lp}$, pour $p\in K\backslash L$, défini plus haut.\\
On a le lemme suivant :
\begin{lem}
Soit $\Fc$ un faisceau pervers défini sur $\CC^{*K\backslash L}\times \CC^L$ relativement à la stratification du croisement normal, alors $\bar{\gamma_{i}}^{-1}(\Fc)$ est un faisceau pervers sur $[0,1] \times \CC^l$ relativement à la stratification produit du croisement normal avec $[0,1]$.
\end{lem}
Ainsi, si $\Fc$ est un faisceau pervers sur $\CC^{*k-l}\times \CC^l$, comme le faisceau est contant le long de $[0,1]$ on a l'automorphisme, noté $\tilde{\gamma_{i}}(\Fc)$ suivant : 
$$\begin{array}{cccc}
\tilde{\gamma}_{i} : \Fc\mid_{\{1\}^{k-l}\times \CC^l} \simeq (\bar{\gamma_{i}}^{-1}\Fc)\mid_{\{0\}\times \CC^l}\simeq (\bar{\gamma_{i}}^{-1}\Fc)\mid_{\{1\}\times \CC^l} \simeq \Fc\mid_{\{1\}^{k-l}\times \CC^l}
\end{array}$$
On note alors $\mu_{KL}$ le foncteur de la catégorie $\Pc erv_{\CC^{*k-l}\times \CC^l}$ dans la catégorie $Rep(\pi_{1}(\CC^{*k-l}), \Pc erv_{\CC^L)}$ défini par :
$$\begin{array}{cccc}
\mu_{KL} : & \Pc erv_{\CC^{*K\backslash L}\CC^L} & \longrightarrow & Rep(\pi_{1}(\CC^{*k-l}), \Pc erv_{\CC^L})\\
 &\Fc & \longmapsto & (\Fc\mid_{\CC^L}, \{\tilde{\gamma_{i}}\}_{i\leq k-l})
\end{array}$$
C'est une équivalence de catégorie. \\
\cqfd
On a alors la proposition suivante :
\begin{prop}
L'image par $R_{\Sigma}^s$ du champ $\PPP_{\CC^n}$ est équivalente à la donnée : 
$$\Pc=\Big\{\{\Pc erv_{\CC^K}\}_{K\subset I}, \{\mu_{KL}\circ \rho_{\CC^K, \CC^{*K\backslash L }\times \CC^L}\} \{f_{KLM}\}\Big\}$$
\end{prop}
\dem
Définissons une équivalence dans $\SSS^s$ entre \linebreak$R_\Sigma^s(\PPP_\CC^n)$ et $\Pc$. Vu la définition de $\Pc$, il suffit de définir pour tout $K \subset \{1, \ldots, n\}$ une équivalence de catégorie :
$$\Xi_{K} : (\PPP_{\CC^n})_{p_{K}} \buildrel\sim\over\longrightarrow \Pc erv_{\CC^K}$$
et pour tout couple $K \subset L $ de parties de $\{1, \ldots, n\}$ une équivalence $\Xi_{KL}$ : 
$$\Xi_{KL} : (i_{L*}i_{L}^{-1}\PPP_{\CC^n})_{p_{K}} \longrightarrow \Pc erv_{\CC^{K\backslash L}\CC^L}$$
et un isomorphisme de foncteurs $\xi_{KL}$ : 
$$\shorthandoff{;:!?}\relax
\xymatrix @!0 @C=1.4cm @R=0,6cm {
(\PPP_{\CC^n})_{p_{K}} \ar[rrr]^{(\eta_{K})_{p_{K}}~~~~~}  \ar[ddddd]_{\Xi_{K}} &&&(i_{K}^{-1}i_{L*}i_{l}^{-1}\PPP_{\CC^n})_{p_{K}} \ar[ddddd]^{\Xi_{KL}} \\
~\\
& & \ar@{=>}[ld]_\sim^{\xi_{KL}}\\
&~\\
\\
\Pc erv_{\CC^K} \ar[rrr]_{\rho_{\CC^K,\CC^{*K\backslash L} \CC^L}}&&& \Pc erv_{\CC^{*K \backslash L} \CC^L}
}$$

Nous n'allons définir ces équivalences et isomorphismes de foncteurs que pour les strates $S_{K}$ et $S_{L}$ suivantes :
$$S_{K}= (\CC^*)^{n-k} \times \{0\}^k \text{~et ~}S_{L}= (\CC^*)^{n-l} \times \{0\}^l.$$
La généralisation à des strates quelconques n'est pas difficile mais nécessite de lourdes notations. \\
Commen\c cons par définir l'équivalence $\Xi_{K}$.\\
Rappelons les notations données précédemment, si $\varepsilon=(\varepsilon_{1}, \ldots, \varepsilon_{n})$ est un n-uplet, on note $\bar{\varepsilon}=(\varepsilon_{1}, \ldots, \varepsilon_{n-k})$, $\tilde{\varepsilon}=(\varepsilon_{n-k+1}, \ldots, \varepsilon_{n-l})$ et $\check{\varepsilon}=(\varepsilon_{n-l+1}, \ldots, \varepsilon_{n})$. De même la notation $B_{\tilde{x}}^{\tilde{\varepsilon}}$ désigne le produit des boules centrées en $x_{i}$ et de rayon $\varepsilon_{i}$. On a alors l'équivalence suivante :
$$\begin{displaystyle}
(\PPP_{\CC^n})_{p_{K}}  \simeq  2\varinjlim_{\varepsilon}\Gamma(B_{\bar{1}}^{\bar{\varepsilon}} \times B_{\tilde{0}}^{\tilde{\varepsilon}} \times B_{\check{0}}^{\check{\varepsilon}}, \PPP_{\CC^n})
\end{displaystyle}$$
Mais comme, pour $\varepsilon_{i}$ assez petit, chaque boule $B_{1}^{\varepsilon_{i}}$ est incluse dans $\CC^*$ et comme l'injection des boules centrées en zéro dans $\CC$ est une équivalence d'homothopie stratifiée la $2$-limite est constante et on a pour $\bar{\varepsilon}$ fixé :
$$\begin{displaystyle}
2\varinjlim_{\varepsilon}\Gamma(B_{\bar{1}}^{\bar{\varepsilon}} \times B_{\tilde{0}}^{\tilde{\varepsilon}} \times B_{\check{0}}^{\check{\varepsilon}}, \PPP_{\CC^n}) \simeq \Gamma(B_{\bar{1}}^{\bar{\varepsilon}} \times \CC^k, \PPP_{\CC^n})
\end{displaystyle}$$
Et d'après le lemme \ref{strictpervers} on a l'équivalence :
$$\begin{displaystyle}
\Gamma(B_{\bar{1}}^{\bar{\varepsilon}} \times \CC^k, \PPP_{\CC^n}) \simeq \Pc erv_{\CC^k}
\end{displaystyle}$$
On définit l'équivalence $\Xi_{K}$ comme la composée de ces trois équivalences.\\
Considérons maintenant la fibre en $p_{k}$ du champ $i_{K}^{-1}i_{L*}i_{L}^{-1}\PPP_{\CC^n}$. Comme ci-dessus on a l'équivalence  : 
$$\begin{displaystyle}
(i_{K}^{-1}i_{L*}i_{L}^{-1}\PPP_{\CC^n})_{p_{K}} \simeq  2\varinjlim_{\varepsilon}\Gamma(B_{\bar{1}}^{\bar{\varepsilon}} \times B_{\tilde{0}}^{\tilde{\varepsilon}} \times B_{\check{0}}^{\check{\varepsilon}}, i_{L*}i_{L}^{-1}\PPP_{\CC^n})
\end{displaystyle}$$
et d'après la définition du foncteur $i_{L*}$ on a l'égalité : 
$$\begin{displaystyle}
2\varinjlim_{\varepsilon}\Gamma(B_{\bar{1}}^{\bar{\varepsilon}} \times B_{\tilde{0}}^{\tilde{\varepsilon}} \times B_{\check{0}}^{\check{\varepsilon}}, i_{L*}i_{L}^{-1}\PPP_{\CC^n})  = 2\varinjlim_{\varepsilon}\Gamma(B_{\bar{1}}^{\bar{\varepsilon}} \times B_{\tilde{0}}^{\tilde{\varepsilon}*} \times \{0\}^{l}, i_{L}^{-1}\PPP_{\CC^n})
\end{displaystyle}$$
où $B_{\tilde{0}}^{\tilde{\varepsilon}*}= B_{\tilde{0}}^{\tilde{\varepsilon}} \cap \CC^{*k-l}$. 
Comme ci-dessus cette $2$-limite est constante car l'injection des boules centrées en zéro est une équivalence homotoique stratifiée, on a donc pour $\bar{\varepsilon}$ fixé :
$$\begin{displaystyle}
2\varinjlim_{\varepsilon}\Gamma(B_{\bar{1}}^{\bar{\varepsilon}} \times B_{\tilde{0}}^{\tilde{\varepsilon}*} \times \{0\}^{l}, i_{L}^{-1}\PPP_{\CC^n}) \simeq \Gamma(B_{\bar{1}}^{\bar{\varepsilon}} \times \CC^{*k-l} \times \{0\}^{l}, i_{L}^{-1}\PPP_{\CC^n})
\end{displaystyle}$$
Par définition de $\Gamma(W, \CCC)$, on a :
$$ \Gamma(B_{\bar{1}}^{\bar{\varepsilon}} \times \CC^{*k-l} \times \{0\}^{l}, i_{L}^{-1}\PPP_{\CC^n}) \simeq  \Gamma(B_{\bar{1}}^{\bar{\varepsilon}} \times \CC^{*k-l} \times \{0\}^{l}, \PPP_{\CC^n})$$
Mais, comme $B_{\bar{1}}^{\bar{\varepsilon}} \times B_{\tilde{0}}^{\tilde{\varepsilon}*} \times \{0\}^{l}$ est fermé dans l'ouvert paracompact $B_{\bar{1}}^{\bar{\varepsilon}} \times B_{\tilde{0}}^{\tilde{\varepsilon}*} \times \CC^{l}$ et d'après la proposition \ref{paracompact} le foncteur canonique suivant est une équivalence : 
$$\begin{displaystyle}
2\varinjlim_{\check{\varepsilon}}\Gamma(B_{\bar{1}}^{\bar{\varepsilon}} \times \CC^{*k-l} \times B_{\check{0}}^{\check{\varepsilon}}, \PPP_{\CC^n}) \buildrel\sim\over\longrightarrow  \Gamma(B_{\bar{1}}^{\bar{\varepsilon}} \times \CC^{*k-l} \times \{0\}^{l}, \PPP_{\CC^n})
\end{displaystyle}$$
Cette $2$-limite est constante et le foncteur canonique est une équivalence :
$$\begin{displaystyle}
\Gamma(B_{\bar{1}}^{\bar{\varepsilon}} \times \CC^{*k-l} \times \CC^l, \PPP_{\CC^n})\buildrel\sim\over\longrightarrow2\varinjlim_{\check{\varepsilon}} \Gamma(B_{\bar{1}}^{\bar{\varepsilon}} \times \CC^{*k-l} \times B_{\check{0}}^{\check{\varepsilon}}, \PPP_{\CC^n}) 
\end{displaystyle}$$
Mais d'après le lemme \ref{strictpervers} on a :
$$\begin{displaystyle}
\Gamma(B_{\bar{1}}^{\bar{\varepsilon}} \times \CC^{*k-l} \times \CC^l, \PPP_{\CC^n})\simeq \Pc erv_{\CC^{k-l}\times \CC^l} 
\end{displaystyle}$$
Ainsi on a l'équivalence, notée $\Xi_{KL}$ :
 $$\Xi_{KL} : (i_{K}^{-1}i_{L*}i_{L}^{-1}\PPP_{\CC^n})_{p_{K}} \simeq \Pc erv_{\CC^{*k-l}\times \CC^l}$$
Pour finir la démonstration du lemme il reste à définir l'isomorphisme :
$$\shorthandoff{;:!?}\relax
\xymatrix @!0 @C=1.4cm @R=0,6cm {
(\PPP_{\CC^n})_{p_{K}} \ar[rrr]^{(\eta_{KL})_{p_{K}}~~~~~}  \ar[ddddd]_{\Xi_{K}} &&&(i_{K}^{-1}i_{L*}i_{l}^{-1}\PPP_{\CC^n})_{p_{K}} \ar[ddddd]^{\Xi_{KL}} \\
~\\
& & \ar@{=>}[ld]_\sim^{\xi_{KL}}\\
&~\\
\\
\Pc erv_{\CC^k} \ar[rrr]&&& \Pc erv_{\CC^{*k-l}\times \CC^l}
}$$
où le foncteur du bas est la restriction usuelle. Considérons l'isomorphisme suivant, il provient de la $2$-fonctorialité  des $2$-limites projectives :
$$
\shorthandoff{;:!?}\relax
\xymatrix @!0 @C=2.5cm @R=0,6cm {
(\PPP_{\CC^n})_{p_{K}} \ar[rrr]^{(\eta_{KL})_{p_{K}}~~~~~}  \ar[ddddd]_{} &&&(i_{K}^{-1}i_{L*}i_{l}^{-1}\PPP_{\CC^n})_{p_{K}} \ar[ddddd]^{} \\
~\\
& & \ar@{=>}[ld]_\sim\\
&~\\
\\
\begin{displaystyle}2\varinjlim_{(\bar{\varepsilon}, \tilde{\varepsilon}, \bar{\varepsilon})}   \end{displaystyle}\Gamma(B_{\bar{1}}^{\bar{\varepsilon}}\times B_{\tilde{0}}^{\tilde{\varepsilon}}\times B_{\check{0}}^{\check{\varepsilon}}, \PPP_{\CC^n})  \ar[rrr]_{\begin{displaystyle}2\varinjlim_{(\bar{\varepsilon}, \tilde{\varepsilon}, \bar{\varepsilon})}   \end{displaystyle} \eta_{KL}}&&& \begin{displaystyle}2\varinjlim_{(\bar{\varepsilon}, \tilde{\varepsilon}, \bar{\varepsilon})}   \end{displaystyle}\Gamma(B_{\bar{1}}^{\bar{\varepsilon}}\times B_{\tilde{0}}^{\tilde{\varepsilon}*}\times \{0\}^l, \PPP_{\CC^n})\Big)
}$$
Mais par définition du foncteur naturel d'adjonction on a l'isomorphisme : 
$$
\shorthandoff{;:!?}\relax
\xymatrix @!0 @C=2.5cm @R=0,6cm {
\Gamma(B_{\bar{1}}^{\bar{\varepsilon}}\times B_{\tilde{0}}^{\tilde{\varepsilon}}\times B_{\check{0}}^{\check{\varepsilon}}, \PPP_{\CC^n}) \ar[rrr]^{\eta_{KL}}   &&& \Gamma(B_{\bar{1}}^{\bar{\varepsilon}}\times B_{\tilde{0}}^{\tilde{\varepsilon}*}\times \{0\}^l, \PPP_{\CC^n})  \\
~\\
& & \\
&\ar@{=>}[ru]^\sim~\\
\\
\Gamma(B_{\bar{1}}^{\bar{\varepsilon}}\times B_{\tilde{0}}^{\tilde{\varepsilon}}\times B_{\check{0}}^{\check{\varepsilon}}, \PPP_{\CC^n}) \ar[uuuuu]^{Id} \ar[rrr]&&& \Gamma(B_{\bar{1}}^{\bar{\varepsilon}}\times B_{\tilde{0}}^{\tilde{\varepsilon}*}\times B_{\check{0}}^{\check{\varepsilon}}, \PPP_{\CC^n}) \ar[uuuuu]
}$$
où le foncteur du bas est le foncteur de restriction. Par la suite nous considérerons  la $2$-limite de cet isomorphisme. On a de plus le morphisme suivant :
$$
\shorthandoff{;:!?}\relax
\xymatrix @!0 @C=2.5cm @R=0,6cm {
\Gamma(B_{\bar{1}}^{\bar{\varepsilon}}\times B_{\tilde{0}}^{\tilde{\varepsilon}}\times B_{\check{0}}^{\check{\varepsilon}}, \PPP_{\CC^n}) \ar[rrr]^{\rho} \ar[ddddd]   &&&\Gamma(B_{\bar{1}}^{\bar{\varepsilon}}\times B_{\tilde{0}}^{\tilde{\varepsilon}*}\times B_{\check{0}}^{\check{\varepsilon}}, \PPP_{\CC^n})  \ar[ddddd]\\
~\\
& & \ar@{=>}[ld]_\sim\\
&~\\
\\
\Pc erv_{\CC^k}  \ar[rrr]_{\rho}&&&  \Pc erv_{\CC^{*k-l}\times \CC^l}
}$$
Ce qui définit par propriété universelle de la $2$-limite inductive, un isomorphisme de foncteur :
$$
\shorthandoff{;:!?}\relax
\xymatrix @!0 @C=2.5cm @R=0,6cm {
\begin{displaystyle}2\varinjlim_{(\bar{\varepsilon}, \tilde{\varepsilon}, \bar{\varepsilon})}   \end{displaystyle}\Gamma(B_{\bar{1}}^{\bar{\varepsilon}}\times B_{\tilde{0}}^{\tilde{\varepsilon}}\times B_{\check{0}}^{\check{\varepsilon}}, \PPP_{\CC^n}) \ar[rrr] \ar[ddddd]   &&&\begin{displaystyle}2\varinjlim_{(\bar{\varepsilon}, \tilde{\varepsilon}, \bar{\varepsilon})}   \end{displaystyle}\Gamma(B_{\bar{1}}^{\bar{\varepsilon}}\times B_{\tilde{0}}^{\tilde{\varepsilon}*}\times B_{\check{0}}^{\check{\varepsilon}}, \PPP_{\CC^n})  \ar[ddddd]\\
~\\
& & \ar@{=>}[ld]_\sim\\
&~\\
\\
\Pc erv_{\CC^k}  \ar[rrr] &&&  \Pc erv_{\CC^{*k-l}\times \CC^l}
}$$
En composant les isomorphismes précédents correctement, on obtient l'isomorphisme $\xi_{KL}$ cherché.
\cqfd

%
%
%
%
%
%
%
%
%
%
%
%
%
%
%
%
%
%
%

\chapter{Catégorie et champ des faisceaux pervers sur $\CC^n$}

Dans ce paragraphe on considère la catégorie et le champ des faisceaux pervers sur $\CC^n$ stratifié par le croisement normal. \\

Dans un premier temps nous rappelons le résultat de M. Granger, A. Galligo et Ph. Maisonobe qui démontrent l'équivalence entre la catégorie des faisceaux pervers sur $\CC^n$ stratifié par le croisement normal et une sous-catégorie pleine de la catégorie des représentation du carquois $Q_{n}$. Puis nous transformons cette équivalence de catégories en équivalence de champs : nous définissons un champ de carquois équivalent au champ des faisceaux pervers. La définition de ce champ utilise la description d'un champ sur un espace stratifié et plus particulièrement sur des propriétés données par la stratification du croisement normal.  \\

Rappelons et donnons quelques notations dont nous nous servons dans ce chapitre.\\
Les strates de $\CC^n$ sont indéxées par les parties de $\{1,\cdots, n\}$. Ainsi on note $S_{K}$, pour $K$ partie de $\{1,\cdots,n\}$, la strate :
$$S_{K}=\prod_{i= 1}^n M_{i} \left \{ \begin{array}{cc}  M_{i}=\{0\} & \text{ si~}i \in K\\
              										M_{i}=\CC^*& \text{ si~} i\notin K
								\end{array} \right.
$$
Soit $i_{K}$ l'injection de $S_{K}$ dans $\CC^n$ :
$$i_{K}: S_{K} \hookrightarrow \CC^n$$
Pour tout $K \subset \{1, \cdots, n\}$ on note $p_{K}$ le point suivant :
$$p_{K}= (\delta_{1}, \cdots, \delta_{n}) \left\{\begin{array}{ccc}
									     \delta_{i}=0 & i \in K\\
									     \delta_{i}=1 & i\notin K
         							             \end{array}\right.$$

et tout $p \notin K$, on note $\gamma_{Kp}$ le chemin de $S_{K}$ défini par : 
$$ \gamma_{Kp}(t)=(x_{1}(t), \cdots, x_{n}(t)) \text{~avec~} \left\{ \begin{array}{lccc}
                                                   x_{i}(t)= 0 &\text{~si~} i \in K\\
                                                   x_{i}(t)= e^{2i\pi t} & i=p\\
                                                   x_{i}(t)= 1 &\text{~sinon~}
                                                   \end{array} \right.\\$$
La famille $\big\{\gamma_{Kp}\big\}$ forment une famille génératrice du groupe fondamental de la strate $S_{K}$.\\
Dans cette section si $K$ est une partie de $\{1, \cdots, n\}$, $k$ désigne son cardinal. \\
On considère la notation suivante. Si $Z,U,F$ sont des ensembles de $\CC$ et $J$, $K$, $L$ est une partition de $\{1, \ldots,n\}$ alors on note $Z^KU^LF^M$ l'ensemble :
$$Z^KU^LF= \prod_{i\in \{1, \ldots, n\}}M_{i}\text{~avec~}\left\{ \begin{array}{llcc}
                                                   M_{i}=Z &\text{~pour~} i \in J\\
                                                   M_{i}= U & \text{~pour~} i \in K\\
                                                   M_{i}= F & \text{~pour~} i \in L\\ 
                                                   \end{array} \right.$$
On donne des notations particulière à certain de ces ensembles :
$$\begin{array}{cclccc}
\CC^K&=&\displaystyle{ \prod_{\substack{i\in K}}M_{i}}\left\{ \begin{array}{lccc}
                                                   M_{i}=\CC & i \in K\\
                                                   M_{i}= \{1\} & i\notin K
                                                   \end{array} \right.\\
\end{array}$$                                                                                                  
Si $K$ est un ensemble de $\{1, \ldots,n\}$ on note $\overline{K}$ son complémentaire dans $\{1, \ldots, n\}$ :
$$\overline{K}=\{1, \ldots, n\} \backslash K$$
                                         
\section{Equivalence de Galligo, Granger, Maisonobe}

Considérons la catégorie $R(c_{n})$ des représentations du carquois $c_{n}$. Un objet de cette catégorie est la donnée, pour toute partie $K$ de $\{1,\cdots,n\}$,  d'un espace vectoriel $E_{K}$ et, pour tout couple de parties $(K, K \cup \{p\})$ de deux applications linéaires :
$$\begin{array}{cccc}
u_{Kp} : & E_{K}&\rightarrow &E_{K\cup p}\\
v_{Kp} :&E_{K\cup p} & \rightarrow & E_{K} 
\end{array}$$
\begin{Def}\label{Ccn}
Soit $\Cc_{n}$ la sous-catégorie pleine de la catégorie $R(Q_{n})$ formée des objets tels que :
\begin{itemize}
\item[(i)] pour tout couple de partie de $\{1,\cdots,n\}$ $K$, $ K \cup \{p\}$ on ait :
$$M_{Kp}= v_{Kp} u_{Kp} + Id \text{~soit inversible,}$$
\item[(ii)] pour tout quadruplet $K$, $K\cup \{p\}$, $K \cup \{q\}$ et $K \cup\{p,q\}$ de parties de $\{1,\cdots ,n\}$, les applications linéaires $u_{Kp}$, $u_{Kq}$, $u_{Kp,q}$, $u_{Kq,p}$ et  $v_{Kp}$, $v_{Kq}$, $v_{Kp,q}$, $v_{Kq,p}$ données par le diagramme :
$$\xymatrix{& E_{K\cup p} \ar@/^/[ld]^{v_{Kp}} \ar@/^/[rd]^{u_{Kpq}}&\\
E_{K}\ar@/^/[ur]^{u_{Kp}} \ar@/^/[dr]^{u_{Kq}}& &E_{K \cup \{p, q\}} \ar@/^/[ul]^{v_{Kpq}} \ar@/^/[dl]^{v_{Kqp}}\\
& E_{K \cup q} \ar@/^/[lu]^{v_{Kq}} \ar@/^/[ru]^{u_{Kqp}}
}$$ 
 vérifient les conditions suivantes : 
$$u_{Kp} u_{Kpq}= u_{Kq} u_{Kqp}, ~v_{Kpq}v_{Kp}=v_{Kqp} v_{Kq}, ~ v_{Kpq}u_{Kqp}=u_{Kp} v_{Kq}$$
\end{itemize}
\end{Def}
 
Dans \cite{GGM}, A.  Galligo, M. Granger et Ph. Maisonobe ont démontré le théorème suivant :
\begin{thm}
La catégorie $\Pc erv_{\CC^n}$ des faisceaux pervers sur $\CC^n$ stratifié par le croisement normal est équivalente à la catégorie $\Cc_{n}$.
\end{thm}
Ils démontrent ce théorème en définissant deux foncteurs $\alpha_{\CC^n}$ et $\beta_{\CC^n}$ quasi-inverse l'un de l'autre. On ne donnera ici que la définition du foncteur $\alpha_{\CC^n}$.
$$\alpha_{\CC^n} : \Pc erv_{\CC^n} \longrightarrow \Cc_{n}$$
Soit $\Fc$ un faisceau pervers sur $\CC^n$, l'image par $\alpha_{\CC^n}$ de $\Fc$ est la donnée : 
\begin{itemize}
\item[$\bullet$] pour toute partie $K$ de $\{1, \cdots, n\}$, de l'espace vectoriel $(R^k\Gamma_{\RR_{-}^{K}\CC\backslash\RR_{-}}\Fc)_{0}$ où $k$ est le cardinal de $K$, 
\item[$\bullet$] pour tout couple $K \subsetneq K \cup \{p\}$ de parties de $\{1, \cdots, n\}$, des fibres en zéro des morphismes naturels  : 
$$\begin{array}{ccccc}
u_{Kp}:&(R^k\Gamma_{\RR_{-}^{K}\CC\backslash\RR_{-}}\Fc)_{0}& \longrightarrow &(R^{k+1}\Gamma_{\RR_{-}^{K\cup p}\CC\backslash\RR_{-}}\Fc)_{0}\\
v_{Kp}:&(R^{k+1}\Gamma_{\RR_{-}^{K\cup p}\CC\backslash\RR_{-}}\Fc)_{0}& \longrightarrow & (R^{k+1}\Gamma_{\RR_{-}^{K}\RR_{-}^{*\{p\}}\CC\backslash\RR_{-}}\Fc)_{0}\\
\end{array}$$
\end{itemize}
Cette donnée est bien un objet de $\Cc_{n}$. La condition $(ii)$ est obtenue par fonctorialité des triangles considérés.\\
Pour démontrer que la condition $(i)$ est vérifiée, A. Galligo, M. Granger et Ph. Maisonobe démontrent les lemmes suivants :
\begin{lem}\label{concentre}
Soit $\Fc$ un faisceau pervers sur $\CC^n$ la restriction du faisceau $R^k\Gamma_{\RR_{-}^{K}\CC^{*}}\Fc$ à $\CC^K\CC^{*}$ est concentré en degré $k$, de plus ce faisceau est constructible relativement à la stratification réelle $\bigcup_{I}S_{I}\cap \RR_{-}^{K}\CC^{*}$. 
\end{lem}
Comme $S_{K} \subset \RR_{-}^{K}\CC^{*}$, la restriction à $S_{K}$ de $R^k\Gamma_{\RR_{-}^{K}\CC^{*}}\Fc$  est un faisceau localement constant. On note $E$ sa fibre en $p_{K}$ et $M_{Kp}$ les monodromies définies avec les chemins $\gamma_{Kp}$.\\
\begin{lem}\label{monodromie}
Pour tout $p \notin K$, on peut identifier la fibre en zéro des faisceaux $R^k\Gamma_{\RR_{-}^{K}\CC\backslash\RR_{-}^{}}\Fc$ et $R^k\Gamma_{\RR^K_{-}\RR_{-}^{*\{p\}}\CC\backslash\RR_{-}}\Fc$ avec la fibre en $p_{K}$ du faisceau $R\Gamma_{\RR_{-}^K\CC^*}\Fc$. De plus cet isomorphisme peut être choisi de fa\c con à ce que le diagramme suivant commute :
$$\xymatrix{
(R^k\Gamma_{\RR_{-}^{K}\CC^{*}}\Fc)_{p_{K}} \eq[d] \ar[r]^{M_{Kp}-Id}& (R^k\Gamma_{\RR_{-}^{K}\CC^{*}}\Fc)_{p_{K}} \eq[d]\\
(R^k\Gamma_{\RR_{-}^{K}\CC\backslash\RR_{-}^{}}\Fc)_{0} \ar[r] & (R^{k+1}\Gamma_{\RR_{-}^{K}\RR_-^{*\{p\}}\CC\backslash\RR_{-}^{}}\Fc)_{0}
}$$
où le morphisme du bas est le morphisme naturel donné par le triangle distingué : 
$$R\Gamma_{\RR_{-}^{K}\RR_-^{*\{p\}}\CC\backslash\RR_{-} }\Fc \longrightarrow R\Gamma_{\RR_{-}^K\CC^{*\{p\}}\CC \backslash\RR_{-}}\Fc \longrightarrow  R\Gamma_{ \RR_{-}^{K}\CC\backslash\RR_{-}}\Fc .$$
\end{lem}
\dem On peut supposer, sans perte de généralité, que $K=\{1, \cdots, k\}$ et $p=k+1$. On note $K'=\{1, \cdots, k+1\}$. On a ainsi :
$$\RR_{-}^K\CC^*=\RR_{-}^k\times \CC^{*n-k} \text{,} p_{K}=(0, \cdots,0,1, \cdots,1) $$
$$ \RR_{-}^{K}\RR_-^{*\{p\}}\CC\backslash\RR_{-} =\RR_{-}^k \times \RR_{-}^*\times\CC\backslash\RR_{-}^{n-k+1}$$
Notons $\Gc$ la restriction du faisceau $R^k\Gamma_{\RR_{-}^{k}\times\CC^{*n-k}}\Fc$ à $\CC^{k}\times\CC^{*n-k}$. 
$$\Gc=R^k\Gamma_{\RR_{-}^{k}\times\CC^{*n-k}}\Fc)|_{\CC^k\times\CC^{*n-k}} $$
Soit $i$ l'injection de $\CC^k\times\CC^{*n-k}$ dans $\CC^n$.
$$i : \CC^k\times \CC^{*n-k} \hookrightarrow \CC^n$$
Démontrons tout d'abord l'existence d'isomorphisme $n_{1}$ et $n_{2}$ tels que le diagramme suivant commute :
$$\xymatrix{
(R^0\Gamma_{\RR^k\times \CC\backslash\RR_{-}^{n-k}}\Fc)_{0} \ar[d]_{n_{1}}^{\sim} \ar[r] & (R^1\Gamma_{\RR^k\times\RR_{-}^* \times\CC\backslash\RR_{-}^{n-k}}\Fc)_{0} \ar[d]_{\sim}^{n_{2}}\\
(R^0\Gamma_{\RR^k\times \CC\backslash\RR_{-}\times \CC^{*n-k}}\Fc)_{p_{K'}} \ar[r] &  (R^1\Gamma_{\RR^k\times\RR_{-}^* \CC^{*n-k}}\Fc)_{p_{K'}}
 }$$
où les morphismes horizontaux sont les morphismes naturels. On a le diagramme commutatif suivant :
\begin{scriptsize}
$$\xymatrix{
(R^k\Gamma_{\RR^k\times \CC\backslash\RR_{-}^{n-k}}\Fc)_{0} \ar[r]  & (R^{k+1}\Gamma_{\RR^k\times\RR_{-}^* \times\CC\backslash\RR_{-}^{n-k}}\Fc)_{0} \\
\Gamma(\{0\}^{k+1}\times \CC\backslash\RR_{-}^{n-k-1}, R^k\Gamma_{\RR^k\times \CC\backslash\RR_{-}^{n-k}}\Fc) \ar[r]  \ar[u]^{\sim} \ar[d]_{\sim}& \Gamma(\{0\}^{k+1}\times \CC\backslash\RR_{-}^{n-k-1}, R^{k+1}\Gamma_{\RR^k\times\RR_{-}^* \times\CC\backslash\RR_{-}^{n-k}}\Fc)  \ar[u]_{\sim} \ar[d]^\sim\\
(R^k\Gamma_{\RR^k\times \CC\backslash\RR_{-}^{n-k}}\Fc)_{p_{K'}} \ar[r] & (R^{k+1}\Gamma_{\RR^k\times\RR_{-}^* \times\CC\backslash\RR_{-}^{n-k}}\Fc)_{p_{K'}}
}$$
\end{scriptsize}Chaque morphisme vertical est un isomorphisme car la restriction des faisceaux $R^k\Gamma_{\RR^k\times \CC\backslash\RR_{-}^{n-k}}\Fc$ et $R^{k+1}\Gamma_{\RR^k\times\RR_{-}^* \times\CC\backslash\RR_{-}^{n-k}}\Fc$ à $\{0\}^{k+1}\times \CC\backslash \RR_{-}^{n-k-1}$ est constante.\\
Considérons alors les morphismes naturels issus des triangles distingués :
$$\begin{array}{cccc}
R^k\Gamma_{\RR^k\times \CC\backslash\RR_{-}\times\CC^{*n-k-1}}\Fc &\longrightarrow & R^k\Gamma_{\RR^k\times \CC\backslash\RR_{-}^{n-k}}\Fc\\
R^{k+1}\Gamma_{\RR^k\times\RR_{-}^* \times\CC^{*n-k}}\Fc & \longrightarrow & R^{k+1}\Gamma_{\RR^k\times\RR_{-}^* \times\CC\backslash\RR_{-}^{n-k}}\Fc
\end{array}$$
Leurs fibres en $p_{K'}$ est un isomorphisme puisque les complémentaires de respectivement $\RR^k\times \CC\backslash\RR_{-}^{n-k}$ et $\RR^k\times\RR_{-}^* \times\CC\backslash\RR_{-}^{n-k}$ dans respectivement $\RR^k\times \CC\backslash\RR_{-}\times\CC^{*n-k-1}$ et $\RR^k\times \CC\backslash\RR_{-}^{n-k}$ ne contiennent pas $p_{K'}$. Les isomorphismes $n_{1}$ et $n_{2}$ sont les composés des ces isomorphismes.\\
D'après le lemme \ref{concentre} on a les isomorphismes naturels :
$$\begin{array}{rclc}
R^k\Gamma_{\RR_{-}^{k}\times\CC\backslash\RR_{-}^{n-k}}\Fc& \simeq & R^0\Gamma_{\CC^k\times \CC\backslash\RR_{-}^{n-k}}Ri_{*}\Gc \\
R^{k+1}\Gamma_{\RR_{-}^{k}\times \RR_{-}^*\times\CC\backslash\RR_{-}^{n-k}}\Fc&\simeq& R^1\Gamma_{\CC^k\times\RR_{-}^* \CC\backslash\times\RR_{-}^{n-k}}Ri_{*}\Gc
\end{array}$$
\\On peut alors appliquer le lemme plus général :
\begin{lem}
Soit $\Fc$ un faisceau localement constant sur $\CC^k\times \CC^{*n-k}$. Notons $E$ sa fibre en $p_{K}$ et $M$ sa monodromie définie par le lacet $\gamma_{k+1}$. Alors on peut identifier les fibres en $p_{K'}$ des faisceaux $R^0\Gamma_{\CC^k\times\CC\backslash\RR_{-}\times\CC^*}Ri_{*}\Fc$ et $R^1\Gamma_{\CC^k\times\RR_{-}^*\times\CC^*}Ri_{*}\Fc$  à $E$ de fa\c con à ce que le diagramme suivant commute :
$$\xymatrix{
(R^0\Gamma_{\CC^k\times\CC\backslash\RR_{-}\times\CC^*}Ri_{*}\Fc)_{p_{K'}} \ar[r] \ar[d]_{\sim} & (R^1\Gamma_{\CC^k\times\RR_{-}^*\times\CC^*}Ri_{*}\Fc)_{p_{K'}}\ar[d]^\sim \\
E \ar[r]_{M-Id} &E
}$$
\end{lem}
\cqfd
Pour démontrer que les applications $u_{Kp}$ et $v_{Kp}$ vérifient la condition $(i)$ de la définition \ref{Ccn}, il suffit alors de rappeler que l'on a le diagramme commutatif suivant :
$$\xymatrix{
(R^{k+1}\Gamma_{\RR_{-}^{K p}\CC\backslash\RR_{-}}\Fc)_{0} \ar[r]^{v_{Kp}} & (R^{k+1}\Gamma_{ \RR_{-}^{K}\RR_-^{*\{p\}}\CC\backslash\RR_{-}}\Fc)_{0} \\
(R^k\Gamma_{ \RR_{-}^{K}\CC\backslash\RR_{-}}\Fc)_{0} \ar[ru] \ar[u]^{u_{Kp}}
}$$

\section{''Stackification'' sur $\CC^n$}
Dans cette partie, nous définissons un champ $\CCC_{n}$ sur $\CC^n$ stratifié par le croisement normal équivalent au champ $\PPP_{\CC^n}$. On utilise pour cela la proposition \ref{stricte}.
\subsection*{Définition du champ $\CCC_{\CC^n}$}
Nous allons donc définir un champ de catégorie de carquois, $\CCC_n$, sur $\CC^n$ strictement constructible relativement à la stratification du croisement normal, on utilise pour cela la proposition \ref{stricte}.  \\
Sur chaque strate $S_{K}$ on se donne la catégorie $\Cc_{k}$ où $k$ est le cardinal de $K$. Pour tout couple $L \subset K$ de parties de $\{1, \cdots, n\}$ il faut donc se donner un  foncteur $F_{LK}$ : 
$$F_{LK} : \Cc_{k} \longrightarrow Rep(\pi_{1}(\CC^{*k-l}), \Cc_{l})$$
Rappelons que la catégorie $Rep(\pi_{1}(\CC^{*k-l}), \Cc_{l})$ est équivalente à la ca\-té\-go\-rie dont les objets sont des objets de $\Cc_{l}$ munis de $k-l$ automorphismes, ainsi un objet de   $Rep(\pi_{1}(\CC^{*k-l}), \Cc_{l})$ est la donnée :
\begin{itemize}
\item[$\bullet$] pour toute partie $J$ de $L$ d'un espace vectoriel noté $E_{J}$
\item[$\bullet$] pour tout couple de parties $J\subsetneq J\cup \{p\}$ de $L$ de deux applications linéaires :
$$\begin{array}{cccccc}
u_{Jp} : &E_{J} & \longrightarrow & E_{Jp}\\
v_{Jp} : &E_{Jp} & \longrightarrow &E_{J}
\end{array}$$
\item[$\bullet$] et pour toute partie $J$ de $L$ de $k-l$ endomorphismes : 
$$M_{Ji} : E_{J} \longrightarrow E_{J}$$
qui commutent aux applications $u_{Jp}$ et $v_{Jp}$.
\end{itemize} 
\begin{Def}                                                                                      
Soient $K$ et $L$ deux parties de $\{1,\cdots, n\}$ telles que $K\supset L$, on note  $F_{LK}$ le foncteur défini par :
$$\begin{array}{cccc}
F_{LK} : & \Cc_{k} &\longrightarrow & Rep(\pi_{1}(\CC^{*k-l}), \Cc_{l})\\
& (\{E_{J}\}_{J \subset K},\{u_{Jp}\},\{v_{Jp}\}) &\longmapsto & (\{E_{J}\}_{J\subset L}, \{u_{Jp}\},\{v_{Jp}\}, \{M_{Jp}\}_{p \notin J})
\end{array}$$
où $k=\vert K\vert$ et  $M_{Jp}$ est l'endomorphisme de $E_{J}$ :
$$M_{Jp}= v_{Jp}\circ u_{Jp}+ Id$$
\end{Def}
Les conditions (ii) de la définition \ref{Ccn} de la catégorie $\Cc_{n}$ assurent que $M_{Jp}$ commute aux applications $u_{Jp}$ et $v_{Jp}$.\\
\textbf{Exemple}\\
Pour $n=2$, $S_{\{1\}}$ et $S_{\{2\}}$ sont les  strates $S_{1}=\{0\} \times \CC^*$ et $S_{2}= \CC^*\times \{0\} $ et $F_{\{12\}\{1\}}$ est le foncteur :
$$\begin{array}{cccc}
F_{\{12\}\{1\}} :& \Cc_{2} & \longrightarrow & Rep(\pi_{1}(\CC^{*}), \Cc_{l}) \\
~\\
 & \xymatrix{E_{1} \ar@/^/[r]^{u_{\{1\}2}} \ar@/^/[d] &E_{12}\ar@/^/[d] \ar@/^/[l]^{v_{\{1\},2}} \\
 E_{\emptyset} \ar@/^/[u] \ar@/^/[r] & E_{2} \ar@/^/[u] \ar@/^/[l]}
 & \longmapsto & \xymatrix{ E_{1}\ar@(dr,ur)[]_{M_{\{1\}2}} \ar@/^/[d]\\
                                                   E_{\emptyset} \ar@(dr,ur)[]_{M_{\emptyset2}} \ar@/^/[u]}
\end{array}$$	
Soient maintenant un triplet, $J\subsetneq K \subsetneq L$ de partie de $\{1, \cdots n\}$  considérons le diagramme suivant : 
$$ \shorthandoff{;:!?}
\xymatrix@!0 @R=3cm @C=8cm{
\Cc_{l} \ar[r]^{F_{KL}} \ar[d]_{F_{JL}} & Rep(\pi_{1}(\CC^{*k-l}),\Cc_{k}) \ar[d]^{Rep(\pi_{1}(\CC^{*k-l}),F_{JK})} \\
Rep(\pi_{1}(\CC^{j-l}), \Cc_{j}) \ar[r]_{\eta~~~~~~~~} & Rep(\pi_{1}(\CC^{*k-l}),Rep(\pi_{1}(\CC^{*j-k}),\Cc_{j}))
}$$	
D'après la remarque \ref{auto}, la catégorie \linebreak $Rep(\pi_{1}(\CC^{*k-l}),Rep(\pi_{1}(\CC^{*j-k}),\Cc_{j}))$ est formée d'objets de $\Cc_{j}$ munis de $l-j$ automorphismes.\\	
Le morphisme du bas est bien un  isomorphisme puisqu'on a bien \linebreak $\CC^{*j-l}=\CC^{*k-l}\times \CC^{*j-k}$. \\
Pour comprendre le foncteur $Rep(\pi_{1}(\CC^{*k-l}),F_{JK})$ explicitons le dans le cas de la dimension deux, avec $J=\emptyset$, $K=\{1\}$ et $L=\{1,2\}$. On a alors $S_{\emptyset}=\CC^* \times \CC^*$, $S_{\{1\}}=\{0\} \times \CC^*$ et $S_{\{1,2\}}=\{0\} \times \{0\}$. Le foncteur $Rep(\pi_{1}(\CC^{*k-l}),F_{JK})$ est alors le foncteur suivant :
$$\begin{array}{cccc}
 Rep(\pi_{1}( \CC^*), \Cc_{1})  & \longrightarrow &Rep(\pi_{1}( \CC^*), Rep(\pi_{1}(\CC^*), \Cc_{l}) \\
~\\
   \xymatrix{ E_{1}\ar@(dr,ur)[]_{M_{\{1\}2}} \ar@/^/[d]\\
                                                   E_{\emptyset} \ar@(dr,ur)[]_{M_{\emptyset2}} \ar@/^/[u]}
 & \longmapsto &\xymatrix{ ~\\
                                                   E_{\emptyset} \ar@(lu,ru)[]^{M_{\emptyset1}}\ar@(dr,ur)[]_{M_{\emptyset2}} }
\end{array}$$	
Il apparaît alors clairement que le diagramme précédent commute, c'est à dire qu'on a l'égalité :
$$\eta \circ F_{JL}=Rep(\pi_{1}(\CC^{*k-l}), F_{JK}) \circ F_{KL}$$ 
ainsi la famille $(\{\Cc_{k}\}, \{F_{LK}\}, \{Id\})$ est un objet de $\SSS^s$. 		
\begin{Def}
On note $\CCC_n$ le champ image par $Q^s$ de   l'objet de $\SSS^s$ défini par la donnée  :
\begin{itemize}
\item pour toute strate $S_{K}$ de $\CC^n$ de la catégorie $\Cc_{k}$, 
\item  pour tout couple de strates, $S_{K}$ et $S_{L}$, telles que $S_{K} \subset \overline{S}_{L}$, du foncteur $F_{LK}$,
\item et pour tout triplet de strates $S_{J}$, $S_{K}$, $S_{L}$ tel que $J \subsetneq K \subsetneq L$ du morphisme de foncteur $Id$.
\end{itemize}
\end{Def}
\begin{thm}
Le champ $\CCC_n$ est équivalent au champ $\PPP_{\CC^n}$
\end{thm}
\dem
Pour démontrer ce théorème nous allons en fait définir une équivalence dans $\SSS^s_{\CC^n}$ entre $R^s(\PPP_{\CC^n})$ et $(\{ \Cc_{K}\}, \{ F_{LK}\}, \{Id\})$. On rappel que   : $$R^s(\PPP_{\CC^n})\simeq\Big\{\{\Pc erv_{\CC^K}\}_{K\subset I}, \{\mu_{KL}\circ \rho_{\CC^K, \CC^{*K\backslash L}\times \CC^L}\}_{L\subset K} \{f_{KLM}\}\Big\}.$$ 
Ainsi il faut définir pour tout $K\subset \{1, \ldots, n\}$ une équivalence entre $\Pc erv_{\CC^K}$ et $\Cc_{k}$. Nous utilisons évidemment les foncteurs $\alpha_{\CC^K}$ définis par Galligo, Granger Maisonobe. Nous devons aussi définir pour tout couple de strates $(S_{K}, S_{L})$ tel que $S_{K}\subset \overline{S}_{L}$, un isomorphisme de foncteur : 
$$
\shorthandoff{;:!?}\relax
\xymatrix @!0 @C=2.5cm @R=0,6cm {
\Pc erv_{\CC^K} \ar[rrr]^{\alpha_{K}} \ar[ddddd]_{\mu_{KL}\circ \rho}   &&&\Cc_{K} \ar[ddddd]^{F_{LK}}\\
~\\
& & \ar@{=>}[ld]^{\sim}_{a_{KL}}\\
&~\\
\\
Rep(\pi_{1}(\CC^{*k-l}), \Pc erv_{\CC^L}) \ar[rrr]_{Rep(\pi_{1}(\CC^{*k-l}), \alpha_{L})} &&& Rep(\pi_{1}(\CC^{*k-l}), \Cc_{L}) }$$
qui vérifient des relations de commutation.\\
Rappelons les définitions des foncteurs $\alpha_{K}$ et $\mu_{KL}$. \\
L'image par $\alpha_{K}$ d'un faisceau pervers $\Fc$ sur $\CC^K$ est la donnée : 
\begin{itemize}
\item[$\bullet$] pour toute partie $J$ de $K$, de l'espace vectoriel $$(R^j\Gamma_{(\RR^J \CC\backslash\RR_-)\cap \CC^K}\Fc)_{p_{K}},$$ 
\item[$\bullet$] pour tout couple $J \subsetneq J \cup \{p\}$ de parties de $K$, des fibres en $p_{K}$ des morphismes naturels  : 
$$\begin{array}{ccccc}
u_{Jp}:&(R^j\Gamma_{(\RR_-^J\CC\backslash\RR_-)\cap\CC^K}\Fc)_{p_{K}}& \longrightarrow &(R^{j+1}\Gamma_{(\RR_-^{J\cup p}\CC\backslash\RR_-)\cap \CC^K}\Fc)_{p_{K}}\\
v_{Jp}:&(R^{j+1}\Gamma_{(\RR_-^{J\cup p}\CC\backslash\RR_-)\cap \CC^K}\Fc)_{p_{K}}& \longrightarrow & (R^{j+1}\Gamma_{ (\RR_-^{J\cup p}\RR_-^{*\{p\}}\CC\backslash\RR_-)\cap \CC^K}\Fc)_{p_{K}}\\
\end{array}$$
\end{itemize}

Soit maintenant $\Fc$ un faisceau pervers sur $\CC^K$. l'image par $\mu_{KL}$ de $\Fc$ est donné par :
$$\begin{array}{cccc}
   \mu_{KL} : & \Pc erv_{\CC^L\CC^{*K\backslash L}} : &\longrightarrow & Rep(\pi_{1}(\CC^{*k-l}),\Pc erv_{\CC^L})\\
   & \Fc & \longmapsto & (\Fc|_{\CC^L}, \{\tilde{\gamma}_{i}\})
  \end{array}$$
où $\tilde{\gamma_{i}}$ sont les isomorphismes définis par :
$$\tilde{\gamma}_{i}\Fc|_{\CC^L}\simeq \bar{\gamma}_{i}^{-1}(\Fc)|_{\{0\}\times\CC^l}\simeq  \bar{\gamma}_{i}^{-1}(\Fc)|_{\{1\}\times\CC^l}\simeq \Fc|_{\CC^L}$$
avec $\bar{\gamma}_{i}$ l'application :
$$\bar{\gamma_{i}} : [0,1]\times \CC^l \buildrel(\gamma_{p},Id)\over\longrightarrow (\CC^*)^{k-l} \times \CC^l \buildrel\sim\over\longrightarrow \CC^{*K\backslash L}\CC^L.$$
et $\gamma_{i}$ le chemin :
$$\begin{array}{cccc}
\gamma_{p} : & [0,1] & \longrightarrow & \CC^{*k-l}\\
 & t & \longmapsto & (\underbrace{1,\cdots,1}_{i-1}, e^{2i\pi t}, 1\cdots,1)
\end{array}$$
Notons que la restriction de $\bar{\gamma}_{i}$ à $[0,1]\times \{0\}^l$ est $\gamma_{Lp}$ pour un $p\in K\backslash L$.\\
Ainsi pour définir les isomorphismes $a_{KL}$ il faut définir un isomorphisme fonctoriel entre les objets :
$$\Big(\big\{(R^j\Gamma_{\RR_{-}^J\CC\backslash\RR_{-}}\Fc)_{0}\big\}_{J\subset L},\big\{u_{Jp}, v_{Jp}\big\}_{\substack{J\subset L\\p\in L\backslash J}}, \big\{M_{Jp}\big\}_{\substack{J\subset L\\p\in K\backslash L}} \Big)$$
$$\txt{~et~}$$
$$\Big(\big\{(R^j\Gamma_{(\RR_{-}^J\CC\backslash\RR_{-})\cap \CC^L}\Fc\mid_{\CC^L})_{p_{K}}\big\}_{J \subset L}, \big\{u'_{Jp},v'_{Jp}\big\}_{\substack{J\subset L\\p\in L\backslash J}}, \big\{M_{Jp}'\big\}_{\substack{J\subset L\\p\in K\backslash L}}\Big)$$
où $u_{Jp}$ et $v_{Jp}$ sont respectivement les fibres en $0$ des morphismes naturels : 
$$R^j\Gamma_{\RR_{-}^J\CC\backslash\RR_{-}}\Fc \longrightarrow R^{j+1}\Gamma_{\RR_{-}^{J\cup p}\CC\backslash\RR_{-}}\Fc$$
$$R^{j+1}\Gamma_{\RR_{-}^{J\cup p}\CC\backslash\RR_{-}}\Fc \longrightarrow R^{j+1}\Gamma_{\RR_{-}^J\RR_{-}^{*\{p\}}\CC\backslash\RR_{-}}\Fc$$
et où $u'_{Jp}$ et $v'_{Jp}$ sont respectivement les fibres en $p_{K}$ des morphismes naturels : 
$$R^j\Gamma_{(\RR_{-}^J\CC\backslash\RR_{-})\cap \CC^L}(\Fc\mid_{\CC^L}) \longrightarrow R^{j+1}\Gamma_{(\RR_{-}^{J\cup p}\CC\backslash\RR_{-})\cap \CC^L}(\Fc\mid_{\CC^L})$$
$$R^{j+1}\Gamma_{(\RR_{-}^{J\cup p}\CC\backslash\RR_{-})\cap \CC^L}(\Fc\mid_{\CC^L}) \longrightarrow R^{j+1}\Gamma_{(\RR_{-}^J\RR_{-}^{*\{p\}}\CC\backslash\RR_{-})\cap \CC^L}(\Fc\mid_{\CC^L})$$
et où $M_{Jp}$ est d'après le lemme \ref{monodromie} la monodromie de la restriction à $S_J$ du faisceau $R^j\Gamma_{\RR^J\CC^*}\Fc$, définie par le chemin $\gamma_{Jp}$, et $M_{Jp}'$ est l'image par $Rep(\pi_{1}(\CC^{k-l}), \alpha_{L})$ de l'isomorphisme $\tilde{\gamma}_{i}$ tel que la restriction à $[0,1]\times \{0\}$ de $\bar{\gamma}_{i}$ soit égale à $\gamma_{Jp}$. Ainsi $M_{Jp}'$ est l'endomorphisme :
$$M_{Jp}'=R^j\Gamma_{(\RR^J\CC\backslash \RR_{-})\cap \CC^L} \tilde{\gamma}_{i}$$

Commen\c cons par identifier fonctoriellement $\big\{(R^j\Gamma_{\RR_-^J\CC\backslash \RR_-}\Fc)_{0}\big\}_{J\subset K}$ et  $\big\{(R^j\Gamma_{\RR_-^J\CC\backslash \RR_-\cap \CC^K}\Fc\mid_{\CC^K})_{p_{K}}\big\}_{J \subset K}$, de cette manière on identifiera aussi $\big\{u_{Jp}\big\}$ avec $ \big\{u'_{Jp}\big\}$. \\
Rappelons que, d'après le lemme \ref{monodromie}, on a l'isomorphisme naturel :
$$(R^j\Gamma_{\RR_-^J\CC\backslash \RR_-}\Fc)_{0} \simeq (R^j\Gamma_{\RR_-^J\CC\backslash \RR_-}\Fc)_{p_{K}}$$
L'isomorphisme que l'on cherche est la composée de cet isomorphisme et de l'isomorphisme donné par le lemme suivant.
\begin{lem}
Soit $\Fc$ un faisceau pervers sur $\CC^n$ relativement au croisement normal. Pour tout couple $K \subset L$  de parties de $\{1, \cdots, n\}$ les complexes suivants sont fonctoriellement isomorphes : 
$$(R\Gamma_{\RR_-^K\CC\backslash \RR_-}\Fc)\mid_{\CC^L}\simeq R\Gamma_{(\RR_-^K\CC\backslash \RR_-)\cap \CC^L}(\Fc\mid_{\CC^L})$$
\end{lem}
\dem
Définissons tout d'abord un morphisme naturel $n_{K}$ de $(R\Gamma_{\RR_-^K\CC\backslash \RR_-}\Fc)\mid_{\CC^L}$ dans $R\Gamma_{(\RR_-^K\CC\backslash \RR_-)\cap \CC^L}(\Fc\mid_{\CC^L})$.
$$n_{K} : (R\Gamma_{\RR_-^K\CC\backslash \RR_-}\Fc)\mid_{\CC^L} \longrightarrow R\Gamma_{(\RR_-^K\CC\backslash \RR_-)\cap \CC^L}(\Fc\mid_{\CC^L})$$
Notons $j$ l'injection de $\CC^L$ dans $\CC^n$ :
$$j: \CC^L \hookrightarrow \CC^n$$
On a l'égalité :
$$R\Gamma_{\RR_{-}^K\CC\backslash\RR_{-}}\Fc =R\Hc om(\underline{\CC}_{\RR_{-}^K\CC\backslash\RR_{-}}, \Fc)$$
où $\underline{\CC}$ est le faisceau constant sur $\CC^n$ de fibre $\CC$. On a alors le morphisme :
$$j^{-1}R\Hc om_{\underline{\CC}}(\underline{\CC}_{\RR^K\CC\backslash \RR_{-}}, \Fc)\longrightarrow R\Hc om_{j^{-1}\underline{\CC}}(j^{-1}(\underline{\CC}_{\RR^K\CC\backslash \RR_{-}}), j^{-1}\Fc)$$
Mais comme $\RR_{-}^K\CC\backslash \RR_{-}$ est localement fermé, on a l'isomorphisme :
$$R\Hc om_{j^{-1}\underline{\CC}}\big(j^{-1}(\underline{\CC}_{\RR_{-}^K\CC\backslash \RR_{-}}), j^{-1}\Fc\big)\simeq R\Hc om_{j^{-1}\underline{\CC}}\big((j^{-1}\underline{\CC})_{(\RR_{-}^K\CC\backslash\RR_{-})\cap \CC^L}, j^{-1}\Fc\big)$$
Comme l'image inverse d'un faisceau constant est un faisceau constant on a bien :
$$R\Hc om_{j^{-1}\underline{\CC}}\big((j^{-1}\underline{\CC})_{(\RR_{-}^K\CC\backslash\RR_{-})\cap \CC^L}, j^{-1}\Fc\big) =R\Gamma_{(\RR_{-}^K\CC\backslash\RR_{-})\cap \CC^L}(\Fc|_{\CC^L})$$
La composition de ces morphismes est le morphisme $n_{K}$.\\
Montrons maintenant que $n_{K}$ est un isomorphisme. \\
Considérons pour cela la restriction à $\CC^L$ du triangle distingué suivant :
$$R\Gamma_{\CC^n\backslash (\CC \backslash \RR_{-}^n)} \Fc \longrightarrow \Fc\longrightarrow R\Gamma_{\CC \backslash \RR_{-}^n}\Fc$$
Comme $\CC^n\backslash (\CC \backslash \RR_{-}^n)$ est la réunion disjointe :
$$\CC^n\backslash (\CC \backslash \RR_{-}^n) = \coprod_{K\neq \emptyset} \RR_{-}^K\CC\backslash \RR_{-}$$
le triangle distingué devient :
$$\bigoplus_{K\neq \emptyset}(R\Gamma_{\RR^K_{-}\CC \backslash \RR_{-}}\Fc )|_{\CC^L}\longrightarrow \Fc|_{\CC^L}\longrightarrow (R\Gamma_{\CC\backslash \RR_{-}^n}\Fc)|_{\CC^L}$$
Vu la définition des morphismes $n_{K}$, le triplet $\begin{displaystyle}\big(\bigoplus_{K\neq \emptyset}n_{K}, Id, n_{\emptyset}\big)\end{displaystyle}$ est un morphisme de triangles :
$$\xymatrix{
\bigoplus_{K\neq \emptyset}(R\Gamma_{\RR^K_{-}\CC \backslash \RR_{-}}\Fc )|_{\CC^L} \ar[r] \ar[d]_{\bigoplus_{K\neq \emptyset}n_{K}} & \Fc|_{\CC^L} \ar[r] \ar[d]_{Id} & (R\Gamma_{\CC\backslash \RR_{-}^n}\Fc)|_{\CC^L} \ar[d]_{n_{\emptyset}}\\
\bigoplus_{K\neq \emptyset}(R\Gamma_{\RR^K_{-}\CC \backslash \RR_{-})\cap \CC^L}(\Fc|_{\CC^L}) \ar[r] & \Fc|_{\CC^L} \ar[r] & R\Gamma_{(\CC\backslash \RR_{-}^n)\cap \CC^L}(\Fc|_{\CC^L})
}$$
Pour montrer que, pour tout $K$, le morphisme $n_{K}$ est un isomorphisme il suffit donc de montrer que le morphisme $n_{\emptyset}$ est un isomorphisme. Mais d'après le lemme \ref{concentre} les complexes $R\Gamma_{\CC\backslash\RR_{-}^n}\Fc$ et $R\Gamma_{\CC\backslash \RR_{-}^n\cap\CC^L}\Fc$ sont concentrés en degré zéro. Il suffit donc de démontrer que $h^0(n_{\emptyset})$ est un isomorphisme. Notons que :
$$\begin{array}{cccc}
R^0\Gamma_{\CC \backslash \RR_{-}^n}\Fc & \simeq & i_{*}i^{-1}h^0(\Fc)\\
R^0\Gamma_{\CC \backslash \RR_{-}^n\cap \CC^L}(\Fc|_{\CC^L}) & \simeq & i'_{*}i^{'-1}h^0(\Fc)
\end{array}$$
où $i$ et $i'$ sont respectivement les injections de $(\CC\backslash\RR_{-})^n$ dans $\CC^n$ et $(\CC\backslash \RR_{-})^n\cap \CC^L$ dans $\CC^L$. On prend la restriction de $h^0(\Fc)$ car $(\CC\backslash \RR_{-})^n$ est inclus dans $S_{\emptyset}$. \\
Supposons maintenant que $L=\{1, \cdots, l\}$. Soit $x=(x_{1}, \cdots,x_{l}, 1, \cdots ,1)$ un point de $\CC^L$. On a 
$$(\Gamma_{(\CC\backslash\RR_{-})^n}\Fc)_{x}\simeq \varinjlim_{\varepsilon}\Gamma(B_{x}^\varepsilon, i_{*}i^{-1}h^0(\Fc))$$
où la limite parcourt les $n$-uplets réels positifs, et où $B_{x}^\varepsilon$ est le polydisque centré en $x$ et de polyrayon $\varepsilon$.  D'après la définition du foncteur $i$ on a :
$$\begin{array}{cclc}
\varinjlim_{\varepsilon}\Gamma(B_{x}^\varepsilon, i_{*}i^{-1}h^0(\Fc))&\simeq &\varinjlim_{\varepsilon}\Gamma(B_{x}^\varepsilon\cap (\CC \backslash\RR_{-})^n, i^{-1}h^0(\Fc))\\
 & \simeq & \varinjlim_{\varepsilon}\Gamma(B_{\tilde{x}}^{\tilde{\varepsilon}}\cap (\CC \backslash\RR_{-})^l\times B_{1}^{\hat{\varepsilon}}, i^{-1}h^0(\Fc))
\end{array}$$
Ainsi comme $h^0(\Fc)$ est localement constant et comme pour $\varepsilon$ assez petit $B_{\tilde{x}}^{\tilde{\varepsilon}}\cap (\CC \backslash\RR_{-})^l\times B_{1}^{\hat{\varepsilon}}$ est une réunion disjointe d'ouverts contractiles la limite est constante et isomorphe à :
$$\Gamma\big(B_{\tilde{x}}^{\tilde{\varepsilon}}\cap (\CC \backslash\RR_{-})^l\times B_{1}^{\hat{\varepsilon}}, h^0(\Fc)\big).$$
Calculons maintenant la fibre en $x$ du faisceau $R^0\Gamma_{\CC \backslash \RR_{-}^n\cap \CC^L}(\Fc|_{\CC^L})$.
$$\begin{array}{ccl}
(R^0\Gamma_{\CC \backslash \RR_{-}^n\cap \CC^L}(\Fc|_{\CC^L}))_{x}&\simeq &\varinjlim_{}\Gamma(B_{\tilde{x}}^{\tilde{\varepsilon}}\times \{1\}^{n-l}, i'_{*}i^{'-1}h^0(\Fc))\\
&\simeq&\varinjlim\Gamma(B_{\tilde{x}}^{\tilde{\varepsilon}}\cap (\CC\backslash\RR_{-}^l)\times \{1\}^{n-l}, h^0(\Fc))
\end{array}$$
Comme, pour tout $\varepsilon$, $B_{\tilde{x}}^{\tilde{\varepsilon}}\cap (\CC\backslash\RR_{-}^l)\times \{1\}^{n-l}$ est un fermé dans un ouvert paracompact on a l'isomrphisme :
$$\varinjlim\Gamma(B_{\tilde{x}}^{\tilde{\varepsilon}}\cap (\CC\backslash\RR_{-}^l)\times \{1\}^{n-l}, h^0(\Fc)) \simeq \varinjlim\Gamma(B_{\tilde{x}}^{\tilde{\varepsilon}}\cap (\CC\backslash\RR_{-}^l)\times B_{1}^{\hat{\varepsilon}}, h^0(\Fc))$$
Et pour les mêmes raisons que ci-dessus, cette limite est constante et isomorphe à :
$$\Gamma\big(B_{\tilde{x}}^{\tilde{\varepsilon}}\cap (\CC \backslash\RR_{-})^l\times B_{1}^{\hat{\varepsilon}}, h^0(\Fc)\big).$$
Ce qui finit la démonstration.

\cqfd

\cqfd

\chapter{Applications aux arrangements génériques de droites}
Dans ce chapitre nous reprenons la construction faite dans le chapitre précédent pour l'appliquer au cas de $\CC^2$ stratifié par n droites en position générique. Cette généralisation est possible car d'une part la stratification considérée est localement un croisement normal et d'autre part les voisinages tubulaires des strates sont homéomorphes à des produits. Une étude des sections globales du champ ainsi construit nous permet alors de démontrer l'équivalence de la catégorie des faisceaux pervers sur $\CC^2$ relativement à cette stratification et d'une catégorie de représentations du carquois associé à cette stratification.\\

Soit $D_{1}, \cdots, D_{n}$, n droites en position générique dans $\CC^2$. On considère la stratification $\Sigma$ de $\CC^2$ associée définie par : 
$$\begin{displaystyle}\begin{array}{c}
\Sigma_{\emptyset}= \CC^2 \backslash (\bigcap_{i=1}^n D_{i})\\
\Sigma_{\{i\}}= D_{i} \backslash (\bigcup_{j \neq i} D_{i}\cap D_{j})\\
\Sigma_{\{ij\}}= D_{i}\cap D_{j}
\end{array} \end{displaystyle}$$
Considérons $c_{\Sigma}$ le carquois associé à cette stratification. Pour deux droites on obtient un croisement normal. Pour trois droites le carquois $c_{\Sigma}$ est le suivant : 
$$
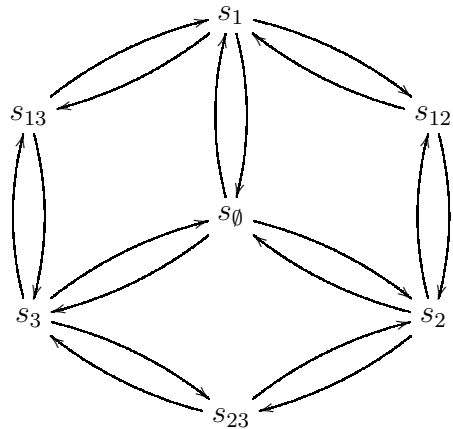
\xymatrix{    & & s_1 \ar@/^/[dll]  \ar@/^/[drr]  \ar@/^/[dd]  \\
                           s_{13}  \ar@/^/[urr]  \ar@/^/[dd]  & &  & & s_{12}   \ar@/^/[ull]  \ar@/^/[dd] \\
                           & & s_{\emptyset} \ar@/^/[uu]  \ar@/^/[dll]  \ar@/^/[drr]   \\
                           s_3  \ar@/^/[uu]  \ar@/^/[urr]  \ar@/^/[drr]   & & & & s_2  \ar@/^/[uu] \ar@/^/[ull]  \ar@/^/[dll]   \\
                           & & s_{23}  \ar@/^/[ull]  \ar@/^/[urr]  \\   }$$ 
On note $\Ic$ l'ensemble $\Big\{\emptyset, \big\{\{i\}\big\}_{i\leq n},\big\{\{i,j\} \big\}_{i<j\leq n}\Big\}$. Pour tout $K \in \Ic$, on note alors $i_{K}$ l'injection de $\Sigma_{K}$ dans $\CC^2$. 
$$i_{K} : \Sigma_{K} \hookrightarrow \CC^2$$                  
Comme pour le croisement normal on a la proposition suivante  :                             
\begin{prop}\label{astrict}
Le champ, $\PPP_{\Sigma}$, des faisceaux pervers relativement à la stratification $\Sigma$ est un champ strictement constructible. 
\end{prop}    
\dem
La démonstration s'appuie sur les mêmes arguments que dans le cas du croisement normal. \\
La restriction à $\Sigma_{ij}$ du champ $\PPP_{\CC^n}$ est constant puisque les strates $\Sigma_{ij}$ sont des points. La restriction à $\Sigma_{\emptyset}$ est elle aussi constante puisque qu'il s'agit du champ des faisceaux localement constants sur $\Sigma_{\emptyset}$. Considérons maintenant la restriction de $\PPP_{\CC^2}$ à $\Sigma_{i}$. On peut supposer sans perte de généralité que $D_{i}= \CC \times \{0\}$ et que la strate $\Sigma_{i}$ est le produit : 
$$\begin{displaystyle}\Sigma_{i} = \CC\backslash (\cup_{j \neq i}p_{j})\times \{0\}\end{displaystyle}$$
où $(p_{j}, 0)=\Sigma_{ij}$.\\
Considérons le voisinage $V_{\varepsilon}$ de $D_{i}$ défini par :
$$V_{\varepsilon}= \big\{(z_{1},z_{2}) \big| ~||z_{2}|| <\varepsilon \big\}$$
\begin{lem}\label{voisinage}
Pour $\varepsilon$ assez petit le voisinage $T_{\varepsilon}= V_{\varepsilon}\backslash(V_{\varepsilon}\cap(\cup_{j\neq i}D_{j})$ est homéomorphe, par un homéophisme stratifié noté $h$, au produit :
$$ \CC\backslash (\cup_{j \neq i}p_{j})\times \CC$$
et la stratification induite est le produit du croisement normal par $ \CC\backslash (\cup_{j \neq i}p_{j})$. 
\end{lem}
On a évidemment l'égalité : 
$$\PPP_{\Sigma} \mid_{\Sigma_{i}}= (\PPP_{\Sigma}\mid_{T_{\varepsilon}}) \mid_{\Sigma_{i}}.$$
Or, comme $h$ est un homéomorphisme stratifié, le champ $\PPP_{\Sigma}\mid_{T_{\varepsilon}}$ est équivalent au champ $h^{-1}(\PPP_{\CC\backslash (\cup_{j \neq i}p_{j}) })$
$$h^{-1}(\PPP_{\CC\backslash (\cup_{j \neq i}p_{j}) }) \simeq \PPP_{\Sigma}\mid_{T_{\varepsilon}}$$
 où $\PPP_{\CC\backslash (\cup_{j \neq i}p_{j}) }$ est le champ des faisceaux pervers sur \linebreak$\CC\backslash (\cup_{j \neq i}p_{j}) \times \CC$ stratifié par la stratification produit du croisement normal et de $\CC\backslash (\cup_{j \neq i}p_{j})$. Mais, d'après le lemme \ref{proj}, $\PPP_{\CC\backslash (\cup_{j \neq i}p_{j}) }$ est lui même équivalent au champ $p^{-1}\PPP_{\CC}$, où $\PPP_{\CC}$ est le champ des faisceaux pervers sur $\CC$, stratifié par le croisement normal, et $p$ est la projection : 
$$p : \CC\backslash (\cup_{i\neq j}p_{j})\times \CC \longrightarrow \CC$$
On a déjà montrer que le champ $\PPP_{\CC}$ était strictement constructible, ainsi d'après le lemme \ref{inv} le champ $p^{-1}(\PPP_{\CC})$ est strictement constructible relativement à la stratification produit. Ce qui démontre que $\PPP_{\Sigma}\mid_{\Sigma_{i}}$ est bien un champ constant. 
\cqfd
On a aussi : 
\begin{lem}
Soit $\Sigma_{K}$ et $\Sigma_{L}$ 2 strates, $K$ et $L$ appartiennent à $\Ic$, telles que $S_{K} \subset \overline{S}_{L}$ (c'est à dire telles que $K\supset L$), on note $k$ et $l$ les cardinaux de respectivement $K$ et $L$. Soit $\CCC_{L}$ un champ constant sur $S_{L}$ de fibre $\Cc$, alors le champ $i_{K}^{-1}i_{L*} \CCC_{L}$ est constant de fibre équivalente à $ Rep(\pi_{1}((\CC^*)^{l-k}), \Cc)$. Cette catégorie est formée d'objets de $\Cc$ munis de $k-l$ automorphismes. 
\end{lem} 
\dem
Localement on est dans le cas du croisement normal, de plus les voisinnages tubulaires des strates sont homéomorphes à des produits de la strate par les boules. On peut ainsi appliquer directement la démonstration du croisement normal.
\cqfd
Comme dans le cas du croisement normal on en déduit la proposition suivante :
\begin{prop}\label{stricte}
La $2$-catégorie des champs strictement constructibles sur $\CC^2$, stratifié par $\Sigma$ est $2$-équivalente à la $2$-catégorie $\SSS_{\Sigma}^{s}$  dont 
\begin{itemize}
\item[$\bullet$] les objets sont donné par : 
\begin{itemize}
\item[-] pour tout $K\subset \Ic$, une catégorie $C_{K}$,
\item[-] pour tout couple $(K,L)$ de parties de $\Ic^2$ tel que $ L \subset K$, un foncteur $F_{lk} : C_{K} \rightarrow Rep(\pi_{1}((\CC^*)^{k-l}), C_{L})$,
\item[-] pour tout triplet $(K,L,M)$ de parties de $\Ic^3$ vérifiant $M\subset L\subset K$ un isomorphisme de foncteurs $\lambda_{KLM}$ : 
$$\shorthandoff{;:!?}\relax
\xymatrix @!0 @C=2,5cm @R=0,6cm {C_{K} \ar[rrr]^{F_{LK}}  \ar[ddddd]_{F_{MK}} &&& Rep(\pi_{1}((\CC^*)^{k-l}),C_{L}) \ar[ddddd]^{Rep(\pi_{1}(\CC^{*k-l}),F_{ML})} \\
~\\
& & \ar@{=>}[ld]_\sim\\
&~\\
\\
Rep(\pi_{1}((\CC^*)^{m-k}),C_{M}) \ar[rrr]_{i^{-1}_k\eta_{lm}~~~~~~~~~} &&& Rep(\pi_{1}((\CC^*)^{k-l}),Rep(\pi_{1}((\CC^*)^{l-m}),C_{M}))}$$
\end{itemize}
tels que, pour tout $k>l>m>p$ les deux morphismes que l'on peut définir entre les foncteurs :
$$ Rep(\pi_{1}(\CC^{*l-k}), Rep(\pi_{1}(\CC^{*m-p}),F_{pm}))\circ Rep(\pi_{1}(\CC^{*l-k}), F_{ml}) \circ F_{lk}$$
$$\text{~et~}$$
$$Rep(\pi_{1}(\CC^{*k-m}), F_{pk})$$
soient égaux.
\item[$\bullet$] pour deux objets $(\{C_{K}\}, \{F_{LK}\}, \{\lambda_{KLM}\})$ et $(\{C'_{K}\}, \{F'_{LK}\}, \{\lambda'_{KLM}\})$ un foncteur est donné par : 
\begin{itemize}
\item[-] pour toute strate $S_{K}$, un foncteur $F_{K} : C_{K}\rightarrow C'_{K}$\\
\item[-] pour tout couple de strates $(S_{K},S_{L})$ tel que $\overline{S}_{L} \supset S_{K}$, un isomorphisme de foncteur $f_{KL}$ :
$$f_{KL} : F'_{LK} \circ F_{K} \buildrel\sim\over\longrightarrow Rep(\pi_{1}((\CC^*)^{k-l}),F_{L}) \circ F_{LK}$$
tels que les deux morphismes que l'on peut définir entre les foncteurs :
$$Rep(\pi_{1}(\CC^{*k-l}), F_{ml}')\circ F_{lk}\circ G_{k}$$
$$\text{et}$$
$$Rep(\pi_{1}(\CC^{*k-m}), G_{m})\circ F_{mk}$$
soit égaux.
\end{itemize}
\item[$\bullet$] Les morphismes entre deux tels foncteurs sont les données pour chaque $S_{K}$ d'un morphisme de foncteurs $\Phi_{k}: G_{k}\rightarrow G'_{k}$ tels que le diagramme suivant commute :
$$\xymatrix{
F'_{kl}\circ G_{k} \ar[d]_{Id\bullet\Phi_{k}} \ar[r]^{g_{kl}~~~~~} & Rep(\pi_{1}(\CC^{*k-l}),G_{l})\circ F_{lk} \ar[d]^{Rep(\pi_{1}(\CC^{*k-l}), \Phi_{k})\bullet Id}\\
F'_{kl}\circ G'_{k}\ar[r]_{g'_{kl}~~~~~} & Rep(\pi_{1}(\CC^{*k-l}), G'_{l})\circ F_{lk}
}$$
\end{itemize}
\end{prop}
Comme dans le cas du croisement normal, on définit un champ $\CCC_{\Sigma}$ équivalent au champ $\PPP_{\Sigma}$. Les données d'un objet de $\SSS_{\Sigma}^s$ étant des données locales, et la stratification étant localement le croisement normal, on peut réutiliser les données définies dans le chapitre précédent. On considère alors $\Cc_{k}$, $F_{LK}$ pour $(K,L) \in \Ic^2$ les catégories et foncteurs définis dans le chapitre précédent. 
\begin{Def}
On note $\CCC_{\Sigma}$ le champ image par $Q^s_{\Sigma}$ de   l'objet de $\SSS_{\Sigma}^s$ défini par la donnée  :
\begin{itemize}
\item pour toute strate $\Sigma_{K}$ de $\CC^2$, de la catégorie $\Cc_{k}$, 
\item  pour tout couple de strates, $(\Sigma_{K},\Sigma_{L})$, telles que $\Sigma_{L} \subset \overline{\Sigma}_{K}$, du foncteur $F_{LK}$,
\item et pour tout triplet de strates $(\Sigma_{J}, \Sigma_{K}, \Sigma_{L})$ tel que $J \subsetneq K \subsetneq L$, du morphisme de foncteur $Id$.
\end{itemize}
\end{Def}
\begin{thm}\label{aeq}
Le champ $\CCC_{\Sigma}$ est équivalent au champ $\PPP_{\Sigma}$.
\end{thm}
\dem
Encore une fois, la démonstration de ce théorème étant locale on peut tout à fait appliquer la preuve du cas à croisement normal.
\cqfd
Considérons la catégorie $R(c_{\Sigma})$ des représentations du carquois $c_{\Sigma}$. Les objets de cette catégorie sont des familles 
$$\Big\{ \{E_{K}\}_{K \in \Ic}, \{u_{Kl}, v_{Kl}\}_{(K, K\cup l)\in \Ic^2}\Big\}$$
où $E_{K}$ est un espace vectoriel, et $u_{Kl}$ et $v_{Kl}$ sont des applications linéaires :
$$\begin{array}{cccc}
u_{Kl} : & E_{K} & \longrightarrow & E_{K\cup l}\\
v_{Kl} : & E_{K\cup l} & \longrightarrow & E_{K}\\
\end{array}$$
\begin{Def}
On note $\Cc_{\Sigma}$ la sous catégorie pleine de la catégorie des représentations du carquois $R(c_{\Sigma})$ dont les objets 
$$\Big\{ \{E_{K}\}_{K \in \Ic}, \{u_{Kl}, v_{Kl}\}\Big\}$$
vérifient les conditions : 
\begin{itemize}
\item[(i)] pour tout couple $(K\cup\{p\},K)$ appartenant à $\Ic^2$  l'application linéaire :
$$M_{Ip}= v_{Ip} u_{Ip} + Id $$
soit inversible,\\
\item[(ii)] pour tout couple $(i,j) \in \{ 1, \ldots, n\}^2$ on ait :
$$u_{\{i\}j}u_{\emptyset i}=u_{\{j\}i}u_{\emptyset j}, ~ v_{\emptyset i}v_{\{i\}j}=v_{\emptyset j}v_{\{j\} i}, ~v_{\{i\}j}u_{\{j\}i}=u_{\emptyset i}v_{\emptyset j},$$
\item[(iii)] pour $i \in \{1, \ldots,n\}$ les applications linéaires $M_{\emptyset i}$ commutent.
\end{itemize}
\end{Def}
L'étude de la catégorie des sections globales de ce champ nous permet de démontrer le théorème suivant :
\begin{thm}
La catégorie $\Pc erv_{\Sigma}$ des faisceaux pervers sur $\CC^2$ stratifié par $\Sigma$ est équivalente à la catégorie $\Cc_{\Sigma}$.
\end{thm}
\dem
 D'après le théorème \ref{aeq}, il suffit de démontrer que la catégorie des sections globales de $\CCC_{\Sigma}$ est équivalente à $\Cc_{\Sigma}$. Les objets de la catégorie des sections globales de $\CCC_{\Sigma}$ sont les familles 
 $$\Big\{ \{P_{K}\}_{K \in \Ic}, \{\omega_{KL}\}_{(K,L) \in \Ic^2}\Big\},$$
 où les $P_{K}$ sont des objets de \linebreak$Rep( \pi_{1}(\Sigma_{K}), \Cc_{k})$, et les $\omega_{KL}$ sont des isomorphismes :
 $$\omega_{KL} : i_{K*}F_{LK}(P_{K})\buildrel\sim\over\longrightarrow i_{K*}\eta_{KL}(P_{L})$$
 qui vérifient des relations de compatibilité. Pour pouvoir expliciter cette catégorie nous avons donc besoin de connaître la topologie de chaque strate. On doit aussi expliciter les  sections globales des foncteurs naturels $i_{K*}\eta_{KL}$. \\
 Les strates $\Sigma_{ij}$ sont des points. Ainsi $Rep(\pi_{1}(\Sigma_{ij}), \Cc_{2})=\Cc_{2}$.\\
Les strates $\Sigma_{i}$ sont isomorphes à $\CC$ privé de $n-1$ points. Pour chacune d'entre elle on fixe un point base et pour tout $j\neq i$ un lacet $\gamma_{ij}$ contournant $\Sigma_{ij}$. La famille $(\gamma_{ij})_{j \neq i}$ est un système de générateurs du groupe fondamental de $\Sigma_{i}$. Les objets de la catégorie $Rep(\pi_{1}(\Sigma_{i}), \Cc_{1})$ sont donc des objets de $\Cc_{2}$ munis de $n-1$ automorphismes. Si $P$ est un objet de $\Cc_{1}$ :
$$P=\xymatrix{
    E  \ar@/^/[r]^{u} & 
    F  \ar@/^/[l]^{v} ,}$$  
un automorphisme de $P$ est donné par deux endomorphismes \linebreak$M : E\rightarrow E$ et $N :F \rightarrow F$   qui commutent aux applications $u$ et $v$. On note $M_{ij}$ et $N_{ij}$ les endomorphismes associés aux lacets $\gamma_{ij}$.\\
Pour la strate $\Sigma_{\emptyset}$ on a le théorème suivant : 
\begin{thm}
Soient $\Ac$ un arrangement d'hyperplans dans $\CC^2$, $p$ un point base du complémentaire de $\Ac$ dans $\CC^2$ et $\{\gamma_{i}\}_{i\leq n}$ une famille de lacets contournant chaque droite de l'arrangement, alors le groupe fondamental du complémentaire de $\Ac$ est le groupe commutatif engendré par la famille $\{\gamma_{i}\}_{i\leq n}$.
\end{thm}
\dem
Voir \cite{Zar} chapitre 3 section 2.
\cqfd
Fixons donc un point base et pour tout $i\leq n$ un lacet $\gamma_{i}$ contournant la droite $D_{i}$. Les objets de $Rep(\pi_{1}(\Sigma_{\emptyset}), Ev)$ sont donc des espaces vectoriels munis de $n$ endomorphismes qui commutent. Comme ci-dessus, on note $M_{i}$ les endomorphismes associés aux lacets $\gamma_{i}$.\\\\
Considérons maintenant les foncteurs naturels d'adjonctions.
\begin{lem}
Pour tout $k\leq n $, la section globale du foncteur d'adjonction $ \eta_{j\emptyset}$ :
$$\begin{array}{ccccccc}
\Gamma( \CC^2,i_{\emptyset*}\eta_{j\emptyset} ):&\Gamma (\CC^2,i_{\emptyset*}\CCC_{\emptyset}) & \longrightarrow & \Gamma(\CC^2,i_{j*}i_{j}^{-1}i_{\emptyset*}\CCC_{\emptyset})\\ 
\end{array}$$
est naturellement isommorphe, pour $\varepsilon$ assez petit, au foncteur suivant :
$$\begin{array}{cccc}
\Gamma(\Sigma_{\emptyset}, \CCC_{\emptyset}) & \longrightarrow & \Gamma(T_{\varepsilon}\cap \Sigma_{\emptyset}, \CCC_{\emptyset})\\
\Big\{E, \{M_{\emptyset i}\}_{i \leq n}\Big\} & \longmapsto & \Big\{E, \{M_{\emptyset i}\}_{i \leq n}\Big\}
\end{array}$$
où $T_{\varepsilon}$ est un voisinage tubulaire de $\Sigma_{j}$.\\
Pour $j<k\leq n$, la section globale du foncteur naturel $i_{j*}\eta_{\{jk\}\{j\}}$ est naturellement isomorphe au suivant :
$$\begin{array}{ccccccc}
 \Gamma( \CC^2,i_{j*}\CCC_{1}) & \longrightarrow &\Gamma( \CC^2, i_{jk*}i_{jk}^{-1}i_{j*}\CCC_{1})\\ 
\Big\{ P,\{M_{jl}, N_{jl}\}_{l\neq j }\Big\} & \longmapsto & \Big\{ P,M_{jk}, N_{jk}\Big\} 
\end{array}$$
\end{lem}
\dem
Montrons la première assertion. On peut supposer sans perte de généralité que $D_{j}=\CC \times \{0\}$. Par définition des $2$-foncteurs $i_{\emptyset*}$ et $i_{j*}$ le foncteur naturel $\eta_{j \emptyset}$ est le foncteur naturel :
$$\eta_{j \emptyset}: \Gamma(\Sigma_{\emptyset}, \CCC_{\emptyset}) \longrightarrow \Gamma(\Sigma_{j}, i_{\emptyset*}\CCC_{\emptyset})$$
Ce morphisme se factorise par la $2$-limite :
$$\begin{displaystyle}
2\varinjlim_{U\supset \Sigma_{j}}\Gamma(U, i_{\emptyset*}\CCC_{\emptyset})
\end{displaystyle}$$
Or d'après la proposition \ref{paracompact}, le foncteur naturel :
$$\begin{displaystyle}
2\varinjlim_{U\supset \Sigma_{j}}\Gamma(U, i_{\emptyset*}\CCC_{\emptyset})\longrightarrow \Gamma(\Sigma_{j}, i_{\emptyset}\CCC_{\emptyset})
\end{displaystyle}$$
est une équivalence. De plus les voisinages tubulaires de $\Sigma_{j}$ forment une base de voisinage de $\Sigma_{j}$ et, pour tout $\varepsilon$, les voisinages $T_{\varepsilon}$  sont homotopiquement équivalents, par une équivalence stratifiée. Ainsi la $2$-limite est constante et, pour $\varepsilon$ assez petit, le foncteur naturel est  une équivalence :
$$\Gamma(T_{\varepsilon}, i_{\emptyset}\CCC_{\emptyset}) \longrightarrow 2\varinjlim_{U\supset \Sigma_{j}}\Gamma(U, i_{\emptyset*}\CCC_{\emptyset})$$ 
et par définition du foncteur $i_{\emptyset}$, on a l'égalité :
$$\Gamma(T_{\varepsilon}, i_{\emptyset*}\CCC_{\emptyset}) =\Gamma(T_{\varepsilon}\cap \Sigma_{\emptyset}, \CCC_{\emptyset})$$
Ainsi, nous avons l'isomorphisme suivant :
$$\shorthandoff{;:!?}\relax
\xymatrix @!0 @C=2,5cm @R=0,6cm {
\Gamma(\Sigma_{\emptyset},\CCC_{\emptyset})\ar[rrr]  &&& \Gamma(\Sigma_{j}, i_{\emptyset*}\CCC_{\emptyset}) \\
~\\
& & \\
&\ar@{=>}[ru]^\sim\\
\\
\Gamma(\Sigma_{\emptyset},\CCC_{\emptyset})\ar[uuuuu]^{Id} \ar[rrr] &&& \begin{displaystyle}\Gamma(T_{\varepsilon}\cap \Sigma_{\emptyset}, \CCC_{\emptyset})\end{displaystyle} \ar[uuuuu]
}$$
où le foncteur du bas est le foncteur de restriction.
Or d'après le lemme \ref{voisinage}, $T_{\varepsilon}\cap \Sigma_{\emptyset}$ est homéomorphe au produit $\CC\backslash (\cup p_{k})\times \CC^*$, ainsi les objets de $ \Gamma(T_{\varepsilon}\cap \Sigma_{\emptyset}, \CCC_{\emptyset})$ sont des espaces vectoriels munis de $n-1$ automorphismes et le foncteur de restriction est bien le foncteur décrit dans le lemme.\\
Considérons maintenant le foncteur suivant :
$$\Gamma( \Sigma_{j},\CCC_{1})  \longrightarrow \Gamma( \CC^2, i_{jk*}i_{jk}^{-1}i_{j*}\CCC_{1})$$
Par un raisonnement analogue au précédent on montre qu'il existe un isomorphisme :
$$\shorthandoff{;:!?}\relax
\xymatrix @!0 @C=2,5cm @R=0,6cm {
\Gamma( \Sigma_{j},\CCC_{1})\ar[rrr]  &&&\Gamma( \Sigma_{jk}, i_{jk}^{-1}i_{j*}\CCC_{1})  \\
~\\
& & \\
&\ar@{=>}[ru]^\sim\\
\\
\Gamma( \CC^2,i_{j*}\CCC_{1})\ar[uuuuu]^{Id} \ar[rrr] &&& \begin{displaystyle}\Gamma( T_{\varepsilon}^{jk} \cap \Sigma_{j}, \CCC_{1}) \end{displaystyle} \ar[uuuuu]
}$$
où le foncteur du bas est la restriction. On retrouve bien le foncteur annoncé.
\cqfd
Un objet de $\Cc_{\Sigma}$ est donc donné par : 
\begin{itemize}
\item un espace vectoriel $E$ muni de $n$ endomorphismes notés $M_{1}, \ldots, M_{n}$ qui commutent,
\item pour tout $i\leq n$, un objet de $\Cc_{1}$ :
$$P_{i}=\xymatrix{
E_{i} \ar@/^/[r]^{u_{i}} & F_{i} \ar@/^/[l]^{v_{i}}
}$$
avec $n-1$ automorphismes, c'est-à-dire $n-1$ couples d'automorphismes notés $(M_{ij}, N_{ij})_{j \neq i}$ : 
$$\begin{array}{cccccc}
M_{ij} : & E_{i} & \buildrel\sim\over\longrightarrow & E_{i}\\
N_{ij} : & F_{i} & \buildrel\sim\over\longrightarrow & F_{i},
\end{array}$$
rappelons que $v_{i} \circ u_{i} +Id=M_{ii}$ est inversible,
\item pour tout ensemble $\{i,j\} \in \Ic$ un objet de $\Cc_{2}$ : 
$$\shorthandoff{;:!?}\relax
\xymatrix @!0 @C=2cm @R=2cm {
F_{ij}^i \ar@/^/[r]^{w_{iji}} \ar@/^/[d]^{v_{iji}} & G_{ij} \ar@/^/[d]^{t_{ijj}} \ar@/^/[l]^{t_{iji}} \\
E_{ij}\ar@/^/[r]^{u_{ijj}} \ar@/^/[u]^{u_{iji}} & F_{ij}^j \ar@/^/[u]^{w_{ijj}} \ar@/^/[l]^{v_{ijj}}
}$$
rappelons que selon la définition de la catégorie $\Cc_{2}$ les applications linéaires suivantes : 
$$v_{ijk}\circ u_{ijk} +Id=M_{ijk} ~~~~t_{ijk} \circ w_{ijk} +Id= N_{ijk}$$
sont inversibles.
\item pour tout $i\leq n $ une application linéaire inversible $\delta_{i} : E \buildrel\sim\over\rightarrow E_{i}$ vérifiant pour tout $j \leq n$ :
\begin{equation}\label{delta} \delta_{i}^{-1}M_{ij}\delta_{i}= M_{i}^{-1}\end{equation}
\item pour tout couple $\{i,j \} \in \Ic $, trois applications linéaires inversibles :
$$\begin{array}{cccc}
\alpha_{ij} : &E_{i} & \buildrel\sim\over\longrightarrow & E_{ij}\\
\beta_{ij} : &F_{i} & \buildrel\sim\over\longrightarrow& F^{i}_{ij}\\
\gamma_{ij} : &E & \buildrel\sim\over\longrightarrow & E_{ij}\\
\end{array}$$
qui vérifient les conditions suivantes pour $k=i$ ou $k=j$ : 
\begin{equation}\label{alpha}\begin{array}{ccc}
\alpha_{ij}^{-1}M_{ijk}\alpha_{ij}= M_{ik}, &
\beta_{ij}^{-1}N_{ijk}\beta_{ij}= N_{ik}, & 
\gamma_{ij}^{-1}M_{ijk}\gamma_{ij}= M_{k}, 
\end{array}\end{equation}
les relations de commutations des isomorphismes de foncteurs donnent la relation : 
$$\alpha_{ij}\circ \delta_{i}= \gamma_{ij}.$$
\end{itemize}
On note $\Delta$ le foncteur de la catégorie des sections globales dans la catégorie $\Cc_{\Sigma}$ :
$$\Delta : \Gamma(\CC^2, \CC_{\Sigma} )\longrightarrow \Cc_{\Sigma}$$
 qui à un objet tel que nous l'avons décrit associe la famille  :
$$\Big\{E, \{E_{\{i\}}\}_{i\leq n}, \{G_{\{i,j\}}\}_{\{i,j\} \in \Ic}, \{u'_{i}, v'_{i}\}_{i \leq n}, \{u_{ij}, v_{ij}\} \Big\}$$
où les applications linéaires $u'_{i}$, $v'_{i}$, $u_{ij}$ et $v_{ij}$ sont définies par : 
$$\begin{array}{llccc}
u'_{i}=u_{i} \circ \delta_{i} & v'_{i}= \delta_{i}^{-1} \circ v_{i} \\
u_{ij} = w_{iji} \circ \beta_{ij} & v_{ij}=\beta_{ij}^{-1}\circ t_{iji} & \text{si $i<j$}\\
u_{ij} = w_{ijj} \circ \beta_{ij} & v_{ij}=\beta_{ij}^{-1}\circ t_{ijj} & \text{si $i>j$}\\
\end{array}$$
Cette famille est bien un objet de $\Cc_{\Sigma}$. En effet les conditions $(i)$ et $(ii)$ de la définition de $\Cc_{\Sigma}$ sont vérifiées par définitions des catégories $\Cc_{1}$ et $\Cc_{2}$. Quant à la conditions $(iii)$, elle est donnée par les relations (\ref{delta}) et (\ref{alpha}) et par la commutation des $M_{i}$.\\
Définissons un foncteur $\Lambda$ quasi-inverse de $\Delta$. Ce foncteur associe à un objet $\Big\{\{E_{K}\}_{K\in \Ic}, \{u_{Kl}, v_{Kl}\}\Big\}$ de $\Cc_{\Sigma}$ l'objet défini par : 
\begin{itemize}
\item pour la strate $\Sigma_{\emptyset}$, l'espace vectoriel $E$ et des $n$ automorphismes $M_{\emptyset i}= v_{i}u_{i}+ Id$,
\item pour la strate $\Sigma_{i}$ de l'objet de $\Cc_{1}$ : 
$$\xymatrix{
E \ar@/^/[r]^{u_{i}} &F_{i} \ar@/^/[l]^{v_{i}}\\
}$$
munis pour des $n-1$ endomorphismes $M_{\emptyset j}$ de $E$, avec $i\leq n$, et des $n-1$ endomorphismes de $F_{i}$, $M_{\{i\}j}=v_{ij}u_{ij}+Id$, avec $i \leq n$.
\item pour la strates $\Sigma_{ij}$,  l'objet de $\Cc_{2}$ :
$$\shorthandoff{;:!?}\relax
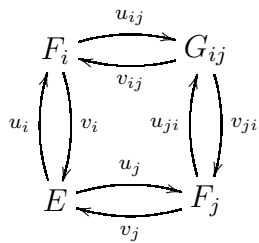
\xymatrix @!0 @C=2cm @R=2cm {
F_{i} \ar@/^/[r]^{u_{ij}} \ar@/^/[d]^{v_{i}} & G_{ij} \ar@/^/[d]^{v_{ji}} \ar@/^/[l]^{v_{ij}} \\
E \ar@/^/[r]^{u_{j}} \ar@/^/[u]^{u_{i}} & F_{j} \ar@/^/[u]^{u_{ji}} \ar@/^/[l]^{v_{j}}
}$$
\item et enfin pour toutes les strates par l'identité.
\end{itemize}
La composée $\Delta\circ \Lambda$ est l'identité de la catégorie $\Cc_{\Sigma}$. \\
Si maintenant $P$ est un objet de la catégorie $\Gamma(\CC^2,\CCC_{\Sigma})$ avec les mêmes notations que ci-dessus, alors les isomorphismes $\delta_{i}$, $\alpha_{ij}$, $\beta_{ij}$ et $\gamma_{ij}$ forment un isomorphisme entre $P$ et $\Lambda \Delta(P)$.
\cqfd

\chapter{Applications aux variétés toriques lisses}
Soit $\Sigma$ un éventail régulier et $X_{\Sigma}$ la variété torique associée. Le but de ce chapitre est de définir un champ sur $X_{\Sigma}$ strictement constructible relativement à la stratification donnée par l'action du tore équivalent au champ des faisceaux pervers sur $X_{\Sigma}$ relativement à la même stratification. Cette construction se basera sur le chapitre précédent et sur le lemme \ref{recouvrement}. En étudiant les sections globales de ce champ nous donnerons alors une sous-catégorie pleine de représentations de carquois équivalente à la catégorie des faisceaux pervers sur $X_{\Sigma}$ stratifiée par l'action du tore.
\section{Rappel sur les variétés toriques}
Dans ce paragraphe nous donnons quelques définitions concernant les variétés toriques. Pour les démonstrations le lecteur pourra se reporter à \cite{Od} \cite{Fu} ou à \cite{Ew}. \\
\subsection{Cônes et variété torique affine}~\\
Soient $V$ un $\RR$-espace vectoriel et $N\subset V$ un réseau, on note $V^*$ le dual de $V$ et $M=Hom(N, \ZZ)$ le dual de $N$, en pratique on prendra toujours $N=\ZZ^n$.
\begin{Def}
Un cône convexe polyédral $\sigma$ est un ensemble :
$$\{\lambda_{1}v_{1}+\cdots+\lambda_{s}v_{s} \in V \mid \lambda_{i}> 0\}$$
où la famille $v_{i}$ est une famille finie de vecteur de $V$, cette famille est une famille  génératrice du cône $\sigma$ et on note :
$$\sigma=<v_{1}, \cdots ,v_{s}>$$
\end{Def}
La dimension de $\sigma$ est la dimension de l'espace affine engendré par le cône $\sigma$.  On note $\check{\sigma}$ le dual du cône $\sigma$ : 
$$\check{\sigma}=\{ u \in V^*\mid <u,v> \geq 0, \forall v \in \sigma\}$$
\begin{Def}
Une face d'un cône convexe polyédral $\sigma$ est soit $\emptyset$, soit $\sigma$, soit un sous-ensemble de $\sigma$ tel qu'il existe un hyperplan $H$ vérifiant :
$$\sigma \cap H \neq \emptyset \text{~et~} \sigma \subset H^+ \text{~ou~} \sigma \subset H^-$$
\end{Def}
On appelle sommet d'un cône une face de dimension zéro, arrête d'un cône une face de dimension un et facette une face de codimension un.
\begin{Def}
On dit qu'un cône convexe polyédral est strictement convexe s'il admet $\{0\}$ comme sommet.
\end{Def}
On a les propriétés suivantes 
\begin{itemize}
\item Une face d'un cône convexe polyédral est aussi un cône convexe polyédral.
\item Toute intersection de faces est une face.
\item Le dual d'un cône convexe polyédral est un cône convexe polyédral. 
\end{itemize}
\begin{Def}
Soit $\sigma=<v_{1}, \ldots,v_{n}>$ un cône strictement convexe polyédral de dimension $n$. Un vecteur $g$ est dit entrant normal à la facette $<v_{1}, \ldots, v_{i-1}, v_{i+1}, \ldots, v_{n}>$ s'il vérifie les conditions :
\begin{itemize}
\item $<g,v_{j}>=0 $ pour $j \neq i$,
\item et $<g,v_{j}> >0$.
\end{itemize}
\end{Def}
Dans tout ce qui suit, si $\sigma=<v_{1}, \ldots,v_{n}>$, un vecteur entrant normal à la facette $<v_{1}, \ldots, v_{i-1}, v_{i+1}, \ldots, v_{n}>$ est noté $g_{i}$.
\begin{prop}\label{cone}
Soient $\sigma=<v_{1}, \ldots,v_{k}>$ un cône polyédral, strictement convexe de dimension $k$, $v_{k+1}, \ldots, v_{n}$ des vecteurs tels que la famille $\{v_{i}\}_{i\leq n}$ soit une base de $\RR^n$ et $g_{1}, \ldots, g_{n}$ des vecteurs entrants normaux aux facettes de $<v_{1}, \ldots,v_{n}>$, alors 
$$\check{\sigma}=<g_{1}, \ldots, g_{k}, \pm g_{k+1},  \ldots \pm g_{n}>$$
\end{prop}
\dem
Cette proposition et sa démonstration sont une généralisation du théorème $2.1$ chapitre $5$ de \cite{Ew} et de sa démonstration.\\
On note $\tilde{\sigma}=<g_{1}, \ldots, g_{k}, \pm g_{k+1},  \ldots \pm g_{n}>$, soit $x \in \tilde{\sigma}$.
$$x=\sum_{i=1}^n\lambda_{i}g_{i}$$ 
avec $\lambda_{i}\geq 0$ pour $i\leq k$. Alors par définition des vecteurs entrants normaux aux facettes on a, pour $i\leq n$ :
$$<x,g_{i}>=\lambda_{i}<g_{i}, v_{i}> \geq 0.$$
Ainsi $\tilde{\sigma} \subset \check{\sigma}$.\\
Soit maintenant $x \in \check{\sigma}$ tel que $x\notin \tilde{\sigma}$. 
$$x=\sum_{i=1}^n \lambda_{i}g_{i}$$
Comme $x \notin \tilde{\sigma}$ il existe $1\leq j\leq k$, tel que $\lambda_{j}<0$, mais dans ce cas 
$$<x,v_{j}>=\lambda_{j}<g_{j},v_{j}> <0$$
d'où la contradiction.
\cqfd
\begin{Def}
On dira qu'un cône $\sigma$ est rationnel si on peut prendre des générateurs de $\sigma$ appartenant à $N$.
\end{Def}
Si $\sigma$ est rationnel, son dual $\check{\sigma}$ l'est aussi. 
\begin{prop}[Lemme de Gordon]
Soit $\sigma$ un cône convexe rationnel, alors le semi-groupe $ \check{\sigma} \cap M$ est engendré par un nombre fini d'éléments.
\end{prop}
\noindent
Pour tout $a=(a_{1}, \cdots, a_{n}) \in \ZZ^n$, on note $$X^a=X_{1}^{a_{1}}\ldots X_{n}^{a_{n}}$$
Soit alors $\CC [\sigma]$ la $\CC-$algèbre définie par : 
$$\CC[\sigma]= \left\{\sum \lambda_{a}X^{a}\mid \lambda_{a}= 0 \text{~si~} a \notin \check{\sigma}\cap M \right\} $$
On note $X_{\sigma}$ la variété affine : 
$$X_{\sigma}=Spec(\CC[\sigma])$$
Choisir un système de générateurs du monoïde $\check{\sigma}\cap \ZZ^n$ permet de plonger la variété $X_{\sigma}$ dans $\CC^k$. En effet si $v_{1}, \ldots, v_{k}$ est une famille génératrices de $\check{\sigma} \cap \ZZ^n$ alors $X^{v_{1}}, \ldots, X^{v_{n}}$ est un système générateur de la $\CC$-algèbre $\CC[\sigma]$. On  a alors l'identification :
$$\CC[\sigma]\simeq \CC[\xi_{1}, \cdots, \xi_{k}]/ I$$
où $I$ est le noyau du morphisme d'algèbre : 
$$\begin{array}{ccc}
 \CC[\xi_{1}, \ldots, \xi_{k}]& \longrightarrow &\CC[X^{v_{1}}, \ldots, X^{v_{k}}]\\ 
\xi_{i} & \longmapsto & X^{v_{i}}
\end{array}$$
\begin{prop}
La variété $X_{\sigma}$ est lisse si et seulement s'il existe un système de générateurs de $\sigma$ que l'on peut compléter en une base de $\ZZ^n$.
\end{prop}
\noindent
\begin{bf} Exemples \end{bf}\\
\begin{itemize}
\item Soit $\sigma=\{0\}$, on a alors $\check{\sigma}= \RR^n$ d'où 
$$\CC[\sigma]= \CC[X_{1}, \ldots, X_{n}, X_{1}^{-1}, \ldots ,X_{n}^{-1}]$$
Ainsi $X_{\sigma}$ est isomorphe à $(\CC^*)^n$. \\
\item Soit $X=\RR^n$ et $\sigma= <e_{1}, e_{2}>$ où $\{e_{1}, e_{2}\}$ est la base canonique de $\RR^2$, alors $\check{\sigma}=\sigma$ et 
$$\CC[\sigma]= \CC[X_{1}, X_{2}]$$
ainsi $X_{\sigma}= \CC^2$.\\
\item Soit $V= \RR^2$ et $\sigma= <(1,0), (1,2)>$, le semi-groupe $\sigma \cap \ZZ^n$ est engendré par $(1,0),(1,1),(1,2)$. La variété $X_{\sigma}$ est alors :
$$X_{\sigma}= \{(u,v,w) \in \CC^3 | ~uw=v^2\}$$
\end{itemize}
Un morphisme de semi-groupe $\check{\sigma} \cap \ZZ^n \rightarrow \check{\sigma'} \cap \ZZ^n$ induit un morphisme de $\CC$-algèbre $\CC[\sigma] \rightarrow \CC[\sigma']$ et ainsi un morphisme de variétés $X_{\sigma'} \rightarrow X_{\sigma}$. En particulier si $\tau$ est une face de $\sigma$ alors $\check{\sigma} \cap \ZZ^n$ est un sous-semi-groupe de $\check{\tau} \cap \ZZ^n$, ce qui induit un morphisme $X_{\tau}\rightarrow X_{\sigma}$. 
\begin{lem}
Si $\tau$ est une face de $\sigma$ l'application induite $X_{\tau}\rightarrow X_{\sigma}$ est une injection et identifie $X_{\tau}$ comme un ouvert de zarisky de $X_{\sigma}$.
\end{lem}

\subsection{Eventail et variétés toriques }
\begin{Def}
Un éventail $\Delta$ est un ensemble fini de cônes strictement convexes, polyédraux et rationnels vérifiant les propriétés suivantes : 
\begin{itemize}
\item[(i)] toute face d'un cône appartient à $\Delta$,
\item[(ii)] l'intersection de deux cônes de $\Delta$ est une face de chacun des deux cônes. 
\end{itemize}
\end{Def}
\begin{bf} Remarque \end{bf} L'ensemble constitué de toutes les faces d'un cône est un éventail. 
\begin{Def}
Soit $\Delta$ un éventail, on note $X_{\Delta}$ la variété algébrique formée de l'union disjointe des variétés $X_{\check{\sigma}}$ pour tout $\sigma \in \Delta$ recollées comme suit : pour deux cônes $\sigma$ et $\tau$ de $\Delta$ on recolle $X_{\check{\sigma}}$ et $X_{\check{\tau}}$ le long de l'ouvert de Zarisky $X_{\check{\sigma \cap \tau}}$.
\end{Def}
La compatibilité des morphismes de recollement est donnée par la correspondance de l'inclusion entre les cônes et les variétés associées. \\
\section{Description du champ  $\PPP_{\Delta}$ et de la catégorie $\Pc erv_{\Delta}$}
Dans le paragraphe qui suit on note $\Delta$ un éventail,  $v_{1}, \cdots, v_{k}$ les vecteurs primitifs qui engendrent les cônes de dimension un appartenant à $\Delta$, et $\Ic$ l'ensemble des parties $I$ de $\{1, \cdots, k\}$ telles que $\sigma_{I}=pos(\{v_{i}\}_{i \in I})$ soit un cône de $\Delta$. On note $X_{I}$ la variété associée au cône $\sigma_{I}$ et $\PPP_{I}$ le champ des faisceaux pervers sur $X_{I}$ relativement à la stratification donnée par l'action du tore.\\
Pour tout couple $(I,J)$ d'ensembles de $\Ic$, on note $h_{IJ}$ l'isomorphisme de recollement des variétés $X_{I}$ et $X_{J}$ :
$$h_{IJ} : X_{I \cap J} \longrightarrow X_{J \cap I}$$
On suppose de plus que $\Delta$ un éventail régulier de $\RR^n$. 
\begin{Def}
Soient $\sigma_{I}$ un cône de $\Delta$, il est maximal dans $\Delta$ si le seul cône  appartenant à $\Delta$ qui le contient est lui-même. On dit alors que la partie $I$ de $\Ic$ est maximale et on note $\Mm \Ic$ l'ensemble des parties maximales de $\Ic$.
\end{Def}
Pour le reste du paragraphe on se fixe, pour tout cône $\sigma_{I}$ maximal de $\Delta$ une base de $\ZZ^n$ qui contient les vecteurs $\{v_{i}\}_{i \in I}$ .
\begin{Def}
Soit $\Delta$ un éventail de $\CC^n$, on note $c_{\Delta}$ le carquois donné par :
\begin{itemize}
\item[$\bullet$] pour tout cône $\sigma_{I}$ de $\Delta$, un sommet $s_{I}$ et $n-j$ flèches de source et de but $s_{I}$, $j$ étant le cardinal maximal des ensembles maximaux de $\Ic$ contenant $I$,
\item[$\bullet$] pour tout couple $(\sigma_{I}, \sigma_{I'})$ de cônes de $\Delta$ tels que $\sigma_{I'}$ soit une face de codimension un de $\sigma_{I}$, deux flèches l'une de source $s_{I}$ et de but $s_{I'}$ et l'autre de source $s_{I'}$ et de but $s_{I}$.
\end{itemize}
\end{Def}
\textbf{Exemples}
\begin{itemize}
\item Si $e_{1}, \ldots, e_{n}$ est une base de $\ZZ^n$ et $\Delta$ est l'éventail associé au cône $<e_{1}, \ldots, e_{n}>$, la variété torique associée est l'espace $\CC^n$ et le carquois associé est un hypercube de dimension $n$, dont les arrêtes contiennent deux flèches en sens inverse. 
\item Si maintenant $\Delta$ est l'éventail associé au cône $<e_{1}, \ldots e_{l}>$\linebreak avec $l\leq n$, la variété torique associée est isomorphe au produit $\CC^{l}\times \CC^{*n-l}$ et le carquois associé est un hypercube de dimension $l$ dont les arrêtes sont formées de deux flèches de sens inverse muni à chaque sommet de $n-l$ flèches ayant pour source et pour but le sommet.
\item Si $\Delta \in \CC^2$ est la réunion des faces des cônes $<e_{1}, e_{2}>$ et $<-e_{1}-e_{2}>$ où $(e_{1},e_{2})$ est la base canonique de $\CC^2$. Le carquois associé est le suivant :
$$
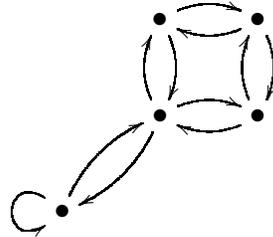
\xymatrix{
& \bullet \ar@/^/[r] \ar@/^/[d] & \bullet \ar@/^/[l] \ar@/^/[d] \\
& \bullet \ar@/^/[u] \ar@/^/[r] \ar@/^/[dl] & \bullet \ar@/^/[u] \ar@/^/[l]\\
\bullet \ar@/^/[ur] \ar@(ul,dl)[]
}$$
\end{itemize}
\begin{Def}\label{CcDelta}
On note $\Cc_{\Delta}$ (ou $\Cc_{X_{\Delta}}$) la sous-catégorie pleine de $R (c_{\Delta})$ formée des objets tels que  :
\begin{itemize}
\item[(i)] pour tout couple $(K\cup\{p\}, K)$ appartenant à $\Ic^2$  l'application linéaire :
$$M_{Kp}= v_{Kp} u_{Kp} + Id $$
soit inversible,
\item[(ii)] pour tout quadruplets $(K, K\cup \{p\}, K\cup \{q\}, K \cup \{p, q\})$ de parties appartenant à $\Ic$ les applications linéaires $u_{Kp}$, $u_{Kq}$, $u_{Kpq}$, $u_{Kqp}$ et $v_{Kp}$, $v_{Kq}$, $v_{Kpq}$, $v_{Kqp}$ données par le diagramme 
$$
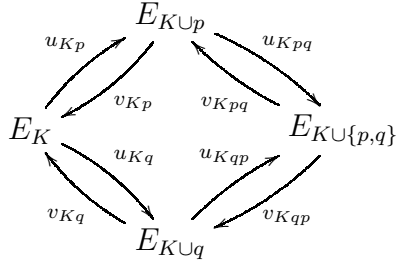
\xymatrix{& E_{K\cup p} \ar@/^/[ld]^{v_{Kp}} \ar@/^/[rd]^{u_{Kpq}}&\\
E_{K}\ar@/^/[ur]^{u_{Kp}} \ar@/^/[dr]^{u_{Kq}}& &E_{K \cup \{p, q\}} \ar@/^/[ul]^{v_{Kpq}} \ar@/^/[dl]^{v_{Kqp}}\\
& E_{K \cup q} \ar@/^/[lu]^{v_{Kq}} \ar@/^/[ru]^{u_{Kqp}}
}$$  
 vérifient les conditions suivantes : 
$$u_{Kp} u_{Kpq}= u_{Kq} u_{Kqp}, ~v_{Kpq}v_{Kp}=v_{Kqp} v_{Kq}, ~ v_{Kpq}u_{Kqp}=u_{Kp} v_{Kq}$$
\item[(iii)] pour tout couple $(K, K')$ maximaux de $\Ic$, considérons l'ensemble des vecteurs $\{e_{i}\}_{i \in \Jc}$ tels que $\{e_{i}\}_{i \in I}$ soit la base fixée de $\ZZ^n$ qui contienne $\{v_{i}\}_{i \in K}$ et $\{e_{i}\}_{i \in I'}$ soit la base de $\ZZ^n$ qui contienne $\{v_{i}\}_{i\in K'}$ et tels que $e_{i}=v_{i}$ pour tout $i \in K\cup K'$; alors pour tout $J \subset K \cap K'$  et $p \in I' \backslash K \cap K'$, l'application linéaire $M_{Jp}$ vérifie  : 
$$M_{Jp}=M_{J, i_{i}}^{\alpha^p_{i_{1}}}\ldots M_{J, i_{m}}^{\alpha^p_{i_{m}}}$$
où $\{i_{1}, \cdots, i_{m}\}\cup J=I$ et $(\alpha^p_{i_{1}}, \cdots, \alpha^p_{i_{n}})$ sont les coordonnées du vecteur $v_{p}$ dans la base $\{e_{i}\}_{i \in I}$.
\end{itemize}
\end{Def}
\noindent
Pour mieux comprendre cette définition donnons quelques exemples : \\
\begin{itemize}
\item[$\bullet$]
Quand la variété torique est $\CC^n$ on retrouve la catégorie $\Cc_{\CC^n}$ définie dans le paragraphe précédent, par exemple la catégorie $\Cc_{\CC^2}$ est la catégorie formée des objets :
$$
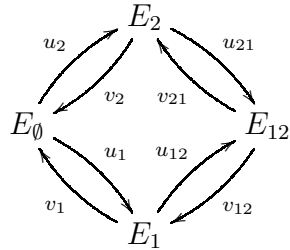
\xymatrix{ & E_{2} \ar@/^/[ld]^{v_{2}} \ar@/^/[rd]^{u_{21}} \\
E_{\emptyset} \ar@/^/[ur]^{u_{2}} \ar@/^/[dr]^{u_{1}}& & E_{12} \ar@/^/[ul]^{v_{21}} \ar@/^/[dl]^{v_{12}}\\
& E_{1} \ar@/^/[lu]^{v_{1}} \ar@/^/[ru]^{u_{12}}
}$$
tels que : 
\begin{itemize}
\item[-] 
$u_{1} u_{12}= u_{2} u_{21}, ~v_{12}v_{1}=v_{21} v_{2}, ~ v_{ji}u_{ij}=u_{j} v_{i}$
\item[-] $M_{i}= v_{i} \circ v_{j} + Id$ et $M_{ij}= v_{ij} \circ u_{ij} +Id$ soient inversibles. 
\end{itemize}
\item[$\bullet$]
Considérons l'éventail $\Delta$ dans $\RR$ formé de trois cônes $\Delta_{\emptyset}=\{0\}$, $\Delta_{1}=\RR_{+}=<e_{1}>$ et $\Delta_{2}=\RR_{-}=<e_{2}>$, la variété torique associée est $\PP^1$.\\
La catégorie $\Cc_{\PP^1}$ est la catégorie formée des objets suivants : 
$$\xymatrix{E_{1}\ar@/^/[r]^{u_{\emptyset1}} & E_{\emptyset} \ar@/^/[l]^{v_{\emptyset1}} \ar@/^/[r]^{v_{\emptyset2}} & E_{2} \ar@/^/[l]^{u_{\emptyset2}}}$$
tels que : 
\begin{itemize}
\item[-] $v_{\emptyset i} \circ u_{\emptyset i} + Id= M_{\emptyset i}$ soit inversible,
\item[-] comme $e_{1}=-e_{2}$ la condition $(iii)$ de la définition \ref{CcDelta} se traduit  par :
$$M_{\emptyset1}=(M_{\emptyset2})^{-1}$$
\end{itemize}
\item[$\bullet$]
La catégorie $\Cc_{\PP^2}$ est la catégorie formée des objets : 
$$
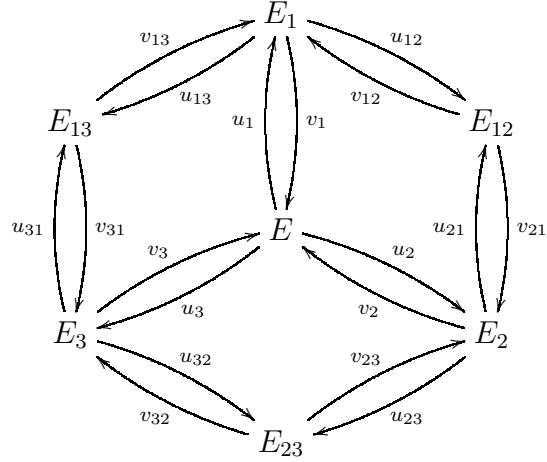
\xymatrix{    & & E_1 \ar@/^/[dll]^{u_{13}}  \ar@/^/[drr]^{u_{12}}  \ar@/^/[dd]^{v_{1}}  \\
                           E_{13}  \ar@/^/[urr]^{v_{13}}  \ar@/^/[dd]^{v_{31}}  & &  & & E_{12}   \ar@/^/[ull]^{v_{12}}  \ar@/^/[dd]^{v_{21}} \\
                           & & E \ar@/^/[uu]^{u_1}  \ar@/^/[dll]^{u_{3}}  \ar@/^/[drr]^{u_{2}}   \\
                           E_3  \ar@/^/[uu]^{u_{31}}  \ar@/^/[urr]^{v_{3}}  \ar@/^/[drr]^{u_{32}}   & & & & E_2  \ar@/^/[uu]^{u_{21}}  \ar@/^/[ull]^{v_{2}}  \ar@/^/[dll]^{u_{23}}   \\
                           & & E_{23}  \ar@/^/[ull]^{v_{32}}  \ar@/^/[urr]^{v_{23}}  \\   }$$  
tels que :
\begin{itemize}
\item[-] 
$u_{i} u_{ij}= u_{j} u_{ji}, ~v_{ij}v_{i}=v_{ji} v_{j}, ~ v_{ji}u_{ij}=u_{j} v_{i}$
\item[-] $M_{\emptyset i}= v_i \circ u_i+ Id$ et $M_{ij}= v_{ij} \circ u_{ij} + Id$ soient inversibles,
\end{itemize}
Regardons maintenant ce que signifie la condition $(iii)$ :
\begin{itemize}
\item[-] considérons par exemple $M_{\emptyset 3}$, $e_{3}$ a pour coordonnée $(-1,-1)$ dans la base $\Bc=\emptyset\cup(e_{1}, e_{2})$, ainsi on a :
$$M_{\emptyset3} =M_{\emptyset1}^{-1} M_{\emptyset2}^{-1},$$
\item[-]  intéressons nous maintenant à l'endomorphisme $M_{13}$, cette fois $\Bc=\{e_{1}\} \cup \{e_{2}\}$ on a ainsi :
$$M_{13}=M_{12}^{-1}$$ 
De même on a pour tout $(i,j,k) \in \{1,2,3\}^3$, tous les trois différents on a :
$$M_{ij}= M_{ik}^{-1}$$
\end{itemize}   
\end{itemize}                                               
Le but de ce paragraphe est de démontrer le théorème suivant : 
\begin{thm}\label{torsection}
Les catégories $\Cc_{\Delta}$ et $\Pc erv_{\Delta}$ sont équivalentes.
\end{thm}
\dem
Construisons tout d'abord un champ sur $X_{\Delta}$ équivalent au champ des faisceaux pervers $\PPP_{\Delta}$. Nous montrons ensuite que la catégorie des sections globales de ce champ est équivalente à la catégorie $\Cc_{\Delta}$.\\

L'ensemble des variétés $\{X_K\}_{K\in \Mm\Ic}$ associées aux cônes maximaux de $\Delta$ forme un recouvrement d'ouverts de $X_{\Delta}$.  Ainsi d'après le lemme \ref{recouvrement} le champ, $\PPP_{\Delta}$, des faisceaux  pervers sur $X_{\Delta}$ stratifié par l'action du tore est équivalent au champ défini par la donnée : 
$$(\{\PPP_K\}_{I \in \Mm\Ic}, \{H_{IJ}\}, \{\Upsilon_{IJK}\})$$
où les $H_{IJ}$ sont les équivalences de champs : 
$$
H_{IJ} : h_{IJ}^{-1}(\PPP_I\mid_{X_{I \cap J}}) \longrightarrow \PPP_{X_{J}}\mid_{X_{I\cap J}}
$$
définies pour un ouvert $U$ de $X_{I}\cap X_{J}$ par : 
$$\begin{array}{cccc}
H_{IJ}(U) : & \Pc erv_{h_{IJ}(U)}  & \longrightarrow & \Pc erv_{U}\\
& \Fc & \longmapsto & h_{IJ}^{-1}(\Fc)
\end{array}$$
et où $\Upsilon_{IJK}$ sont les isomorphismes de foncteurs donnés par la relation :
$$h_{IJ}\mid_{X_{IJK}}\circ h_{JK}\mid_{X_{IJK}}=h_{IK}\mid_{X_{IJK}}$$
où $X_{IJK}= X_{I \cap J \cap K}$.\\
Notons que les variétés $X_{I}$ étant lisses elles sont isomorphes au produit $\CC^k\times (\CC^*)^{n-k}$ où $k$ est le cardinal de $I$. De plus, cet isomosphisme envoie la stratification définie par l'action du tore sur le croisement normal. On considère alors les données : 
$$(\{\CCC_{I}\}_{I \in \Mm \Ic}, \{\Phi_{IJ}\}, \{\theta_{IJK}\})$$
où 
\begin{itemize}
\item
les champs $\CCC_{I}$ sont les champs, sur $X_{I}$, définis au chapitre $4$ , 
\item
les $\Phi_{IJ}$ sont les équivalences de champs : 
$$\Phi_{IJ} : \CCC_{I}\mid_{X_{I \cap J}} \longrightarrow \CCC_{J}\mid_{X_{I \cap J}}$$
définies par la composée des foncteurs :
$$\Phi_{IJ}=\boldsymbol{\alpha}_{J}\mid_{X_{I \cap J}}\circ H_{IJ} \circ \boldsymbol{\beta}_{I} \mid_{X_{I \cap J}}$$
où $\boldsymbol{\alpha}_{I}$ est l'équivalence définie au chapitre $4$ et $\boldsymbol{\beta}_{I}$ est un quasi-inverse fixé de $\boldsymbol{\alpha}_{I}$.
Notons qu'il existe des isomorphismes de foncteurs $a_{IJ}$ : 
$$\shorthandoff{;:!?}
\xymatrix@!0 @C=1,5cm @R=1cm{ 
\PPP_{I}\mid_{X_{I\cap J}} \ar[ddd]_{\boldsymbol{\alpha}_{I}\mid_{X_{I \cap J}}} \ar[rrr]^{H_{IJ}}&& & \PPP_{J}\mid_{X_{I \cap J}} \ar[ddd]^{\boldsymbol{\alpha}_{J}\mid_{X_{I \cap J}}}\\
& &\ar@{=>}[ld]^{a_{IJ}}_{\sim}\\
~&~&\\
\CCC_{I}\mid_{X_{I \cap J}} \ar[rrr]_{\Phi_{IJ}} && & \CCC_{J}\mid_{X_{I \cap J}} }$$
\item
et les $\theta_{IJK}$ sont des isomorphismes de foncteurs : 
$$\theta_{IJK} : \Phi_{IJ}\circ \Phi_{JK} \buildrel\sim\over\longrightarrow \Phi_{IK}$$
tels que le diagramme suivant commute : 
$$\xymatrix{
\Phi_{IJ} \circ \Phi_{JK} \circ \Phi_{KL} \ar[r]^{\Phi_{IJ} (\theta_{JKL})} \ar[d]_{\theta_{IJK}\circ \Phi_{KL}} & \Phi_{IJ} \circ \Phi_{JL} \ar[d]^{\theta_{IJL}}\\
\Phi_{IK} \circ \Phi_{KL} \ar[r]_{\theta_{IKL}} & \Phi_{IL} 
}$$
\end{itemize}
Alors la donnée $\boldsymbol{\alpha}=(\{\boldsymbol{\alpha}_{I}\}_{I \in \Mm \Ic}, \{a_{IJ}\})$ est une équivalence de \linebreak $(\{\PPP_{I}\}_{I \in \Mm \Ic}, \{H_{IJ}\}, \{\Upsilon_{IJK}\})$ dans $(\{\CCC_{I}\}_{I \in \Mm \Ic},\{\Phi_{IJ}\}, \{\theta_{IJK}\})$. \\
Ainsi on a le théorème suivant :
\begin{thm}
Le champ $\PPP_{\Delta}$ est équivalent au champ $\CCC_{\Delta}$ défini par les données : 
$$(\{\CCC_{I}\}_{I \in \Mm \Ic},\{\Phi_{IJ}\}, \{\theta_{IJK}\})$$
\end{thm}
Pour démontrer le théorème \ref{torsection}, il reste à montrer que la catégorie des sections globales de $\CCC_{\Delta}$ est équivalente à la catégorie $\Cc_{\Delta}$. Un objet des sections globales de $\CCC_{\Delta}$ est donné par une famille :
$$\Big\{\{P_{I}\}_{I \in \Mm\Ic}, \{\omega_{IJ}\}_{(I,J) \in \Mm\Ic^2}\Big\}$$
où $P_{I}$ est un objet de $\Gamma(X_{I},\CCC_{I})$ et $\omega_{IJ}$ est un isomorphisme :
$$\omega_{IJ} : \Phi_{IJ}(P_{I}\mid_{X_{I\cap J}})\buildrel\sim\over\longrightarrow P_{J}\mid_{X_{I \cap J}}$$
Ainsi pour expliciter les objets de la catégorie $\Gamma(X_{\Delta},\CCC_{\Delta})$ il faut expliciter les objets des catégories $\Gamma(X_{I},\CCC_{I})$ et la section globale des équivalences $\Phi_{IJ}$.
\begin{lem}\label{rep}
Pour tout $I \in \Ic$, la catégorie $\Gamma(X_{I},\CCC_{I})$ est équivalente à la catégorie $Rep(\pi_{1}(S_{I}), \Cc_{i})$.
\end{lem}
\dem
On peut supposer sans perte de généralité que $X_{I}=\CC^{*n-i}\times \CC^{i}$ et donc que $S_{I}=\CC^{*n-i}\times \{0\}^{i}$. Pour tout i-uplet réel $\varepsilon=(\varepsilon_{1}, \ldots, \varepsilon_{i})$, le foncteur de restriction suivant, où $B_{\varepsilon}$ est le polydisque ouvert centré en zéro et de polyrayon $\varepsilon$ :
$$\Gamma(\CC^{*n-i}\times \CC^{i} , \CCC_{I}) \longrightarrow \Gamma(\CC^{n-i}\times B_{\varepsilon} , \CCC_{I})$$
est une équivalence de catégorie. De même si $\varepsilon'$ est un autre i-uplets tel que $\varepsilon_{i}'<\varepsilon_{i}$ le foncteur de restriction :
$$\Gamma(\CC^{*n-i}\times B_{\varepsilon} , \CCC_{I}) \longrightarrow \Gamma(\CC^{n-i}\times B_{\varepsilon'} , \CCC_{I})$$
est aussi une équivalence de plus elle commute avec la restriction de $\CC^{*n-i}\times \CC^{i}$ dans $\CC^{n-i}\times B_{\varepsilon'}$. Ainsi le foncteur naturel suivant est une équivalence :
$$\Gamma(\CC^{*n-i}\times \CC^{i} , \CCC_{I}) \buildrel\sim\over\longrightarrow 2\varinjlim_{\varepsilon}\Gamma(\CC^{n-i}\times B_{\varepsilon} , \CCC_{I})$$
Or la strate $S_{I}$ est un fermé dans l'ouvert paracompact $\CC^{*n-i}\times \CC^i$ donc, d'après le lemme \ref{paracompact}, le le foncteur naturel suivant est une équivalence de catégorie, 
$$2\varinjlim_{\varepsilon}\Gamma(\CC^{*n-i}\times B_{\varepsilon} , \CCC_{I}) \longrightarrow \Gamma(S_{I}, \CCC_{I})$$
de plus la restriction  de $\CCC_{I}$ à la strate $S_{I}$ est constant donc on a bien l'équivalence annoncée.
\cqfd
\'Etudions maintenant les sections globales des équivalences $\Phi_{K,K'}$ où $(K,K')$ est un couple de parties maximales de $\Ic$. On considère l'ensemble des vecteurs $\{e_{i}\}_{i \in \Jc}$ tels que $\{e_{i}\}_{i \in I}$ soit la base fixée de $\ZZ^n$ qui contienne $\{v_{i}\}_{i \in K}$ et $\{e_{i}\}_{i \in I'}$ soit la base de $\ZZ^n$ qui contienne $\{v_{i}\}_{i\in K'}$ et tels que $e_{i}=v_{i}$ pour tout $i \in K\cup K'$. \\
D'après le lemme \ref{rep} un objet de $\Gamma(X_{K \cap K'}, \CCC_{K \cap K'})$ est donné par :
\begin{itemize}
\item[$\bullet$] pour tout $J \subset K \cap K'$, un espace vectoriel $E_{J}$,
\item[$\bullet$] pour tout couple $(J, J \cup p) $ de sous-ensemble de $K \cap K'$, deux applications linéaires $u_{Jp}$ et $v_{Jp}$ :
$$\begin{array}{cccc}
u_{Jp} :& E_{J}& \longrightarrow& E_{J\cup p}\\
v_{Jp} :& E_{J\cup p}& \longrightarrow& E_{J}\\
\end{array}$$
\item[$\bullet$] pour tout $J \subset K\cup K'$ et tout $q \in I\backslash K \cap K'$ un endomorphisme $M_{Jq}$.
\end{itemize}
\begin{Def}\label{Lambda}
On note $\Upsilon$ le foncteur qui a un objet $$\begin{displaystyle}\big\{ \{E_{J}\}_{J \subset K \cap K'}, \{u_{Jp}, v_{Jp}\}_{\substack{J\subset K\cap K'\\ p\in K\cap K'\backslash J}}, \{M_{Jq}\}_{\substack{J \subset K\cap K'\\ q \in I\backslash K \cap K'}}\big\}\end{displaystyle}$$ de $\Gamma(X_{K \cap K'}, \CCC_{K \cap K'})$ associe l'objet :
$$\begin{displaystyle}\big\{ \{E_{J}\}_{J \subset K \cap K'}, \{u_{Jp}, v_{Jp}\}, \{M_{Jq'}\}_{\substack{J \subset K\cap K'\\ q' \in I'\backslash K \cap K'}}\big\} \end{displaystyle}$$ 
où si $J \cup \{i_{1}, \ldots, i_{m}\}=I$ et $(\alpha_{1}^{(q')},\ldots, \alpha_{n}^{(q')})$ sont les coordonnées de $e_{q'}$ dans la base $\{e_{i}\}_{i \in I}$ :
$$M_{Jq'}=M_{J,i_{1}}^{\alpha^{(q')}_{i_{1}}}\cdots M_{J,i_{m}}^{\alpha^{(q')}_{i_{m }}}.$$ 
\end{Def}

\begin{prop}
Le foncteur $\Upsilon$ et les sections globales de $\Phi_{K,K'}$ sont isomorphes.
\end{prop}
\dem
On note $l$ le cardinal de $K\cap K'$. Rappelons la définition de $\Phi_{KK'}$ et de $\boldsymbol{\alpha}_{X_{K\cap K'}}$ : 
$$\Phi_{KK'}=\boldsymbol{\alpha}_{K}\mid_{X_{K \cap K'}}\circ H_{KK'} \circ \boldsymbol{\beta}_{K'} \mid_{X_{K \cap K'}}$$
où $H_{KK'}$ est l'équivalence :
$$\begin{array}{cccc}
H_{KK'} : & \Pc erv_{h_{KK'}(U)} & \longrightarrow & \Pc erv_{U}\\
  & \Fc & \longrightarrow & h^{-1}_{KK'}(\Fc)
\end{array}$$
Le foncteur $\boldsymbol{\alpha}_{X_{K\cap K'}}$ est la composée du foncteur $\mu$ :
$$\xymatrix{
\mu :\Pc erv_{X_{K\cap K'}} \ar[r]& Rep(\pi_{1}(S_{K\cap K'}), \Pc erv_{\CC^l}) 
}$$
avec le foncteur $Rep(\pi_{1}(S_{K\cap K'}), \alpha_{l})$ :
$$ Rep(\pi_{1}(S_{K\cap K'}), \Pc erv_{\CC^l}) \longrightarrow Rep(\pi_{1}(S_{K\cap K'}), \Cc_{l})$$
Ainsi on a le diagramme suivant :
$$\xymatrix{
\Pc erv_{X_{K\cap K'}} \ar[r]^{H_{KK'}} \ar[d]_{\mu} & \Pc erv_{X_{K\cap K'}} \ar[d]^{\mu}\\
Rep(\pi_{1}(S_{K\cap K'}), \Pc erv_{\CC^l})  \ar[d]_{Rep(\pi_{1}(S_{K\cap K'}), \alpha_{l})} & Rep(\pi_{1}(S_{K\cap K'}), \Pc erv_{\CC^l}) \ar[d]^{Rep(\pi_{1}(S_{K\cap K'}), \alpha_{l})} \\
Rep(\pi_{1}(S_{K\cap K'}), \Cc_{l}) \ar[r]_{\Upsilon} & Rep(\pi_{1}(S_{K\cap K'}), \Cc_{l})
}$$
Pour démontrer la proposition il suffit de montrer l'existence d'un isomorphisme :
$$\Upsilon \circ Rep(\pi_{1}(S_{K\cap K'}), \alpha_{l}) \circ \mu \simeq Rep(\pi_{1}(S_{K\cap K'}), \alpha_{l}) \circ \mu \circ H_{KK'}$$
Commen\c cons par rappeler la définition de $Rep(\pi_{1}(S_{K\cap K'}), \alpha_{l}) \circ \mu$.
On peut supposer sans perte de généralité que $K\cap K'=\{1,\cdots,l\}$. Dans ce cas $X_{K\cap K'}= \CC^l\times \CC^{*n-l}$.
Soit $\Fc$ un faisceau pervers sur $\CC^l\times \CC^{*n-l}$, l'image par $Rep(\pi_{1}(S_{K\cap K'}), \alpha_{l}) \circ \mu$ est donné par :
$$\Big\{ \big\{(R^j\Gamma_{\RR_{-}^J\CC\backslash\RR_{-}^{K\cap K'\backslash J}\{1\}}\Fc|_{\CC^l\times\{1\}^{n-l}})_{p_{J}}\big\}_{J\subset K\cap K'}, \big\{u_{Jp}, v_{Jp}\big\}, \{M_{Jp}\big\}\Big\} $$
où $u_{Jp}$ et $v_{Jp}$ sont respectivement les fibres en $p_{J}$ des morphismes naturels : 
$$R^j\Gamma_{(\RR_{-}^J\CC\backslash\RR_{-}^{K\cap K'\backslash J} \times\{1\}}(\Fc\mid_{\CC^L}) \longrightarrow R^{j+1}\Gamma_{(\RR_{-}^{J\cup p}\CC\backslash\RR_-^{K\cap K' \backslash J}\{1\}}(\Fc\mid_{\CC^L})$$
$$R^{j+1}\Gamma_{(\RR_{-}^{J\cup p}\CC\backslash\RR_{-})\cap \CC^L}(\Fc\mid_{\CC^L}) \longrightarrow R^{j+1}\Gamma_{(\RR_{-}^J\RR_{-}^{*\{p\}}\CC\backslash\RR_{-})\cap \CC^L}(\Fc\mid_{\CC^L})$$
et d'après le lemme \ref{commono}, les automorphisme $M_{Jp}$ sont les monodromies des faisceaux localement constants $(R^j\Gamma_{\RR^J\CC^{*}}\Fc)|_{S_{J}}$ définies par le chemin $\gamma_{Jp}$.\\
Considérons maintenant la composée :
$$Rep(\pi_{1}(S_{K\cap K'}), \alpha_{l}) \circ \mu \circ H_{KK'}$$
Vu la définition de $H_{K'K}$, il faut expliciter les isomorphismes de recollements $h_{KK'}$.
\begin{lem}
Si l'on considère les mêmes notations que dans la définition \ref{Lambda}, l'isomorphisme $h_{KK'}$ est donné par :
$$\begin{array}{ccccccc}
X_K  \supset& X_{K\cap K'} &\longrightarrow &  X_{K\cap K'} & \subset X_{K'}\\
 & (x_{1}, \ldots, x_{n}) & \longmapsto & (x_{1}^{\alpha_{1}^{(1)}}\cdots x_{n}^{\alpha_{1}^{(n)}}, \ldots,x_{1}^{\alpha_{n}^{(1)}}\cdots x_{n}^{\alpha_{n}^{(n)}} )
 \end{array}$$
\end{lem} 
\dem
Pour démontrer cette proposition nous avons besoin de ce lemme : 
\begin{lem}\label{entrant}
Soit $\sigma$ un cône entier régulier de $\RR^n$ engendré par des vecteurs primitifs $\{v_{1}, \ldots, v_{n}\}$, tels que le déterminant de $M$, la matrice de colonnes $v_{i}$, soit égal à 1.\\
Soit $G$ la matrice $(^tM)^{-1}$ et $\{g_{1}, \ldots, g_{n}\}$ ses vecteurs colonnes, alors le cône $\check{\sigma}$ est engendré par les vecteurs $\{g_{1}, \ldots, g_{n}\}$.
\end{lem}
\dem
D'après la proposition \ref{cone}, pour démontrer ce lemme il suffit de montrer que les vecteurs $g_{i}$ sont les vecteurs normaux entrants aux faces engendrées par $\{v_{1}, \cdots, v_{i-1}, v_{i+1}, \ldots, v_{n}\}$. C'est à dire que :
\begin{itemize}
\item[(i)] $<g_{i},v_{j}>=0$ pour tout $j\ne i$,
\item[(ii)] $\det M_{i} >0$, où $M_{i}$ est la matrice $(v_{1} \cdots v_{i-1}~ g_{i} ~ v_{i+1} \cdots v_{n})$.
\end{itemize}
La première assertion est donnée par l'égalité $^{t}M G=I_{n}$. \\
Cette égalité donne de plus $<g_{i},v_{i}>=1 $, ainsi la matrice $^tGM_{i}$ est égale à :
$$
  \left(
     \raisebox{0.5\depth}{%
       \xymatrixcolsep{1ex}%
       \xymatrixrowsep{1ex}%
       \xymatrix{
         1 \ar @{-}[dddrrr] & & & <g_{1}, g_{i}> \ar@{-}[ddd]& & & \\
         \\
         \\
         &&& \Vert g_{i} \Vert \ar@{-}[dddrrr] \ar@{-}[ddd]&&& \\
          \\
          \\
          && &<g_{n}, g_{i}> & & & 1
       }%
     }
   \right)
  $$
  donc son determinant est égal à $\Vert g_{i}\Vert$, or $\det ^tG M_{i}= \det M_{i}$, ce qui démontre la deuxième assertion.
\cqfd
On peut supposer sans perte de généralité que $e_{1}, \ldots,e_{n}$ est la base canonique. Les vecteurs entrants normaux aux facettes de $\sigma_{I}$ sont alors encore la famille $e_{1}, \ldots, e_{n}$. On note $g_{i'_{1}}, \ldots, g_{i'_{n}}$ les vecteurs entrants normaux aux facettes de $\Delta_{I'}$. D'après la proposition \ref{cone} et comme les variétés $X_{I}$ et $X_{I'}$ sont lisses, l'isomorphisme de recollement provient du passage de la base $g_{i'_{1}}, \ldots, g_{i'_{n}}$ de $\ZZ^n$ à la base $e_{1}, \ldots, e_{n}$. La matrice de changement de base est l'inverse de la matrice dont les colonnes sont les vecteurs $g_{i'_{n}}$ mais d'après le lemme \ref{entrant}, c'est la matrice dont les lignes sont les coordonnées des $\varepsilon_{i}$. D'après la définition d'un morphisme torique on retrouve bien l'isomorphisme cherché.
\cqfd
Toujours en supposant que $K\cap K'=\{1,\cdots,l\}$ on a, pour tout $i<l$,  $e_{i}=\varepsilon_{i}$  et donc :
$$\alpha_{j}^{(i)}=\delta_{ij}$$
Ainsi l'application $h_{KK'}$ est la suivante :
$$\begin{array}{cccc}
h_{KK'} : &\CC^l\times \CC^{*n-l} & \longrightarrow & \CC^l\times \CC^{*n-l}\\
 & (x_{1}, \ldots, x_{n}) & \longmapsto & (y_{1}, \ldots, y_{n})
\end{array}$$
où
\begin{itemize}
\item pour $i\leq l$, $y_{i}=x_{i}x_{l+1}^{\alpha_{i}^{(l+1)}}\cdots x_{n}^{\alpha_{i}^{(l+1)}}$ 
\item pour $i>l$, $y_{i}=x_{l+1}^{\alpha_{i}^{(l+1)}}\cdots x_{n}^{\alpha_{i}^{(l+1)}}$. 
\end{itemize}
Revenons à la composée $Rep(\pi_{1}(S_{K\cap K'}), \alpha_{l}) \circ \mu \circ H_{KK'}$.  Considérons l'image de $\Fc$ par $\mu \circ H_{KK'}$. 
$$\mu \circ H_{KK'}(\Fc)= \big( (h^{-1}_{KK'}\Fc)|_{\CC^l\times\{1\}^{n-l}}, \{\widetilde{\gamma_{p}\circ h_{KK'}}\}_{p\in I'\backslash K\cap K'}\big)$$
où $\widetilde{\gamma_{p}\circ h_{KK'}}$ est défini comme l'était $\tilde{\gamma}_{p}$.  Mais la restriction de $h_{KK'}$ à $\CC^l\times \{1\}^{n-l}$ est l'identité donc :
$$ (h^{-1}_{KK'}\Fc)|_{\CC^l\times\{1\}^{n-l}}\simeq\Fc|_{\CC^l\times\{1\}^{n-l}}$$
Appliquons maintenant le foncteur $Rep(\pi_{1}(S_{K\cap K'}), \alpha_{l})$ à cet objet. On obtient une famille isomorphe à la famille suivante :
$$\Big\{ \big\{(R^j\Gamma_{\RR_{-}^J\CC\backslash\RR_{-}^{K\cap K'\backslash J}\{1\}}\Fc|_{\CC^l\times\{1\}^{n-l}})_{p_{J}}\big\}_{J\subset K\cap K'}, \big\{u_{Jp}, v_{Jp}\big\}, \{M'_{Jp}\big\}\Big\} $$
où $u_{Jp}$ et $v_{Jp}$ sont respectivement les fibres en $p_{J}$ des morphismes naturels : 
$$R^j\Gamma_{(\RR_{-}^J\CC\backslash\RR_{-}^{K\cap K'\backslash J} \times\{1\}}(\Fc\mid_{\CC^L}) \longrightarrow R^{j+1}\Gamma_{(\RR_{-}^{J\cup p}\CC\backslash\RR_-^{K\cap K' \backslash J}\{1\}}(\Fc\mid_{\CC^L})$$
$$R^{j+1}\Gamma_{(\RR_{-}^{J\cup p}\CC\backslash\RR_{-})\cap \CC^L}(\Fc\mid_{\CC^L}) \longrightarrow R^{j+1}\Gamma_{(\RR_{-}^J\RR_{-}^{*\{p\}}\CC\backslash\RR_{-})\cap \CC^L}(\Fc\mid_{\CC^L})$$
et d'après le lemme \ref{commono}, les automorphismes $M'_{Jp}$ sont les monodromies des faisceaux localement constants $(R^j\Gamma_{\RR^J\CC^{*}}\Fc)|_{S_{J}}$ définies par le chemin $\gamma_{Jp}\circ h_{KK'}|_{S_{J}}$.\\
On peut supposer sans perte de généralité que $J=\{1, \ldots, j\}$, avec $j<l$. La composition de ces applications est alors donnée par :
$$\begin{array}{cccc}
\gamma_{Jp}\circ h_{KK'} : &[0,1] & \longrightarrow & S_{J}\\
& t & \longmapsto & (\underbrace{0, \ldots,0}_{j}, e^{\alpha_{j+1}^{(p)}2i\pi t}, \ldots, e^{\alpha_{n}^{(p)}2i\pi t})
\end{array}$$
On a donc bien $M'_{Jp}=M_{J,j+1}^{\alpha_{j+1}^{(p)}} \cdots M_{J,n}^{\alpha_{n}^{(p)}}$.\\ 
Ce qui montre que :
$$\Upsilon \circ Rep(\pi_{1}(S_{K\cap K'}), \alpha_{l}) \circ \mu \simeq Rep(\pi_{1}(S_{K\cap K'}), \alpha_{l}) \circ \mu \circ H_{KK'}$$
\cqfd
Un objet de $\Gamma(X_{\Delta}, \CCC_{\Delta})$ est donc donné par : 
\begin{itemize}
\item pour tout ensemble, $K$, maximal de $\Ic$, une famille :
$$\Big\{\{E^{K}_{J}\}_{J\subset I}, \{u^{K}_{Jp}, v^{K}_{Jp}\}_{\substack{J\subset I\\ p \in I\backslash J}}, \{M^{K}_{Jq}\}_{\substack{J\subset I\\ j+1\leq q\leq n }}\Big\}$$
\item pour toute intersection, $K\cap K'$, d'ensemble maximal de $\Ic$ et pour tout ensemble $J\subset K\cap K'$ un isomorphisme :
$$ \delta^{KK'}_{J} : E_{J}^{K} \longrightarrow E_{J}^{K'}$$
tels que $$\delta_{J\cup p}^{KK'-1} u^{K'}_{Jp}\delta_{J}^{KK'}= u^{K}_{Jp} \text{,~}\delta_{J}^{KK'-1} v^{K'}_{Jp}\delta_{J\cup p}^{KK'}= v^{K}_{Jp}$$
$$\text{et}$$
$$\delta_{J}^{K'K-1}M'_{Jp}\delta_{J}^{K'K}=M_{J,i_{1}}^{\alpha^{(q')}_{i_{1}}}\cdots M_{J,i_{m}}^{\alpha^{(q')}_{i_{m }}}$$ 
où $J \cup \{i_{1}, \ldots, i_{m}\}=I$ et $(\alpha_{1}^{(q')},\ldots, \alpha_{n}^{(q')})$ sont les coordonnées de $e_{q'}$ dans la base $\{e_{i}\}_{i \in I}$
\end{itemize}
Comme dans le chapitre précédent nous définissons une équivalence, $\Lambda$, entre la catégorie $\Gamma(X_{\Delta}, \CCC_{\Delta})$ et la catégorie $\Cc_{\Delta}$. 
$$\Lambda : \Gamma(X_{\Delta}, \CCC_{\Delta}) \longrightarrow \Cc_{\Delta}$$
Pour cela on note $K_{1}, \ldots, K_{m}$, les parties maximales de $\Ic$.
On définit alors $\Lambda$ comme le foncteur qui à un objet décrit ci-dessus associe la famille donnée par 
\begin{itemize}
\item pour toute partie maximale, $K_{i}$, de $\Ic$, pour toute partie $J$ de $K_{i}\backslash \bigcup_{j<i}K_{j}\cup K_{i}$, l'espace vectoriel $E^{K_{i}}_{J}$ et les applications linéaires $M^{K_{i}}_{Jq}$,  pour tout couple $(J, J\cup p)$ de parties de $K_{i}\backslash \bigcup_{j<i}K_{j}\cup K_{i}$, les applications linéaires $u_{Jp}$ et $v_{Jp}$,
\item pour tout couple $(J, J\cup p)$, tel que $J \cup p\subset K_{i} \backslash \bigcup_{j<i}K_{j}\cup K_{i}$ et $J\subset K_{i}\cap K_{k}$ avec $k<j$, les applications linéaires :
$$u^{K_{i}}_{Jp}\circ \delta_{J}^{K_{k}K_{j}} \text{~~~et~~~} \delta_{J}^{K_{k}K_{j}-1}\circ v_{Jp}^{K_{i}}$$
\end{itemize}
Cette famille est bien un objet de $\Cc_{\Delta}$, les relations de commutation que vérifient les isomorphisme $\delta_{J}^{KK'}$  assurent notamment la condition $(iii)$. Comme au chapitre précédent le foncteur qui à un objet
$$\Big\{ \{E_{J}\}_{J \in \Ic}, \{u_{Jp}, v_{Jp}\}_{\substack{(J,J\cup p) \in \Ic^2}} \{M_{Ji}\}_{1\leq i\leq j}\Big\}$$
 de $\Cc_{\Delta}$ associe l'objet de $\Gamma(X_{\Delta}, \CCC_{\Delta})$ donné par 
 \begin{itemize}
 \item pour tout ensemble maximal $K$ de $\Ic$, la famille :
 $$\Big\{ \{E_{J}\}_{J \subset K}, \{u_{Jp}, v_{Jp}\}_{\substack{(J,J\cup p) \subset K^2}} \{M_{Ji}\}_{\substack{J\subset K\\1\leq i\leq j}}\Big\}$$
 \item pour toute intersection d'ensemble maximal $K\cap K'$ de $\Ic$ et pour tout ensemble $J\subset K\cap K'$ l'identité.
 \end{itemize}
Les conditions de la définition \ref{CcDelta} assurent que cette donnée est bien un objet de $\Gamma(X_{\Delta}, \CCC_{\Delta})$.\\
Ces deux foncteurs sont quasi-inverse l'un de l'autre.
\cqfd

\bibliographystyle{plain}
\bibliography{biblio}
\newpage
~\\
\thispagestyle{empty}
\newpage
\thispagestyle{empty}

\thispagestyle{empty}

\small \subsection*{R\'esum\'e}
Le but de cette thèse est de montrer comment des descriptions élémentaires et explicites de la catégorie des faisceaux pervers peuvent être recollées en une description de la catégorie globale des faisceaux pervers. 
Pour cela nous généralisons une équivalence entre la catégorie des faisceaux pervers sur $\CC^n$ stratifié par le croisement normal et une catégorie de représentation d'un certain carquois démontrée par A. Galligo, M. Granger, Ph. Maisonobe (GGM) en une équivalence de champs. \\
Dans un premier temps nous nous intéressons à la fa\c con de définir un champ sur un espace stratifié. Ainsi nous démontrons la $2$-équivalence entre la $2$-catégorie des champs et une $2$-catégorie dont les objets sont les données de champs sur chacune des strates et de conditions de recollement. Ceci nous permet, dans le cas où l'espace considéré est $\CC^n$ stratifié par le croisement normal de définir assez simplement un champ de carquois strictement constructible. \\
Dans un deuxième temps nous montrons que ce champ est bien équivalent au champ des faisceaux pervers. La difficulté majeure consiste alors à démontrer que les équivalences de catégories de GGM se recollent correctement en une équivalence de champs.\\
Enfin les deux derniers chapitres sont des applications de cette technique. On obtient ainsi de nouvelles descriptions en terme de re\-pré\-sen\-ta\-tion de carquois de la catégorie des faisceaux pervers d'une part sur $\CC^2$ stratifié par un arrangement générique de droites et d'autre part sur les variétés toriques lisses.

\end{document}